\documentclass[12pt,twoside]{article}   	
\usepackage[margin=1in]{geometry}                		
\usepackage[colorinlistoftodos]{todonotes}
\usepackage[colorlinks=true, allcolors=blue]{hyperref}

\usepackage{fancyhdr}
\usepackage{graphicx}
\usepackage{subfigure}
\usepackage{amssymb}
\usepackage{amsmath}
\usepackage{amsfonts}
\usepackage{latexsym}
\usepackage{theorem}
\usepackage{mathrsfs}
\usepackage{mathtools}
\usepackage{empheq}
\usepackage{bm}
\usepackage{bbm}
\usepackage{cases}
\usepackage{xcolor}
\usepackage{setspace}
\usepackage{epstopdf}
\usepackage{float}
\usepackage{tikz}
\usetikzlibrary{calc,intersections}
\usepackage{cite}
\usepackage{nicefrac}
\usepackage{url}
\usepackage{booktabs} 
\usepackage{array} 
\usepackage{paralist} 
\usepackage{verbatim} 
\usepackage{subfigure} 
\usepackage{comment}

\theoremstyle{plain}			
\newtheorem{thm}{Theorem}[section]

\newtheorem{rmk}[thm]{Remark}

\newcommand*\xbar[1]{%
  \hbox{%
    \vbox{%
      \hrule height 0.5pt 
      \kern0.4ex
      \hbox{%
        \kern-0.05em
        \ensuremath{#1}%
        \kern-0.00em
      }%
    }%
  }%
}

\newcommand\eref[1]{(\ref{#1})}

\allowdisplaybreaks[1]

\numberwithin{equation}{section}
\numberwithin{figure}{section}
\numberwithin{table}{section}

\newcommand{\dt}{\Delta t}
\newcommand{\dx}{\Delta x}
\newcommand{\dy}{\Delta y}
\newcommand{\hf}{{\frac{1}{2}}}
\newcommand{\jph}{{j+\frac{1}{2}}}
\newcommand{\jmh}{{j-\frac{1}{2}}}
\newcommand{\kph}{{k+\frac{1}{2}}}
\newcommand{\kmh}{{k-\frac{1}{2}}}
\newcommand{\Sc}{\mathrm{Sc}} 

\graphicspath{{Figures/}}

\title{A Diffuse-Domain Based Numerical Method for a Chemotaxis-Fluid Model}
\author{Chenxi Wang\thanks{Department of Mathematics, Harbin Institute of Technology, Harbin, 150001, China and Department of Mathematics,
Southern University of Science and Technology (SUSTech), Shenzhen, 518055, China; \tt{wangcx2017@mail.sustech.edu.cn}}, Alina
Chertock\thanks{Department of Mathematics, North Carolina State University, Raleigh, NC 27695, USA; {\tt chertock@math.ncsu.edu}},~ Shumo
Cui\thanks{Department of Mathematics and International Center for Mathematics, Southern University of Science and Technology (SUSTech),
Shenzhen, 518055, China; {\tt cuism@sustech.edu.cn}}, ~ Alexander Kurganov\thanks{Department of Mathematics, International Center for
Mathematics and Guangdong Provincial Key Laboratory of Computational Science and Material Design, Southern University of Science and
Technology (SUSTech), Shenzhen 518055, China; {\tt alexander@sustech.edu.cn}},\\ and Zhen Zhang\thanks{Department of Mathematics,
Guangdong Provincial Key Laboratory of Computational Science and Material Design and National Center for Applied Mathematics (Shenzhen),
Southern University of Science and Technology (SUSTech), Shenzhen 518055, China; \tt{zhangz@sustech.edu.cn}}}

\begin{document}
\date{}
\maketitle

\begin{abstract}
In this paper, we consider a coupled chemotaxis-fluid system that models self-organized collective behavior of oxytactic bacteria in a
sessile drop. This model describes the biological chemotaxis phenomenon in the fluid environment and couples a convective chemotaxis system
for the oxygen-consuming and oxytactic bacteria with the incompressible Navier–Stokes equations subject to a gravitational force, which is
proportional to the relative surplus of the cell density compared to the water density.

We develop a new positivity preserving and high-resolution method for the studied chemotaxis-fluid system. Our method is based on the
diffuse-domain approach, which we use to derive a new chemotaxis-fluid diffuse-domain (cf-DD) model for simulating bioconvection in complex
geometries. The drop domain is imbedded into a larger rectangular domain, and the original boundary is replaced by a diffuse interface with
finite thickness. The original chemotaxis-fluid system is reformulated on the larger domain with additional source terms that approximate
the boundary conditions on the physical interface. We show that the cf-DD model converges to the chemotaxis-fluid model asymptotically as
the width of the diffuse interface shrinks to zero. We numerically solve the resulting cf-DD system by a second-order hybrid
finite-volume finite-difference method and demonstrate the performance of the proposed approach on a number of numerical experiments that
showcase several interesting chemotactic phenomena in sessile drops of different shapes, where the bacterial patterns depend on the droplet
geometries.
\end{abstract}



\section{Introduction}\label{sec1}
In this paper, we study the following coupled chemotaxis-fluid system in a sessile drop \cite{tuval2005bacterial}:
\begin{equation}
\begin{aligned}
&n_t+\bm{u}\!\cdot\!\nabla n+\chi\nabla\!\cdot\![nr(c)\nabla c]=D_n\Delta n,\\
&c_t+\bm{u}\!\cdot\!\nabla c=D_c\Delta c-n\kappa r(c),\\
&\rho\left(\bm{u}_t+\bm{u}\!\cdot\!\nabla\bm{u}\right)+\nabla p=\eta\Delta\bm{u}-n\nabla\Phi,\\
&\nabla\!\cdot\!\bm{u}=0,
\end{aligned}
\label{1.1}
\end{equation}
where $n$ and $c$ are the concentrations of bacteria and oxygen, respectively, $\kappa$ is the oxygen consumption rate, and
$\bm{u}=(u,v)^\top$ is the velocity field of a fluid flow governed by the incompressible Navier-Stokes equations with density $\rho$,
pressure $p$ and viscosity $\eta$. In the fluid equation, $\nabla\Phi:=V_bg(\rho_b-\rho)\bm{z}$ describes the gravitational force exerted by
a bacterium onto the fluid along the upwards unit vector $\bm{z}$ proportional to the volume of the bacterium $V_b$, the gravitation
acceleration $g=9.8\,m/s^{-2}$, and the density of bacteria is $\rho_b$ (bacteria are about 10$\%$ denser than water).

In \eref{1.1}, both the bacteria and oxygen are convected by the fluid and diffuse with their respective diffusion coefficients $D_n$ and
$D_c$. The bacteria are active as long as a sufficient oxygen supply is available: this is measured by a dimensionless cut-off function
$r(c)$, which can be modeled, for instance, by
\begin{equation}
r(c)=\begin{cases}
1,&c\ge c^*,\\
0,&c<c^*,
\end{cases}
\label{rc}
\end{equation}
where $c^*$ is an inactivity threshold. The active bacteria both consume the oxygen and, in a chemotactic response, are directed towards a
higher oxygen concentration with the rate proportional to the chemotactic sensitivity $\chi$.

\begin{figure}[ht!]
\centerline{\includegraphics[trim=2.3cm 1.2cm 1.8cm 1.9cm, clip,width=8cm]{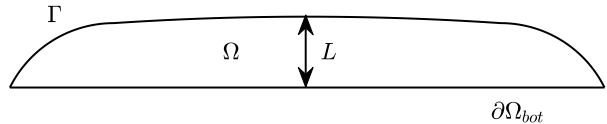}}
\caption{\sf A sketch of the sessile drop domain $\Omega$.\label{OD}}
\end{figure}

A typical shape of the sessile drop is depicted in Figure \ref{OD}. We stress that the boundary conditions on $n$, $c$ and $\bm{u}$ are
essential to ensure that the solutions of \eref{1.1} match well the experimental observations. We shall consider the following mixed
boundary conditions: the boundary conditions on the top interface $\Gamma$ describe the fluid-air surface, which is stress-free, allows no
cell flux, and has saturated air oxygen concentration $c_{air}$:
\begin{equation}
\bm\nu\!\cdot\!\bm u=0,\quad\bm\nu\!\cdot\!\nabla(\bm u\!\cdot\!\bm\tau)=0,\quad
\left(\chi n\nabla c-D_n\nabla n\right)\!\cdot\!\bm\nu=0,\quad c=c_{air},\quad\forall(x,y)\in\Gamma,
\label{1.2}
\end{equation}
where $\bm{\nu}$ and $\bm{\tau}$ are the unit outward normal and tangential vectors on $\Gamma$. A no-slip boundary condition is
applied on the bottom surface ($\partial\Omega_{bot}$) and there is no flux of cells or oxygen through $\partial\Omega_{bot}$:
\begin{equation}
\bm u=\bm0,\quad\nabla n\!\cdot\!\bm\nu=0,\quad\nabla c\!\cdot\!\bm\nu=0,\quad\forall(x,y)\in\partial\Omega_{bot}.
\label{1.3}
\end{equation}

In \cite{chertock2012sinking}, the system \eref{1.1} was studied in a simplified, rectangular shaped domain subject to the same top and
bottom boundary conditions as in \eref{1.2} and \eref{1.3}, respectively, and periodic boundary conditions in the horizontal direction.
Several phenomena of sinking, merging and stationary plumes were discovered in \cite{chertock2012sinking} by numerically solving \eref{1.1}
using a high-resolution hybrid finite-volume finite-difference method. In \cite{deleuze2016numerical}, an upwind finite-element method was
developed and used to investigate the pattern formation and the hydrodynamical stability of the system \eref{1.1} for the same simplified
setup. In \cite{huang2021fully}, a fully decoupled, linear and positivity preserving finite-element method for solving the
chemotaxis-Stokes equations has been recently developed for a similar setup. In \cite{duarte2021numerical}, the chemotaxis-fluid model
without the discontinuous oxygen cut-off function $r(c)$ has been considered, for which a finite-element method has been constructed,
optimal error estimates have been established, and convergence towards regular solutions has been proved. In
\cite{ivanvcic2019free,ivanvcic2020bacterial}, a generalized chemotaxis-diffusion-convection model, which includes the dynamic free surface
and appropriate boundary conditions, has been proposed together with a numerical method, which uses a time dependent grid and incorporates
surface tension and a dynamic contact line.

Stability analysis and dynamics of the chemotaxis-fluid system \eref{1.1} with a deformed free-surface in a shallow chamber were studied
in \cite{chakraborty2018stability}. In particular, a detailed linear stability analysis of a steady-state cell and oxygen concentration
distribution was performed. The chemotaxis-fluid system \eref{1.1}, but without the discontinuous oxygen cut-off function $r(c)$, has been
recently studied in \cite{braukhoff2020global}, where it has been proved that in one or two space dimensions, the system has a unique global
classical solution. In the three-dimensional case, the existence of a global weak solution in a drop shaped domain has been shown and a
uniform in time energy bound has been established.

The main goal of this paper is to develop a robust and accurate numerical method for the chemotaxis-fluid system in the sessile drop domain.
To this end, we extend a diffuse-domain approach to the system \eref{1.1} and construct a new chemotaxis-fluid diffuse-domain (cf-DD) model,
which we numerically solve using a second-order hybrid finite-volume finite-difference method.

The diffuse-domain method was proposed in \cite{li2009solving} following the idea of the smoothed boundary method previously introduced in
\cite{bueno2006spectral2,bueno2006spectral} as a powerful numerical tool for solving diffusion equations with no-flux boundary conditions
imposed at irregular boundaries within the computational domain. The diffuse-domain method can be applied to a variety of PDEs in both
stationary and moving complex geometries with Dirichlet, Neumann or Robin boundary conditions. The key idea of the method is to place the
complex geometry into a larger rectangular domain, introduce a smoothed characteristic function of the original domain, and reformulate the
original PDE(s) on the extended domain with the help of additional source terms, which reflect the contribution of the original boundary
conditions. It has been shown in \cite{guo2021diffuse,li2009solving,lervaag2015analysis,teigen2009diffuse,yu2012extended,yu2020higher} that
the reformulated diffuse-domain model asymptotically converges to the original PDE(s) as the thickness of the diffuse-domain interface tends
to zero. The main advantage of the diffuse-domain method is that the reformulated model can be solved using standard numerical methods even
for very complex domains (with moving boundaries). For example, the diffuse-domain method has been successfully applied to several quite
sophisticated two-phase flow models; see \cite{aland2010two,teigen2011diffuse,aland2011continuum,xu2018level}.

We first follow the diffuse-domain approach and derive a cf-DD model, for which we perform an asymptotic analysis and show that it converges
to the original chemotaxis-fluid model as the thickness of the diffuse-domain interface shrinks to zero. We then use the proposed cf-DD
model to simulate bioconvection in complex droplet geometries using a numerical method, which is derived as follows. The modified cell
density equation is numerically solved by a semi-discrete second-order finite-volume upwind method (introduced in
\cite{chertock2012sinking}) combined with a second-order strong stability-preserving multistep ODE solver, which can be found in, e.g.,
\cite{gottlieb2001strong}. The resulting fully discrete scheme is shown to preserve the positivity of cell density. The modified
Navier-Stokes and oxygen concentration equations are discretized using a second-order projection finite-difference method, combined with the
second-order BDF-like method for the time evolution. The proposed numerical method produces results which, in the middle part of the
considered droplets, qualitatively similar to those reported in \cite{chertock2012sinking}. Using the new method, we were able to capture
complicated dynamics of the bacteria cells including emergence of plumes and their evolution in complex droplet geometries.


The rest of the paper is organized as follows. In \S\ref{sec2}, we describe a non-dimensional version of the coupled chemotaxis-fluid system
\eref{1.1} and introduce typical values of the scaling parameters to be used in our numerical simulations. In \S\ref{sec3}, we present the
reformulated cf-DD model. In \S\ref{sec4}, we introduce the numerical method for the cf-DD system and discuss its implementation. In
\S\ref{sec5}, we report several numerical experiments illustrating a superb performance of the proposed diffuse-domain based numerical
method. Finally, in \S\ref{sec6}, we give few concluding remarks and discuss perspectives of our future work.

\section{Scaling and Setup}\label{sec2}
We denote by $L$ a characteristic length (we may choose, for instance, $L$ to be the maximum height of the drop; see Figure \ref{OD}) and
the characteristic cell density by $n_r$. Rescaling the variables as in \cite{chertock2012sinking,tuval2005bacterial},
\begin{equation}
\bm x'=\frac{\bm x}{L},\quad t'=\frac{D_n}{L^2}t,\quad c'=\frac{c}{c_{air}},\quad n'=\frac{n}{n_r},\quad p'=\frac{L^2}{\eta D_n}p,\quad
\bm u'=\frac{L}{D_n}\bm u,
\label{2.1}
\end{equation}
leads to the five dimensionless parameters $\alpha$, $\beta$, $\gamma$, $\delta$ and the Schmidt number $\Sc$:
\begin{equation}
\alpha:=\frac{\chi c_{air}}{D_n},\quad\beta:=\frac{\kappa n_rL^2}{c_{air}D_n},\quad\gamma:=\frac{V_bn_rg(\rho_b-\rho)L^3}{\eta D_n},\quad
\delta:=\frac{D_c}{D_n},\quad \Sc:=\frac{\eta}{D_n\rho},
\label{2.2}
\end{equation}
which characterize the system \eref{1.1}--\eref{1.3}. Three of the parameters in \eref{2.2}, namely $\alpha$, $\delta$ and $\Sc$, are
determined by the properties of bacteria, fluid and air. Typical values for Bacillus subtilis in water are $\alpha=10$, $\delta=5$ and
$\Sc=500$; see, e.g., \cite{tuval2005bacterial}. The remaining two parameters $\beta$ and $\gamma$ depend also on the chosen length scale $L$
and the reference cell density $n_r$, and thus will be varied in the numerical examples reported in \S\ref{sec5}.

Dropping the primes from the dimensionless quantities in \eref{2.1} yields the following non-dimensional version of the governing
chemotaxis-fluid system:
\begin{align}
&n_t+\nabla\!\cdot\!(\bm{u}n)+\alpha\nabla\!\cdot\![r(c)n\nabla c]=\Delta n,\label{2.3a}\\
&c_t+\bm u\!\cdot\!\nabla c=\delta\Delta c-\beta r(c)n,\label{2.3b}\\
&\bm u_t+\bm u\!\cdot\!\nabla\bm u+\Sc \, \nabla p=\Sc\,\Delta\bm u-\Sc\,\gamma n\bm z,\label{2.3c}\\
&\nabla\!\cdot\!\bm u=0.\label{2.3d}
\end{align}
This system is considered on a sessile drop domain $\Omega$ subject to the initial data
\begin{equation}
n(x,y,0)=n_0(x,y),\quad c(x,y,0)=c_0(x,y),\quad\bm{u}(x,y,0)=\bm{u}_0(x,y)
\label{2.7}
\end{equation}
and the following boundary conditions:
\begin{align}
&\bm\nu\!\cdot\!\bm u=0,\quad\bm\nu\!\cdot\!\nabla(\bm u\!\cdot\!\bm\tau)=0,\quad
\left(\alpha n\nabla c-\nabla n\right)\!\cdot\!\bm{\nu}=0,\quad c=1,&&\forall(x,y)\in\Gamma,\label{2.8}\\
&\bm u=\bm0,\quad\nabla n\!\cdot\!\bm\nu=0,\quad\nabla c\!\cdot\!\bm\nu=0,&&\forall(x,y)\in\partial\Omega_{bot}.
\label{2.9}
\end{align}

\section{Diffuse-Domain Reformulation}\label{sec3}
\subsection{Chemotaxis-Fluid Diffuse-Domain Model}\label{sec3.1}
\begin{figure}[ht!]
\centerline{\includegraphics[trim=0.7cm 0.9cm 0.1cm 1.3cm, clip,width=11cm]{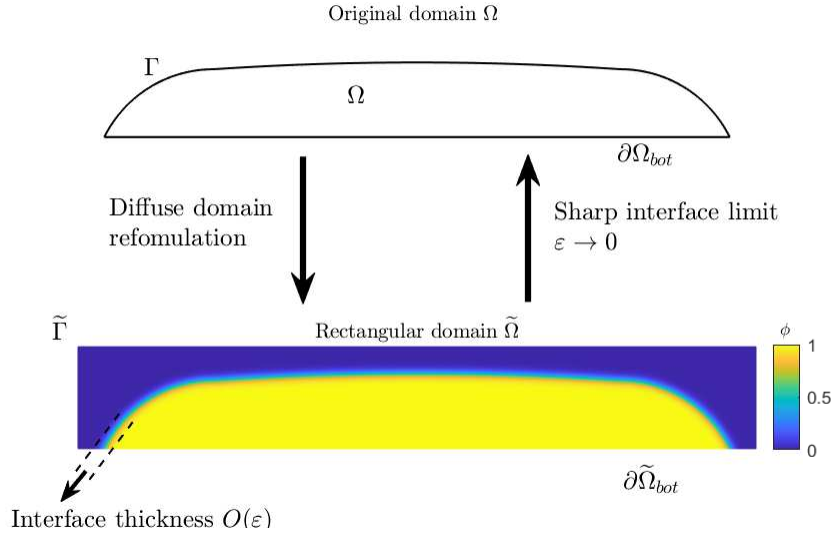}}
\caption{\sf Schematic representation of the diffuse-domain method. The original domain $\Omega$ is embedded in a larger, rectangular domain
$\widetilde\Omega$, where a diffuse-domain function $\phi$ approximates the characteristic function of $\Omega$. The boundary conditions on
$\partial\widetilde\Omega_{bot}$ are the same as those prescribed for the original system on $\partial\Omega_{bot}$, while the boundary
conditions on $\widetilde\Gamma$ are chosen to be consistent with those on $\Gamma$.\label{DD}}
\end{figure}

In order to numerically solve the coupled chemotaxis-fluid system \eref{2.3a}--\eref{2.9} in the drop domain, we propose a diffuse domain
approximation of the chemotaxis-fluid model in a larger rectangular domain $\widetilde\Omega$ outlined in Figure \ref{DD}. The cf-DD model
reads as
\begin{align}
&\phi n_t+\nabla\!\cdot\!(\phi\bm u n)+\alpha\nabla\!\cdot\![r(c)\phi n\nabla c]=\nabla\!\cdot\!(\phi\nabla n)+{\cal B}_n,\label{3.1}\\
&\phi c_t+\phi\bm u\!\cdot\!\nabla c=\delta\nabla\!\cdot\!(\phi\nabla c)-\beta r(c)\phi n+{\cal B}_c,\label{3.2}\\
&\phi\bm u_t+\phi\bm u\!\cdot\!\nabla\bm u+\Sc\,\phi\nabla p=\Sc\,\nabla\!\cdot\!(\phi\nabla\bm u)-\Sc\,\gamma\phi n\bm z+\bm{{\cal B}}_{\bm u},
\label{3.3}\\
&\nabla\!\cdot\!(\phi\bm u)=0,\label{3.4}
\end{align}
with the boundary conditions
\begin{align}
&\bm\nu\!\cdot\!\nabla\bm u=\bm0,\quad\left(\alpha n\nabla c-\nabla n\right)\!\cdot\!\bm\nu=0,\quad c=1,&&\forall(x,y)\in\widetilde{\Gamma},
\label{3.6}\\
&\bm u=\bm0,\quad\nabla n\!\cdot\!\bm\nu=0,\quad\nabla c\!\cdot\!\bm\nu=0,&&\forall(x,y)\in\partial\widetilde\Omega_{bot}.\label{3.5}
\end{align}
We set the diffuse-domain function $\phi$ to be the following approximation of the characteristic function of the original domain $\Omega$:
\begin{equation}
\phi(\bm x)=\hf\left[1-\tanh\left(\frac{3d(\bm x)}{\varepsilon}\right)\right],
\label{3.7}
\end{equation}
where $d(\bm x)$ is the signed distance function to $\Gamma$ ($d<0$ inside $\Omega$) and $\varepsilon$ is the thickness of the diffuse
domain boundary as shown in Figure \ref{DD}. Note that the function $\phi$ is independent of time since the domain $\Omega$ is fixed.
Finally, the terms
\begin{equation}
{\cal B}_n=0,\quad {\cal B}_c=-\frac{1-\phi}{\varepsilon^3}(c-1),\quad\bm{{\cal B}}_{\bm u}=\bm0
\label{3.8}
\end{equation}
are added to enforce the original boundary conditions \eref{2.8} on $\Gamma$, and these terms have been selected following the idea
introduced in \cite{li2009solving}. In \S\ref{sec3.2}, we will show that the cf-DD system \eref{3.1}--\eref{3.8} asymptotically converges to
the original chemotaxis-fluid system \eref{2.3a}--\eref{2.3d} with the boundary conditions \eref{2.8} and \eref{2.9} as $\varepsilon\to0$.

\subsection{Asymptotic Analysis}\label{sec3.2}
We now use the method of matched asymptotic expansions (see, e.g., \cite{AMW,BO,guo2021diffuse,Hol,Peg}) to analyze the cf-DD system
\eref{3.1}--\eref{3.8}. In particular, we expand $n$, $c$, $u$, $v$ and $p$ with respect to the small parameter $\varepsilon$ (representing
the interface thickness according to \eref{3.7}) in regions close to the interface (inner region) and far from the interface (outer region),
which are defined as follows:
$$
\mbox{inner region}=\big\{\bm{x}: |d(\bm x)|<\varepsilon^{s_1}\big\},\quad
\mbox{outer region}=\big\{\bm{x}: |d(\bm x)|>\varepsilon^{s_2}\big\},
$$
where $0<s_1<s_2$. When $\varepsilon$ is small, the inner and outer regions overlap, and the two expansions are to be matched in the
overlapping region $=\big\{\bm x: \varepsilon^{s_2}<|d(\bm x)|<\varepsilon^{s_1}\big\}$.

For the purpose of asymptotic analysis, we consider a smoothed $r(c)$, while in the numerical experiments reported in \S\ref{sec5} the
original formula \eref{rc} has been utilized.

\subsubsection{Outer Expansions}\label{sec3.2.1}
We introduce the vector $\bm{w}:=(n,c,u,v,p)^\top$ and expand it in $\varepsilon$ in the outer region on each side of the interface
$\Gamma$. We denote these formal outer expansions by $\bm{w}^+(\bm{x})$ for $\bm{x}\in\Omega$ at which $\phi(\bm{x})\approx1$ (inside
$\Omega$) and $\bm w^-(\bm x)$ for $\bm x\not\in\Omega$ at which $\phi(\bm x)\approx0$ (outside $\Omega$):
\begin{equation*}
\bm w^\pm=\bm w_0^\pm+\varepsilon\bm w_1^\pm+\varepsilon^2\bm w_2^\pm+\cdots.
\end{equation*}

We first substitute $\bm{w}^+$ into \eref{3.1}--\eref{3.4}. Taking into account that $\phi\approx1$ inside $\Omega$, we combine the leading
${\cal O}(1)$ terms in the resulting expansions and obtain
\begin{equation}
\begin{aligned}
&(n_0^+)_t+\nabla\!\cdot\!(\bm u_0^+n_0^+)+\alpha\nabla\!\cdot\![r(c_0^+)n_0^+\nabla c_0^+]=\Delta n_0^+,\\
&(c_0^+)_t+\bm u_0^+\!\cdot\!\nabla c_0^+=\delta\Delta c_0^+-\beta r(c_0^+)n_0^+,\\
&(\bm u_0^+)_t+\bm u_0^+\!\cdot\!\nabla\bm u_0^++\Sc\,\nabla p_0^+=\Sc\,\Delta\bm u_0^+-\Sc\,\gamma n_0^+\bm z,\\
&\nabla\!\cdot\!\bm u_0^+=0,
\end{aligned}
\label{3.36}
\end{equation}
so that $n_0^+$, $c_0^+$, $\bm u_0^+$ and $p_0^+$ satisfy the chemotaxis-fluid system \eref{2.3a}--\eref{2.3d}. Moreover, substituting
$\bm w^+$ into the boundary condition \eref{3.5}, we can easily see that $n_0^+$, $c_0^+$ and $\bm{u}_0^+$ satisfy the boundary condition
\eref{2.9}. We then perform similar analysis for $\bm w^-$, which results in $c_0^-\equiv1$ since $\phi\approx0$ outside $\Omega$.

\subsubsection{Inner Expansions}\label{sec3.2.2}
We now consider the expansions of $\bm{w}$ in the inner region. To this end, we first use the divergence-free condition \eref{3.4} and
rewrite the system \eref{3.1}--\eref{3.4}, \eref{3.8} in the following equivalent form:
\begin{align}
&\phi n_t+\phi\bm u\!\cdot\!\nabla n+\alpha\nabla\!\cdot\![r(c)\phi n\nabla c]=\nabla\!\cdot\!(\phi\nabla n),\label{3.14}\\
&\phi c_t+\phi\bm u\!\cdot\!\nabla c=\delta\nabla\!\cdot\!(\phi\nabla c)-\beta r(c)\phi n-(1-\phi)\varepsilon^{-3}(c-1),\label{3.15}\\
&\phi\bm u_t+\phi\bm u\!\cdot\!\nabla\bm u+\Sc\,\phi\nabla p=\Sc\,\nabla\!\cdot\!(\phi\nabla\bm u)-\Sc\,\gamma\phi n\bm z,\label{3.16}\\
&\nabla\!\cdot\!(\phi\bm u)=0.\label{3.17}
\end{align}

Next, we introduce a rescaled variable $\xi:=d(\bm x)/\varepsilon$ and a local coordinate system near the interface $\Gamma$:
\begin{equation*}
\bm x(s,\xi;\varepsilon)=\bm X(s)+\varepsilon\xi\bm\nu(s),
\end{equation*}
where $\bm X(s)$ is a parametrization of $\Gamma$, and $s$ is the arc length parameter. We use the notation
$\omega(\bm x)=\widehat{\omega}(s,\xi)$ for any function $\omega$, and notice that the following identities hold:
\begin{align}
&\nabla\omega=(1+\varepsilon\xi\kappa)^{-1}\nabla_\Gamma\widehat{\omega}+\varepsilon^{-1}\bm\nu\widehat{\omega}_\xi,\label{3.10}\\
&\Delta\omega=(1+\varepsilon\xi\kappa)^{-1}\nabla_\Gamma\!\cdot\!\big((1+\varepsilon\xi\kappa)^{-1}\nabla_\Gamma\widehat{\omega}\big)+
\varepsilon^{-1}\kappa(1+\varepsilon\xi\kappa)^{-1}\widehat{\omega}_\xi+\varepsilon^{-2}\widehat{\omega}_{\xi\xi},\label{3.11}
\end{align}
where, the $\nabla_\Gamma$ and $\nabla_\Gamma\cdot$ stand for the curve gradient and curve divergence operators, respectively. In
\eref{3.11}, we have used the facts that $\nabla_\Gamma\!\cdot\!\bm{\nu}=\kappa$, where $\kappa$ is the mean curvature of the interface, and
$\nabla d=\bm\nu$.

We then substitute a formal expansion, which is valid in the inner region,
\begin{equation}
\widehat{\bm w}=\widehat{\bm w}_0+\varepsilon\widehat{\bm w}_1+\varepsilon^2\widehat{\bm w}_2+\cdots.
\label{3.13}
\end{equation}
into the system \eref{3.14}--\eref{3.17}, use relations \eref{3.10}--\eref{3.11} and collect the like powers of $\varepsilon$. At the
leading order term, ${\cal O}(\varepsilon^{-3})$, we obtain
\begin{equation}
(1-\phi)(\widehat c_0-1)=0.
\label{3.18}
\end{equation}
Next, equating the ${\cal O}(\varepsilon^{-2})$ terms in equations \eref{3.14}--\eref{3.16} results in
\begin{align}
&\left[\phi\alpha r(\widehat c_0)\widehat n_0(\widehat c_0)_\xi-\phi(\widehat n_0)_\xi\right]_\xi=0,\label{3.19}\\
&\delta\left[\phi(\widehat c_0)_\xi\right]_\xi-(1-\phi)\widehat c_1=0,\label{3.20}\\
&\left[\phi(\widehat{\bm u}_0)_\xi\right]_\xi=0.\label{3.21}
\end{align}
Finally, balancing the ${\cal O}(\varepsilon^{-1})$ terms in equations \eref{3.14}--\eref{3.17} leads, after some simplifications, to
\begin{align}
&\begin{aligned}
\phi\widehat{\bm u}_0\!\cdot\!\bm\nu(\widehat n_0)_\xi&+\left[\phi\alpha r(\widehat c_0)\left\{\widehat n_0(\widehat c_1)_\xi+
\widehat n_1(\widehat c_0)_\xi\right\}+\phi\alpha r'(\widehat c_0)\widehat c_1\widehat n_0(\widehat c_0)_\xi-
\phi(\widehat n_1)_\xi\right]_\xi\\
&+\kappa\phi\left[\alpha r(\widehat c_0)\widehat n_0(\widehat c_0)_\xi-(\widehat n_0)_\xi\right]=0,
\end{aligned}\label{3.22}\\
&\phi\widehat{\bm u}_0\!\cdot\!\bm\nu(\widehat c_0)_\xi+\delta\left[\phi(\widehat c_1)_\xi\right]_\xi+\delta\kappa\phi(\widehat c_0)_\xi-
(1-\phi)\widehat c_2=0,\label{3.23}\\
&\phi\widehat{\bm u}_0\!\cdot\!\bm\nu(\widehat{\bm u}_0)_\xi+\Sc\big\{\phi(\widehat p_0)_\xi\bm\nu-
\left[\phi(\widehat{\bm{u}}_1)_\xi\right]_\xi-(\phi\bm\nu\!\cdot\!\nabla_\Gamma\widehat{\bm u}_0)_\xi
-\nabla_\Gamma\!\cdot\!(\phi\bm\nu)(\widehat{\bm u}_0)_\xi\big\}=\bm0,\label{3.24}\\
&\bm\nu\!\cdot\!(\phi\widehat{\bm u}_0)_\xi=0.\label{3.25}
\end{align}

\subsubsection{Matching Conditions in the Overlapping Region}\label{sec3.2.3}
In what follows, we derive the boundary conditions on $\Gamma$ by matching the outer and inner expansions in the overlapping region. To this
end, the following matching conditions at $\bm x=\bm X(s)\in\Gamma$ needs to be satisfied (see \cite{Abels2012,guo2021diffuse}):
\begin{align}
&\lim_{\xi\to{\mp\infty}}\widehat{\bm w}_0(s,\xi)=\bm w_0^\pm(\bm x\pm),\label{3.34}\\
&\lim_{\xi\to{\mp\infty}}\big(\widehat{\bm w}_1(s,\xi)\big)_\xi=\bm\nu\!\cdot\!\nabla\bm w_0^\pm(\bm x\pm),
\label{3.30}
\end{align}
where $\bm w_0^\pm(\bm x\pm)$ denote the limits $\lim_{h\to0^-}\bm w_0^\pm(\bm x\pm h\bm\nu)$.

First, we note that \eref{3.18} implies $\widehat c_0(s,\xi)\equiv1$, which together with the matching condition \eref{3.34} imply
\begin{equation}
c_0^+({\bm x}+)=\widehat c_0\equiv1~~\mbox{at}~~\bm x\in\Gamma.
\label{3.29}
\end{equation}
We then use \eref{3.29} to rewrite \eref{3.19} as
\begin{equation*}
\left[\phi(\widehat n_0)_\xi\right]_\xi=0,
\end{equation*}
which implies
\begin{equation*}
\phi(\widehat n_0)_\xi=C(s),
\end{equation*}
and since $\lim_{\xi\to\infty}\phi=0$, we conclude that $C(s)\equiv0$, and hence
\begin{equation*}
(\widehat n_0)_\xi=0.
\end{equation*}
This together with \eref{3.29} allows us to rewrite \eref{3.22} as
\begin{equation*}
\left[\phi\big(\alpha\widehat n_0(\widehat c_1)_\xi-(\widehat n_1)_\xi\big)\right]_\xi=0,
\end{equation*}
and therefore,
\begin{equation*}
\alpha\widehat n_0(\widehat c_1)_\xi-(\widehat n_1)_\xi=0,
\end{equation*}
which, using the matching conditions \eref{3.34} and \eref{3.30}, reduces to the following condition on $n_0^+$ and $c_0^+$:
\begin{equation}
(\alpha n_0^+\nabla c_0^+-\nabla n_0^+)\!\cdot\!\bm\nu=0.
\label{3.31}
\end{equation}

Similarly, we use \eref{3.21} to obtain
\begin{equation*}
(\widehat{\bm u}_0)_\xi=\bm0,
\end{equation*}
which allows us to rewrite \eref{3.24} as
\begin{equation*}
\phi(\widehat p_0)_\xi\bm\nu-\left[\phi(\widehat{\bm u}_1)_\xi\right]_\xi=\bm0.
\end{equation*}
After applying the orthogonal projection operator $P_\Gamma:=I-\bm\nu\otimes\bm\nu:~\widetilde\Omega\to\Gamma$, where $I$ is the
identity matrix, the last equation further reduces to
\begin{equation*}
\left[\phi(\widehat{\bm u}_1\!\cdot\!\bm\tau)_\xi\right]_\xi=0.
\end{equation*}
We then proceed with the arguments similar to those used to derive \eref{3.31} and conclude with
\begin{equation}
\bm\nu\!\cdot\!\nabla(\bm u_0^+\!\cdot\!\bm\tau)=0.
\label{3.35}
\end{equation}
Also note that equation \eref{3.25} and the matching condition \eref{3.34} lead to
\begin{equation}
\bm{\nu}\!\cdot\!\bm{u}_0^+=0.
\label{3.33}
\end{equation}

Finally, \eref{3.29}--\eref{3.33} together with \eref{3.36} imply that the cf-DD system \eref{3.1}--\eref{3.8} asymptotically converges to
the chemotaxis-fluid system \eref{2.3a}--\eref{2.3d} with the boundary conditions \eref{2.8} and \eref{2.9} as $\varepsilon\to0$.

\section{Hybrid Finite-Volume Finite-Difference Method}\label{sec4}
Recall that the advantage of the cf-DD system \eref{3.1}--\eref{3.8} is that it is posed on a simple, rectangular domain and thus it can be
numerically solved in a much easier way compared to the original chemotaxis-fluid system. In this section, we provide a detailed description
of the hybrid finite-volume finite-difference numerical method used to solve the studied cf-DD system.

The cell density equation \eref{3.1} will be solved using a semi-discrete second-order finite-volume upwind scheme combined with a
second-order strong stability-preserving multistep ODE solver for the temporal discretization. The oxygen concentration equation \eref{3.2}
and the Navier-Stokes fluid equations will be discretized using a semi-discrete finite-difference central scheme combined with a
second-order BDF-like method for time integration and a second-order projection method needed to enforce the divergence-free property of the
fluid velocity.

\subsection{Finite-Volume Upwind Scheme for the Cell Density Equation}\label{sec41}
We first define $m:=\phi n$, substitute \eref{3.8} into the cell density equation \eref{3.1} and rewrite it in an equivalent coordinate
form:
\begin{equation}
m_t+\big[(u+\alpha r(c)c_x)m\big]_x+\big[(v+\alpha r(c)c_y)m\big]_y=
\Big[\phi\Big(\frac{m}{\phi}\Big)_x\Big]_x+\Big[\phi\Big(\frac{m}{\phi}\Big)_y\Big]_y.
\label{4.2}
\end{equation}
We then discretize equation \eref{4.2} in space using the semi-discrete second-order finite-volume upwind scheme from
\cite{chertock2012sinking}.

To this end, we divide the computational domain $\widetilde\Omega$ into the cells $I_{j,k}:=[x_\jmh,x_\jph]\times[y_\kmh,y_\kph]$ centered
at $(x_j,y_k)=\big(\nicefrac{(x_\jmh+x_\jph)}{2},\nicefrac{(y_\kmh+y_\kph)}{2}\big)$ with $j=1,\ldots,N$ and $k=1,\ldots,M$. For simplicity,
we use a uniform mesh with $x_\jph-x_\jmh\equiv\dx$ and $y_\kph-y_\kmh\equiv\dy$, where $\dx$ and $\dy$ are small spatial scales. We then
denote the cell averages of $m$ by
\begin{equation*}
\xbar m_{j,k}(t)\approx\frac{1}{\dx\dy}\iint\limits_{I_{j,k}}m(x,y,t)\,{\rm d}x\,{\rm d}y,
\end{equation*}
and integrate equation \eref{4.2} over cell $I_{j,k}$ to obtain
\begin{equation*}
\begin{aligned}
\dx\dy\frac{\rm d}{{\rm d}t}\,\xbar m_{j,k}(t)&+\int\limits_{y_\kmh}^{y_\kph}\big(u+\alpha r(c)c_x\big)m\bigg|_{x_\jmh}^{x_\jph}\,{\rm d}y+
\int\limits_{x_\jmh}^{x_\jph}\big(v+\alpha r(c)c_y\big)m\bigg|_{y_\kmh}^{y_\kph}\,{\rm d}x\\
&=\int\limits_{y_\kmh}^{y_\kph}\phi\Big(\frac{m}{\phi}\Big)_x\bigg|_{x_\jmh}^{x_\jph}\,{\rm d}y+
\int\limits_{x_\jmh}^{x_\jph}\phi\Big(\frac{m}{\phi}\Big)_y\bigg|_{y_\kmh}^{y_\kph}\,{\rm d}x.
\end{aligned}
\end{equation*}
Applying the midpoint rule to the above integrals and dividing by $\dx\dy$ results in
\begin{equation}
\begin{aligned}
\frac{\rm d}{{\rm d}t}\,\xbar m_{j,k}=&-\frac{\big(u+\alpha r(c)c_x\big)m\Big|_{(x_\jph,y_k)}-
\big(u+\alpha r(c)c_x\big)m\Big|_{(x_\jmh,y_k)}}{\dx}\\
&-\frac{\big(v+\alpha r(c)c_y\big)m\Big|_{(x_j,y_\kph)}-\big(v+\alpha r(c)c_y)m\Big|_{(x_j,y_\kmh)}}{\dy}\\
&+\frac{\phi\Big(\dfrac{m}{\phi}\Big)_x\Big|_{(x_\jph,y_k)}-\phi\Big(\dfrac{m}{\phi}\Big)_x\Big|_{(x_\jmh, y_k)}}{\dx}+
\frac{\phi\Big(\dfrac{m}{\phi}\Big)_y\Big|_{(x_j,y_\kph)}-\phi\Big(\dfrac{m}{\phi}\Big)_y\Big|_{(x_j,y_\kmh)}}{\dy}.
\end{aligned}
\label{4.5}
\end{equation}
We note that $\,\xbar m_{j,k}$ as well as many other indexed quantities in \eref{4.5} and below depend on time $t$, but from now on we omit
this dependence for the sake of brevity.

The construction of the scheme will be completed once the fluxes at the cell interfaces in \eref{4.5} are approximated numerically. The
semi-discrete finite-volume upwind scheme can then be written as the following system of ODEs:
\begin{equation}
\frac{{\rm d}}{{\rm d}t}\,\xbar m_{j,k}=-\frac{F_{\jph,k}^x-F_{\jmh,k}^x}{\dx}-\frac{F_{j,\kph}^y-F_{j,\kmh}^y}{\dy}+
\frac{G_{\jph,k}^x-G_{\jmh,k}^x}{\dx}+\frac{G_{j,\kph}^y-G_{j,\kmh}^y}{\dy},
\label{4.6}
\end{equation}
where
\begin{equation}
F_{j\pm\hf,k}^x\approx\big(u+\alpha r(c)c_x\big)m\Big|_{(x_{j\pm\hf},y_k)}\quad\mbox{and}\quad
F_{j,k\pm\hf}^y\approx\big(v+\alpha r(c)c_y\big)m\Big|_{(x_j,y_{k\pm\hf})}
\label{4.7}
\end{equation}
are numerical convection-chemotaxis fluxes, and
\begin{equation}
G_{j\pm\hf,k}^x\approx\phi\Big(\frac{m}{\phi}\Big)_x\bigg|_{(x_{j\pm\hf},y_k)}\quad\mbox{and}\quad
G_{j,k\pm\hf}^y\approx\phi\Big(\frac{m}{\phi}\Big)_y\bigg|_{(x_j,y_{k\pm\hf})}
\label{4.8}
\end{equation}
are centered numerical diffusion fluxes.

In order to ensure stability of the scheme \eref{4.6}--\eref{4.8}, we use an upwind approximation of the convection-chemotaxis fluxes, which
can be written in the following form:
\begin{equation}
F_{\jph,k}^x=\left\{
\begin{aligned}&a_{\jph,k}m_{j,k}^{\rm E}&&\mbox{if}~~a_{\jph,k}\ge0,\\&a_{\jph,k}m_{j+1,k}^{\rm W}&&\mbox{if}~~a_{\jph,k}<0,\end{aligned}
\right.\qquad F_{j,\kph}^y=\left\{
\begin{aligned}&b_{j,\kph}m_{j,k}^{\rm N}&&\mbox{if}~~b_{j,\kph}\ge0,\\&b_{j,\kph}m_{j,k+1}^{\rm S}&&\mbox{if}~~b_{j,\kph}<0.\end{aligned}
\right.
\label{4.9}
\end{equation}
Here, $m_{j,k}^{\rm E,W,N,S}$ are the point values of the piecewise linear reconstruction consisting of the following linear pieces on every
interval $I_{j,k}$:
\begin{equation}
\widetilde m_{j,k}(x,y)=\,\xbar m_{j,k}+(m_x)_{j,k}(x-x_j)+(m_y)_{j,k}(y-y_k),\quad(x,y)\in I_{j,k},
\label{4.10}
\end{equation}
at the points $(x_\jph,y_k)$, $(x_\jmh,y_k)$, $(x_j,y_\kph)$, and $(x_j,y_\kmh)$, respectively. Namely, we have
\begin{equation}
\begin{aligned}
&m_{j,k}^{\rm E}=\widetilde m_{j,k}(x_\jph,y_k)=\,\xbar m_{j,k}+\frac{\dx}{2}(m_x)_{j,k},\\
&m_{j,k}^{\rm W}=\widetilde m_{j,k}(x_\jmh,y_k)=\,\xbar m_{j,k}-\frac{\dx}{2}(m_x)_{j,k},\\
&m_{j,k}^{\rm N}=\widetilde m_{j,k}(x_j,y_\kph)=\,\xbar m_{j,k}+\frac{\dy}{2}(m_y)_{j,k},\\
&m_{j,k}^{\rm S}=\widetilde m_{j,k}(x_j,y_\kmh)=\,\xbar m_{j,k}-\frac{\dy}{2}(m_y)_{j,k}.
\end{aligned}
\label{4.11}
\end{equation}
The second order of accuracy will be guaranteed provided the numerical derivatives $(m_x)_{j,k}$ and $(m_y)_{j,k}$ are to be (at least)
first-order approximations of the corresponding exact derivatives $m_x(x_j,y_k,t)$ and $m_y(x_j,y_k,t)$. In our numerical experiments, we
have used the central-difference approximations,
\begin{equation}
(m_x)_{j,k}=\frac{\,\xbar m_{j+1,k}-\,\xbar m_{j-1,k}}{2\dx}\quad\mbox{and}\quad
(m_y)_{j,k}=\frac{\,\xbar m_{j,k+1}-\,\xbar m_{j,k-1}}{2\dy},
\label{4.13}
\end{equation}
throughout the computational domain except for the cells, where the linear approach \eref{4.13} leads to the appearance of negative
reconstructed values of $m$ in \eref{4.11}. In the cells, where either $m_{j,k}^{\rm E}$ or $m_{j,k}^{\rm W}$ is negative, we replace
\eref{4.13} with a nonlinear minmod2 reconstruction (see, e.g., \cite{NT,Swe,vLeV}):
\begin{equation}
(m_x)_{j,k}={\rm minmod}\left(2\,\frac{\,\xbar m_{j,k}-\,\xbar m_{j-1,k}}{\dx},\,\frac{\,\xbar m_{j+1,k}-\,\xbar m_{j-1,k}}{2\dx},\,
2\,\frac{\,\xbar m_{j+1,k}-\,\xbar m_{j,k}}{\dx}\right),
\label{4.14}
\end{equation}
which guarantees that no negative values of $m$ emerge in \eref{4.11}. We then recalculate the reconstructed values $m_{j,k}^{\rm E}$ and
$m_{j,k}^{\rm W}$. Similarly, if either $m_{j,k}^{\rm N}$ or $m_{j,k}^{\rm S}$ is negative, we set
\begin{equation}
(m_y)_{j,k}={\rm minmod}\left(2\,\frac{\,\xbar m_{j,k}-\,\xbar m_{j,k-1}}{\dy},\,\frac{\,\xbar m_{j,k+1}-\,\xbar m_{j,k-1}}{2\dy},\,
2\,\frac{\,\xbar m_{j,k+1}-\,\xbar m_{j,k}}{\dy}\right),
\label{4.15}
\end{equation}
and recalculate the reconstructed values $m_{j,k}^{\rm N}$ and $m_{j,k}^{\rm S}$. The minmod function used in \eref{4.14} and \eref{4.15} is
defined as
\begin{equation}
{\rm minmod}(z_1,z_2,\ldots):=\left\{
\begin{aligned}&\min_j\{z_j\}&&\mbox{if}~~z_j>0~\forall j,\\&\max_j\{z_j\}&&\mbox{if}~~z_j<0~\forall j,\\&0&&\mbox{otherwise.}
\end{aligned}\right.
\label{4.16}
\end{equation}
The description of the numerical convection-chemotaxis fluxes in \eref{4.9} will be completed once the local speeds in the $x$- and
$y$-directions, $a_{\jph,k}$ and $b_{j,\kph}$, are specified. Since all of the solution components are expected to be smooth, the local
speeds can be approximated using the centered differences and averages as
\begin{equation*}
a_{\jph,k}=u_{\jph,k}+\alpha r(c_{\jph,k})(c_x)_{\jph,k}\quad\mbox{and}\quad b_{j,\kph}=v_{j,\kph}+\alpha r(c_{j,\kph})(c_y)_{j,\kph},
\end{equation*}
where
$$
\begin{aligned}
(c_x)_{\jph,k}&=\frac{c_{j+1,k}-c_{j,k}}{\dx},&u_{\jph,k}&=\hf\left(u_{j,k}^{\rm E}+u_{j+1,k}^{\rm W}\right),&
c_{\jph,k}&=\hf\left(c_{j,k}^{\rm E}+c_{j+1,k}^{\rm W}\right),\\
(c_y)_{j,\kph}&=\frac{c_{j,k+1}-c_{j,k}}{\dy},&v_{j,\kph}&=\hf\left(v_{j,k}^{\rm N}+v_{j,k+1}^{\rm S}\right),&
c_{j,\kph}&=\hf\left(c_{j,k}^{\rm N}+c_{j,k+1}^{\rm S}\right).
\end{aligned}
$$
Here, the point values $c_{j,k}^{\rm E,W,N,S}$, $u_{j,k}^{\rm E,W}$ and $v_{j,k}^{\rm N,S}$ are obtained using the same piecewise linear
reconstruction, which was used to compute the corresponding values of $m$ in \eref{4.11}, but now applied to the point values
$c_{j,k}\approx c(x_j,y_k,t)$, $u_{j,k}\approx u(x_j,y_k,t)$ and $v_{j,k}\approx v(x_j,y_k,t)$, respectively.

Finally, the centered numerical diffusion fluxes in \eref{4.8} are approximated by
\begin{equation}
G_{\jph,k}^x=\frac{\phi_{\jph,k}}{\dx}\left(\frac{\,\xbar m_{j+1,k}}{\phi_{j+1,k}}-\frac{\,\xbar m_{j,k}}{\phi_{j,k}}\right)\quad\mbox{and}
\quad G_{j,\kph}^y=\frac{\phi_{j,\kph}}{\dy}\left(\frac{\,\xbar m_{j,k+1}}{\phi_{j,k+1}}-\frac{\,\xbar m_{j,k}}{\phi_{j,k}}\right),
\label{4.20}
\end{equation}
where $\phi_{j,k}:=\phi(x_j,y_k)$, $\phi_{\jph,k}:=\phi(x_\jph,y_k)$, and $\phi_{j,\kph}:=\phi(x_j,y_\kph)$.

\paragraph{Time Discretization.} The semi-discretization \eref{4.6} results in the system of time-dependent ODEs, which we integrate using
the second-order strong stability-preserving (SSP) three-step method \cite{gottlieb2001strong}. This results in
\begin{equation}
\begin{aligned}
\xbar m_{j,k}^{\,\ell+1}&=\frac{3}{4}\,\xbar m_{j,k}^{\,\ell}-\frac{3}{2}\lambda\left(F_{\jph,k}^{x,\ell}-F_{\jmh,k}^{x,\ell}\right)-
\frac{3}{2}\mu\left(F_{j,\kph}^{y,\ell}-F_{j,\kmh}^{y,\ell}\right)\\
&+\frac{3}{2}\lambda\left(G_{\jph,k}^{x,\ell}-G_{\jmh,k}^{x,\ell}\right)+
\frac{3}{2}\mu\left(G_{j,\kph}^{y,\ell}-G_{j,\kmh}^{y,\ell}\right)+\frac{1}{4}\,\xbar m_{j,k}^{\,\ell-2},
\end{aligned}
\label{4.22}
\end{equation}
where $\dt$ is the time step, $\lambda:=\dt/\dx$, $\mu:=\dt/\dy$, $t^\ell:=\ell\dt$, $\,\xbar m_{j,k}^{\,\ell}:=\,\xbar m_{j,k}(t^\ell)$,
$F_{\jph,k}^{x,\ell}:=F_{\jph,k}^x(t^\ell)$, $F_{\jph,k}^{y,\ell}:=F_{\jph,k}^y(t^\ell)$, $G_{j,\kph}^{x,\ell}:=G_{j,\kph}^x(t^\ell)$, and
$G_{j,\kph}^{y,\ell}:=G_{j,\kph}^y(t^\ell)$.

The resulting fully discrete scheme \eref{4.22} is positivity preserving in the sense that $\,\xbar m_{j,k}^{\,\ell+1}\ge0$ for all $j,k$
provided $\,\xbar m_{j,k}^{\,\ell}\ge0$ and $\,\xbar m_{j,k}^{\,\ell-2}\ge0$ for all $j,k$ and $\dt$ is sufficiently small. In order to
prove this, we first note that the convection-chemotaxis numerical fluxes \eref{4.9} can be rewritten as
\begin{equation}
F_{\jph,k}^{x,\ell}=a_{\jph,k}\left(\frac{1+{\rm sign}(a_{\jph,k})}{2}\,m_{j,k}^{\rm E}+
\frac{1-{\rm sign}(a_{\jph,k})}{2}\,m_{j+1,k}^{\rm W}\right)
\label{4.23}
\end{equation}
and
\begin{equation}
F_{j,\kph}^{y,\ell}=b_{j,\kph}\left(\frac{1+{\rm sign}(b_{j,\kph})}{2}\,m_{j,k}^{\rm N}+
\frac{1-{\rm sign}(b_{j,\kph})}{2}\,m_{j,k+1}^{\rm S}\right),
\label{4.24}
\end{equation}
and by the conservation property of the piecewise-linear reconstruction \eref{4.10} the identity
\begin{equation}
\xbar m_{j,k}^{\,\ell}=\frac{1}{8}\left(m_{j,k}^{\rm E}+m_{j,k}^{\rm W}+m_{j,k}^{\rm N}+m_{j,k}^{\rm S}\right)+
\frac{1}{2}\,\xbar m_{j,k}^{\,\ell}
\label{4.16a}
\end{equation}
holds. Note that the quantities $a_{\jph,k}$, $b_{j,\kph}$ and $m_{j,k}^{\rm E,W,N,S}$ in \eref{4.23}--\eref{4.16a} are evaluated at time
level $t=t^\ell$. We then substitute \eref{4.20} and \eref{4.23}--\eref{4.16a} into \eref{4.22} to obtain
\begin{equation}
\begin{aligned}
\xbar m_{j,k}^{\,\ell+1}=\frac{3}{4}\Bigg\{&\bigg[\frac{1}{8}-\lambda\big|a_{\jph,k}\big|\big(1+{\rm sign}(a_{\jph,k})\big)\bigg]
m_{j,k}^{\rm E}+\lambda\big|a_{\jph,k}\big|\big(1-{\rm sign}(a_{\jph,k})\big)m_{j+1,k}^{\rm W}\\
+&\bigg[\frac{1}{8}-\lambda\big|a_{\jmh,k}\big|\big(1-{\rm sign}(a_{\jmh,k})\big)\bigg]m_{j,k}^{\rm W}+
\lambda\big|a_{\jmh,k}\big|\big(1+{\rm sign}(a_{\jmh,k})\big)m_{j-1,k}^{\rm E}\\
+&\bigg[\frac{1}{8}-\mu\big|b_{j,\kph}\big|\big(1+{\rm sign}(b_{j,\kph})\big)\bigg]m_{j,k}^{\rm N}+
\mu\big|b_{j,\kph}\big|\big(1-{\rm sign}(b_{j,\kph})\big)m_{j,k+1}^{\rm S}\\
+&\bigg[\frac{1}{8}-\mu\big|b_{j,\kmh}\big|\big(1-{\rm sign}(b_{j,\kmh})\big)\bigg]m_{j,k}^{\rm S}+
\mu\big|b_{j,\kmh}\big|\big(1+{\rm sign}(b_{j,\kmh})\big)m_{j,k-1}^{\rm N}\Bigg\}\\
+\frac{3}{2}\Bigg\{&\bigg(\frac{1}{4}-\dt\bigg[\frac{\phi_{\jph,k}+\phi_{\jmh,k}}{\phi_{j,k}(\dx)^2}+
\frac{\phi_{j,\kph}+\phi_{j,\kmh}}{\phi_{j,k}(\dy)^2}\bigg]\bigg)\,\xbar m_{j,k}^{\,\ell}\\
+&\frac{\dt}{(\dx)^2}\bigg(\frac{\phi_{\jph,k}}{\phi_{j+1,k}}\,\xbar m_{j+1,k}^{\,\ell}+
\frac{\phi_{\jmh,k}}{\phi_{j-1,k}}\,\xbar m_{j-1,k}^{\,\ell}\bigg)\\
+&\frac{\dt}{(\dy)^2}\bigg(\frac{\phi_{j,\kph}}{\phi_{j,k+1}}\,\xbar m_{j,k+1}^{\,\ell}+
\frac{\phi_{j,\kmh}}{\phi_{j,k-1}}\,\xbar m_{j,k-1}^{\,\ell}\bigg)\Bigg\}+\frac{1}{4}\,\xbar m_{j,k}^{\,\ell-2}.
\end{aligned}
\label{4.26}
\end{equation}
As one can see from \eref{4.26}, the new values $\{\,\xbar m_{j,k}^{\,\ell+1}\}$ are linear combinations of the non-negative cell averages
$\{\,\xbar m_{j,k}^{\,\ell}\}$, $\{\,\xbar m_{j,k}^{\,\ell-2}\}$ and the reconstructed point value $\{m_{j,k}^{\rm E,W,N,S}\}$, which are
also non-negative since they are computed using the positivity preserving piecewise linear reconstruction \eref{4.11}--\eref{4.16}. Thus, as
long as the following CFL condition is satisfied:
\begin{equation}
\begin{aligned}
\dt&\le\frac{1}{16}\min\bigg\{\frac{\dx}{a_{\max}},\,\frac{\dy}{b_{\max}},\,
\frac{4\,\phi_{j,k}(\dx)^2(\dy)^2}{\big(\phi_{j,\kph}+\phi_{j,\kmh}\big)(\dx)^2+\big(\phi_{\jph,k}+\phi_{\jmh,k}\big)(\dy)^2}\bigg\},
\end{aligned}
\label{4.27}
\end{equation}
where
\begin{equation}
a_{\max}:=\max_{j,k}\left\{\big|a_{\jph,k}\big|\right\},\quad b_{\max}:=\max_{j,k}\left\{\big|b_{j,\kph}\big|\right\},
\label{4.20h}
\end{equation}
the linear combination in \eref{4.26} is a convex combination, which implies the non-negativity of $\,\xbar m_{j,k}^{\,\ell+1}$ for all
$j,k$.

Finally, since $m=\phi n$ and $\phi(\bm x)>0$ for all $\bm x$, we conclude that $\,\xbar n_{j,k}^{\,\ell+1}\ge0$ for all $j,k$.
\begin{rmk}\label{rmk41}
It should be observed that the inequality \eref{4.27} should be satisfied at every time level $t=t^\ell$, but since we use the three-step
time discretization method, we have to choose a fixed $\dt$ at time $t=0$, when the data require to be used to evaluate the maxima in
\eref{4.20h} are not available yet. We therefore replace $a_{\max}$ and $b_{\max}$ in \eref{4.27} with their a-priori upper bounds, which
should be valid for all $t$ and can be obtained for any problem at hand.
\end{rmk}
\begin{rmk}
We note that we obtain $\,\xbar n_{j,k}^{\,1}$ and $\,\xbar n_{j,k}^{\,2}$ at the first two time steps using the first-order forward Euler
time discretization.
\end{rmk}

\subsection{Second-Order Projection Finite-Difference Method for the Navier-Stokes and Oxygen Equations}\label{sec42}
Equipped with the obtained values $\,\xbar n_{j,k}^{\,\ell+1}$, we now construct a second-order projection finite-difference method for
equations \eref{3.2}--\eref{3.4} by following the approach from \cite{guo2017mass,shen2012modeling} proposed in the context of the
Cahn-Hilliard-Navier-Stokes system.

We begin with the second-order time discretization of \eref{3.2}--\eref{3.4}, which is based on the projection method and the BDF method
with Adams-Bashforth extrapolation. Assuming that $n^{\ell+1}\approx n(x,y,t^{\ell+1})$, $c^\ell\approx c(x,y,t^\ell)$,
$\bm u^\ell\approx\bm u(x,y,t^\ell)$ and $p^\ell\approx p(x,y,t^\ell)$ are available, we obtain $c^{\ell+1}$, $\bm u^{\ell+1}$ and
$p^{\ell+1}$ by solving the following equations:
\begin{align}
&\phi\,\frac{3\widetilde{\bm u}^{\ell+1}-4\bm u^\ell+\bm u^{\ell-1}}{2\dt}+\phi\bm u^*\!\cdot\!\nabla\bm u^*+\Sc\,\phi\nabla p^\ell=
\Sc\,\nabla\!\cdot\!(\phi\nabla{\widetilde{\bm u}}^{\ell+1})+\Sc\,\gamma\phi n^{\ell+1}\bm z,\label{4.32}\\
&\frac{3(\bm u^{\ell+1}-\widetilde{\bm u}^{\ell+1})}{2\dt}+\Sc\,\nabla\psi^{\ell+1}=\bm0,\label{4.33}\\
&\nabla\!\cdot\!(\phi\bm u^{\ell+1})=0,\label{4.34}\\
&p^{\ell+1}=\psi^{\ell+1}+p^\ell,\label{4.35}\\
&\begin{aligned}
\phi\,\frac{3c^{\ell+1}-4c^\ell+c^{\ell-1}}{2\dt}&+\phi\bm u^{\ell+1}\!\cdot\!\nabla c^*\\
&=\delta\nabla\!\cdot\!(\phi\nabla c^{\ell+1})-\beta r(c^*)\phi n^{\ell+1}-\frac{1}{\varepsilon^3}(1-\phi)(c^{\ell+1}-1),
\end{aligned}\label{4.36}
\end{align}
where $\bm u^*:=2\bm u^\ell-\bm u^{\ell-1}$, $c^*=2c^\ell-c^{\ell-1}$, and $\psi^{\ell+1}$ is an auxiliary variable.

The scheme \eref{4.32}--\eref{4.36} is implemented in the following way. First, we solve the elliptic equation \eref{4.32} for
$\widetilde{\bm u}^{\ell+1}$ subject to the boundary conditions specified in \eref{3.6} and \eref{3.5} for $\bm u$. We multiply both sides
of \eref{4.33} by $\phi$, take the divergence of the result, and use the divergence-free condition \eref{4.34} to obtain the elliptic
equation on $\psi^{\ell+1}$,
\begin{equation}
\nabla\!\cdot\!(\phi\nabla\psi^{\ell+1})=\frac{3}{2\,\Sc\,\dt}\nabla\!\cdot\!(\phi\widetilde{\bm{u}}^{\ell+1}),
\label{4.26h}
\end{equation}
which is solved subject to the homogeneous Neumann boundary condition $\nabla\psi^{\ell+1}\!\cdot\!\bm{\nu}=0$ prescribed on
$\partial\widetilde{\Omega}$. Next, we substitute the computed $\psi^{\ell+1}$ into \eref{4.33} and \eref{4.35} to obtain $\bm{u}^{\ell+1}$
and $p^{\ell+1}$. Finally, we substitute $\bm{u}^{\ell+1}$ and $n^{\ell+1}$ into \eref{4.36} and obtain the elliptic equation for
$c^{\ell+1}$, which is solved subject to the boundary condition specified in \eref{3.6} and \eref{3.5}.
\begin{rmk}
It is shown in \cite{guermond2004error} that the scheme \eref{4.32}--\eref{4.35} without the source term $\Sc\,\gamma\phi n^{\ell+1}\bm z$ is
unconditionally stable.
\end{rmk}
\begin{rmk}
We note that we obtain $\bm u^1$, $p^1$ and $c^1$ at the first time step using the following first-order time discretization, which is based
on the projection method and the backward Euler method:
$$
\begin{aligned}
&\phi\,\frac{\widetilde{\bm{u}}^{\,\ell+1}-\bm{u}^{\,\ell}}{\dt}
+\phi\bm{u}^\ell\!\cdot\!\nabla\bm{u}^\ell+\Sc\,\phi\nabla p^\ell=
\Sc\,\nabla\!\cdot\!(\phi\nabla\widetilde{\bm{u}}^{\,\ell+1})+\Sc\,\gamma\phi n^{\ell+1}\bm{z},\\
&\frac{\bm{u}^{\ell+1}-\widetilde{\bm{u}}^{\,\ell+1}}{\dt}+\Sc\,\nabla(p^{\ell+1}-p^\ell)=\bm{0},\\
&\nabla\!\cdot\!(\phi\bm{u}^{\ell+1})=0,\\
&\phi\,\frac{c^{\ell+1}-c^\ell}{\dt}
+\phi\bm{u}^{\ell+1}\!\cdot\!\nabla c^\ell=\delta\nabla\!\cdot\!(\phi\nabla c^{\ell+1})-
\beta r(c^\ell)\phi n^{\ell+1}-\frac{1}{\varepsilon^3}(1-\phi)(c^{\ell+1}-1).
\end{aligned}
$$
\end{rmk}

\paragraph{Spatial Discretization.} We now denote the point values of $\bm u$, $\psi$, $p$ and $c$ at the cell centers $(x_j,y_k)$ at time
level $t=t^\ell$ by $\bm u_{j,k}^\ell$, $\psi_{j,k}^\ell$, $p_{j,k}^\ell$ and $c_{j,k}^\ell$, respectively and apply the second-order
central difference approximations to construct a fully discrete scheme.

First, we discretize equation \eref{4.32} and use the cell averages $\,\xbar n_{j,k}^{\,\ell+1}$ obtained in \S\ref{sec41} to update
$\widetilde u_{j,k}^{\,\ell+1}$ and $\widetilde v_{j,k}^{\,\ell+1}$ by solving the linear systems
\begin{equation*}
\begin{aligned}
\phi_{j,k}&\left[\frac{3\widetilde u_{j,k}^{\,\ell+1}-4u_{j,k}^\ell+u_{j,k}^{\ell-1}}{2\dt}
+u_{j,k}^*\frac{u^*_{j+1,k}-u^*_{j-1,k}}{2\dx}+v_{j,k}^*\frac{u^*_{j,k+1}-u^*_{j,k-1}}{2\dy}
+\Sc\,\frac{p_{j+1,k}^\ell-p_{j-1,k}^\ell}{2\dx}\right]\\
&=\Sc\,\frac{\phi_{\jph,k}(\widetilde u_{j+1,k}^{\,\ell+1}-\widetilde u_{j,k}^{\,\ell+1})-\phi_{\jmh,k}(\widetilde u_{j,k}^{\,\ell+1}-
\widetilde u_{j-1,k}^{\,\ell+1})}{(\dx)^2}\\
&+\Sc\,\frac{\phi_{j,\kph}(\widetilde u_{j,k+1}^{\,\ell+1}-\widetilde u_{j,k}^{\,\ell+1})-
\phi_{j,\kmh}(\widetilde u_{j,k}^{\,\ell+1}-\widetilde u_{j,k-1}^{\,\ell+1})}{(\dy)^2},
\end{aligned}
\end{equation*}
and
\begin{equation*}
\begin{aligned}
&\begin{aligned}
\phi_{j,k}\left[\frac{3\widetilde v_{j,k}^{\,\ell+1}-4v_{j,k}^\ell+v_{j,k}^{\ell-1}}{2\dt}\right.
&+u_{j,k}^*\frac{v^*_{j+1,k}-v^*_{j-1,k}}{2\dx}+v_{j,k}^*\frac{v^*_{j,k+1}-v^*_{j,k-1}}{2\dy}\\
&\left.+\Sc\,\frac{p_{j,k+1}^\ell-p_{j,k-1}^\ell}{2\dy}-\Sc\,\gamma n_{j,k}^{\ell+1}\right]
\end{aligned}\\
&\quad~=\Sc\,\frac{\phi_{\jph,k}(\widetilde v_{j+1,k}^{\,\ell+1}-\widetilde v_{j,k}^{\,\ell+1})-\phi_{\jmh,k}(\widetilde v_{j,k}^{\,\ell+1}-
\widetilde v_{j-1,k}^{\,\ell+1})}{(\dx)^2}\\
&\quad~+\Sc\,\frac{\phi_{j,\kph}(\widetilde v_{j,k+1}^{\,\ell+1}-\widetilde v_{j,k}^{\,\ell+1})-\phi_{j,\kmh}(\widetilde v_{j,k}^{\,\ell+1}-
\widetilde v_{j,k-1}^{\,\ell+1})}{(\dy)^2},
\end{aligned}
\end{equation*}
for $\{\widetilde u_{j,k}^{\,\ell+1}\}$ and $\{\widetilde v_{j,k}^{\,\ell+1}\}$, respectively. We then discretize equation \eref{4.26h} and
obtain $\psi_{j,k}^{\ell+1}$ by solving the linear system
\begin{equation*}
\begin{aligned}
&\frac{\phi_{\jph,k}(\psi_{j+1,k}^{\ell+1}-\psi_{j,k}^{\ell+1})-\phi_{\jmh,k}(\psi_{j,k}^{\ell+1}-\psi_{j-1,k}^{\ell+1})}{(\dx)^2}\\
&+\frac{\phi_{j,\kph}(\psi_{j,k+1}^{\ell+1}-\psi_{j,k}^{\ell+1})-\phi_{j,\kmh}(\psi_{j,k}^{\ell+1}-\psi_{j,k-1}^{\ell+1})}{(\dy)^2}\\
&=\frac{3}{2\Sc\,\dt}\left(\frac{\phi_{j+1,k}\widetilde u_{j+1,k}^{\,\ell+1}-\phi_{j-1,k}\widetilde u_{j-1,k}^{\,\ell+1}}{2\dx}+
\frac{\phi_{j,k+1}\widetilde v_{j,k+1}^{\,\ell+1}-\phi_{j,k-1}\widetilde v_{j,k-1}^{\,\ell+1}}{2\dy}\right).
\end{aligned}
\end{equation*}
Next, we find $u_{j,k}^{\ell+1}$, $v_{j,k}^{\ell+1}$ and $p_{j,k}^{\ell+1}$ by discretizing \eref{4.33} and \eref{4.35} as follows:
$$
\begin{aligned}
&u_{j,k}^{\ell+1}=\widetilde u_{j,k}^{\,\ell+1}-\frac{2}{3}\Sc\,\dt\frac{\psi_{j+1,k}^{\ell+1}-\psi_{j-1,k}^{\ell+1}}{2\dx},\quad
v_{j,k}^{\ell+1}=\widetilde v_{j,k}^{\,\ell+1}-\frac{2}{3}\Sc\,\dt\frac{\psi_{j,k+1}^{\ell+1}-\psi_{j,k-1}^{\ell+1}}{2\dy},\\
&p_{j,k}^{\ell+1}=
\psi_{j,k}^{\ell+1}+p_{j,k}^\ell.
\end{aligned}
$$
Finally, we discretize \eref{4.36} and update $c_{j,k}^{\ell+1}$ by solving the linear system
\begin{equation*}
\begin{aligned}
&\phi_{j,k}\left[\frac{3c_{j,k}^{\ell+1}-4c_{j,k}^\ell+c_{j,k}^{\ell-1}}{2\dt}+u_{j,k}^{\ell+1}\frac{c^*_{j+1,k}-c^*_{j-1,k}}{2\dx}+
v_{j,k}^{\ell+1}\frac{c^*_{j,k+1}-c^*_{j,k-1}}{2\dy}\right]\\
&=\delta\frac{\phi_{\jph,k}(c_{j+1,k}^{\ell+1}-c_{j,k}^{\ell+1})-\phi_{\jmh,k}(c_{j,k}^{\ell+1}-c_{j-1,k}^{\ell+1})}{(\dx)^2}\\
&+\delta\frac{\phi_{j,\kph}(c_{j,k+1}^{\ell+1}-c_{j,k}^{\ell+1})-\phi_{j,\kmh}(c_{j,k}^{\ell+1}-c_{j,k-1}^{\ell+1})}{(\dy)^2}
-\beta r(c_{j,k}^*)\phi_{j,k}n_{j,k}^{\ell+1}-\frac{1}{\varepsilon^3}(1-\phi_{j,k})(c_{j,k}^{\ell+1}-1).
\end{aligned}
\label{4.46}
\end{equation*}

\subsection{Numerical Boundary Conditions}\label{sec43}
The boundary conditions on $\partial\widetilde\Omega_{bot}$ are given by \eref{3.5}, which is implemented using $m=\phi n$ and the ghost
cell technique as follows:
\begin{equation*}
\xbar m_{j,0}^{\,\ell}=\frac{\phi_{j,0}}{\phi_{j,1}}\,\xbar m_{j,1}^{\,\ell},\quad u_{j,0}^\ell=v_{j,0}^\ell=0,\quad
c_{j,0}^\ell=c_{j,1}^\ell,\quad\forall j,\ell.
\label{4.47}
\end{equation*}

The boundary conditions on $\widetilde\Gamma$ are given by \eref{3.6}. We first rewrite the second equation in \eref{3.6} as
$\partial(\ln n)/\partial\bm\nu=\alpha\partial c/\partial\bm\nu$, which can be easily integrated on each of the three sides of
$\widetilde\Gamma$. We then use $m=\phi n$ and the ghost cell technique to end up with the following boundary conditions for the
corresponding three sides:
\begin{equation*}
\begin{aligned}
&\xbar m_{0,k}^{\,\ell}=\frac{\phi_{0,k}}{\phi_{1,k}}\,\xbar m_{1,k}^{\,\ell}e^{\alpha(1-c_{1,k})},&&u_{0,k}^\ell=u_{1,k}^\ell,&&
v_{0,k}^\ell=v_{1,k}^\ell,&&c_{0,k}^\ell=1,&&\forall k,\ell,\\
&\xbar m_{N+1,k}^{\,\ell}=\frac{\phi_{N+1,k}}{\phi_{N,k}}\,\xbar m_{N,k}^{\,\ell}e^{\alpha(1-c_{N,k})},&&u_{N+1,k}^\ell=u_{N,k}^\ell,&&
v_{N+1,k}^\ell=v_{N,k}^\ell,&&c_{N+1,k}^\ell=1,&&\forall k,\ell,\\
&\xbar m_{j,M+1}^{\,\ell}=\frac{\phi_{j,M+1}}{\phi_{j,M}}\,\xbar m_{j,M}^{\,\ell}e^{\alpha(1-c_{j,M})},&&u_{j,M+1}^\ell=u_{j,M}^\ell,&&
v_{j,M+1}^\ell=v_{j,M}^\ell,&&c_{j,M+1}^\ell=1,&&\forall j,\ell.
\end{aligned}
\end{equation*}

\section{Numerical Examples}\label{sec5}
In this section, we apply our new high-resolution method to simulate the bio-convection patterns of the oxygen-driven swimming bacteria in
different sessile drops. In all of the examples, we use a uniform mesh with $\dx=\dy=0.01$ and $\dt=6.25\times10^{-6}$, which is chosen
according to \eref{4.27} and Remark \ref{rmk41}. We fix the thickness of the diffuse-domain boundary to be $\varepsilon=0.01$. We follow
\cite{tuval2005bacterial} and choose the cut-off function $r(c)$ being \eref{rc} with $c^*=0.3$, and the following parameters: $\alpha=10$,
$\delta=5$, and $\Sc=500$. The values of $\beta$ and $\gamma$ will vary and will be specified below.

\subsection{Stable Stationary Plumes}\label{sec51}
In this section, we consider four sessile drops of different shapes determined by a given function $f(x,y)$ representing the original domain
$\Omega=\{(x,y)~|~f(x,y)>0,y>0\}$, for which we compute the signed distance function $d(x,y)$ to $\Gamma=\{(x,y)~|~f(x,y)=0,y>0\}$ needed
to obtain the diffuse-domain function $\phi(x,y)$; see \eref{3.7}. In order to implement the proposed diffuse-domain based method, $\Omega$
is imbedded into a larger domain $\widetilde\Omega$, which is taken either $\widetilde\Omega=[-5,5]\times[0,1.5]$ (Examples 1 and 2) or
$\widetilde\Omega=[-7.5,7.5]\times[0,1.5]$ (Examples 3 and 4).

In Examples 1--4, we take the parameters $\beta=10$ and $\gamma=1000$.

\begin{figure}[ht!]
\centerline{\hspace*{0.1cm}\includegraphics[trim=2.4cm 0.8cm 2.3cm 0.4cm,clip,width=0.49\textwidth]{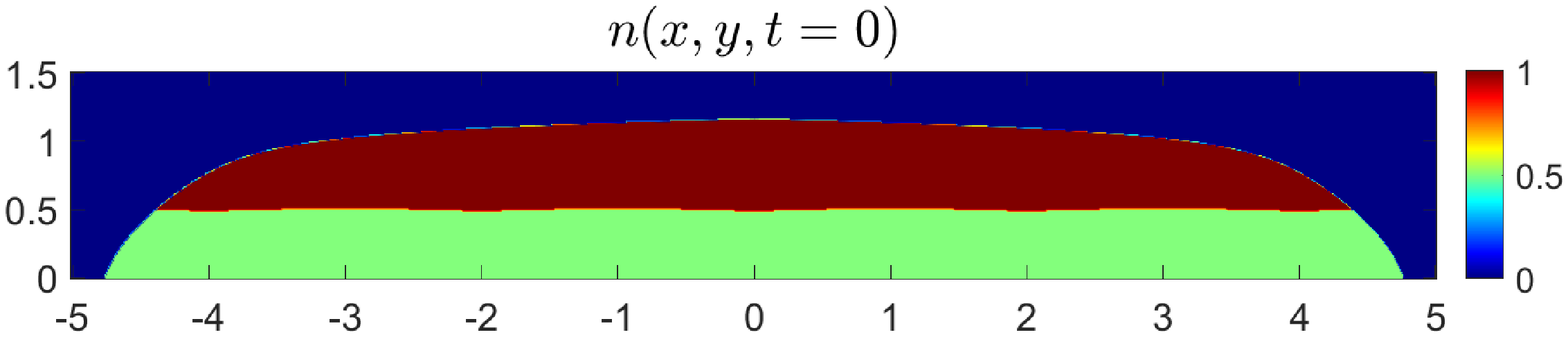}\hspace*{0.25cm}
            \includegraphics[trim=2.4cm 0.8cm 2.3cm 0.4cm,clip,width=0.49\textwidth]{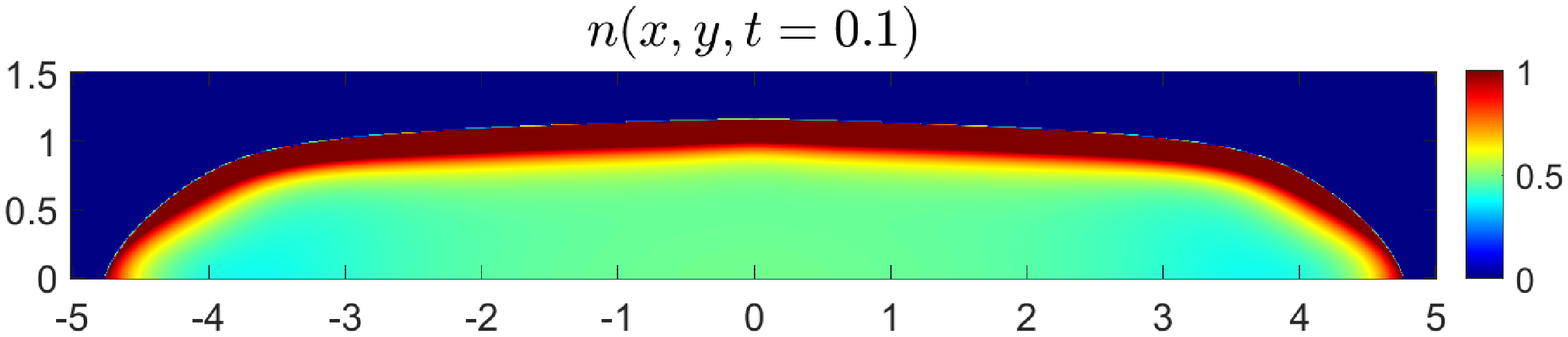}}
\vskip6pt
\centerline{\hspace*{0.1cm}\includegraphics[trim=2.4cm 0.8cm 2.3cm 0.4cm,clip,width=0.49\textwidth]{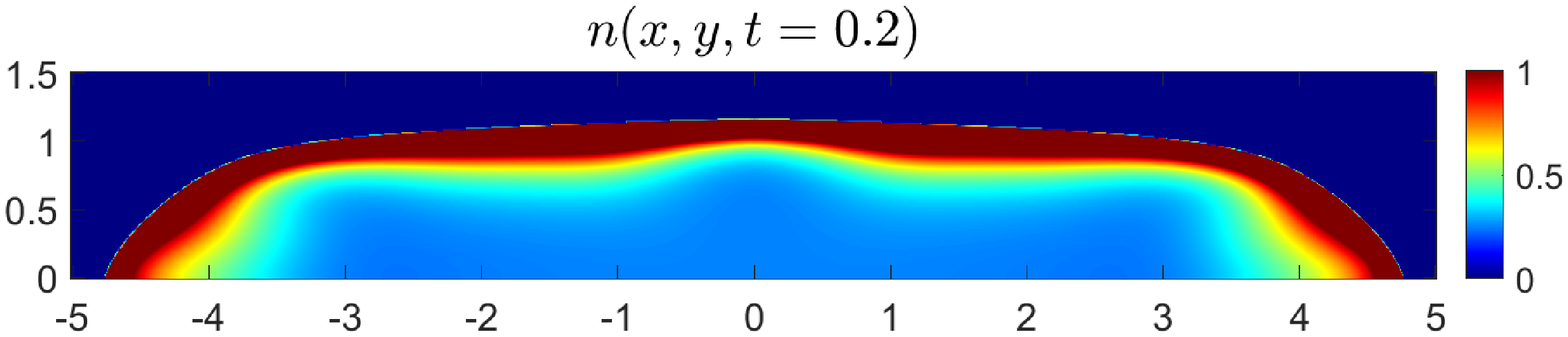}\hspace*{0.25cm}
            \includegraphics[trim=2.4cm 0.8cm 2.3cm 0.4cm,clip,width=0.49\textwidth]{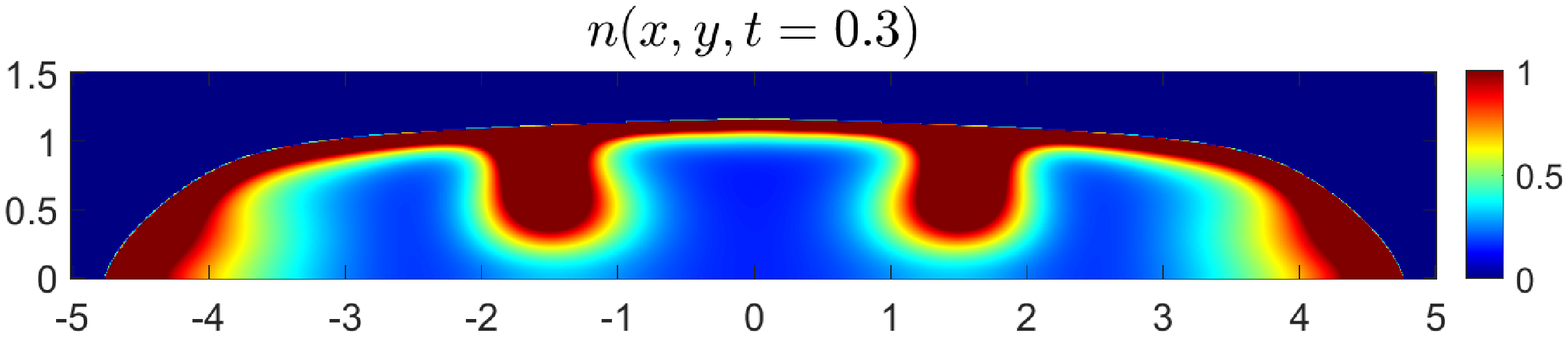}}
\vskip6pt
\centerline{\hspace*{0.1cm}\includegraphics[trim=2.4cm 0.8cm 2.3cm 0.4cm,clip,width=0.49\textwidth]{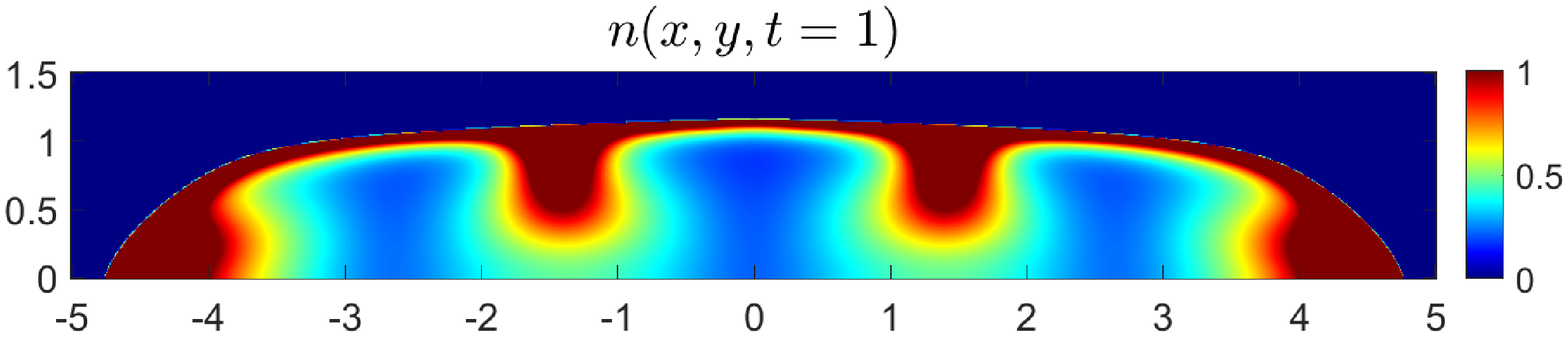}\hspace*{0.25cm}
            \includegraphics[trim=2.4cm 0.8cm 2.3cm 0.4cm,clip,width=0.49\textwidth]{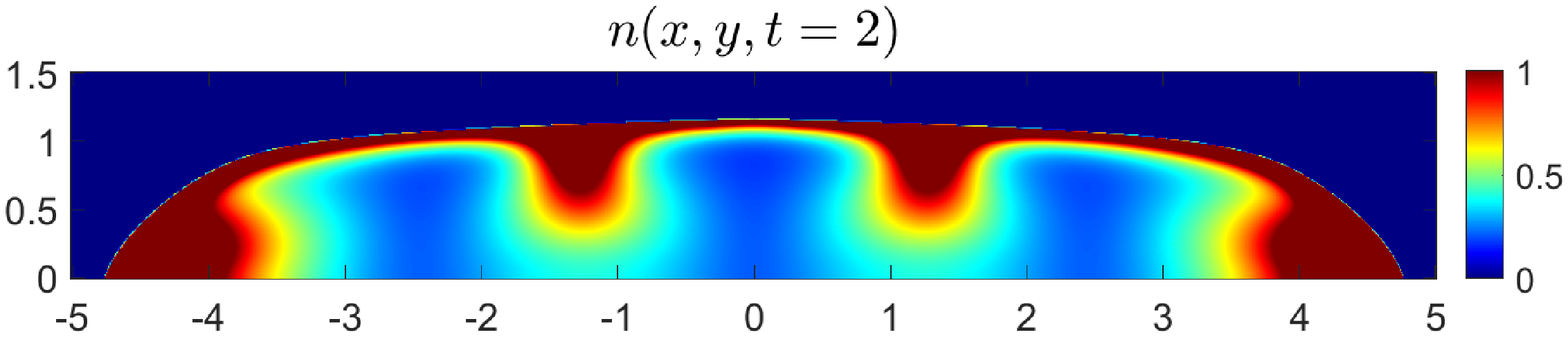}}
\vskip6pt
\centerline{\hspace*{0.1cm}\includegraphics[trim=2.4cm 0.8cm 2.3cm 0.4cm,clip,width=0.49\textwidth]{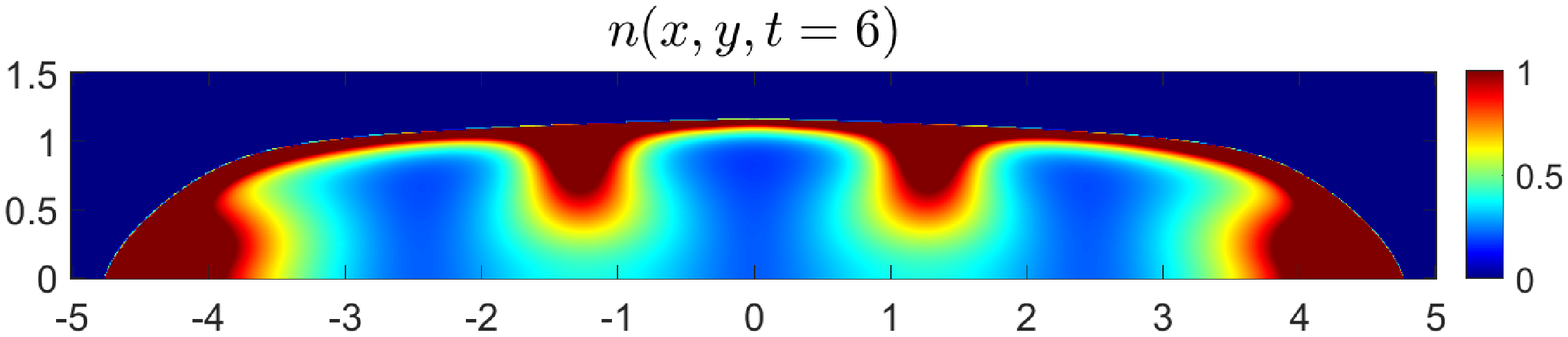}\hspace*{0.25cm}
            \includegraphics[trim=2.4cm 0.8cm 2.3cm 0.4cm,clip,width=0.49\textwidth]{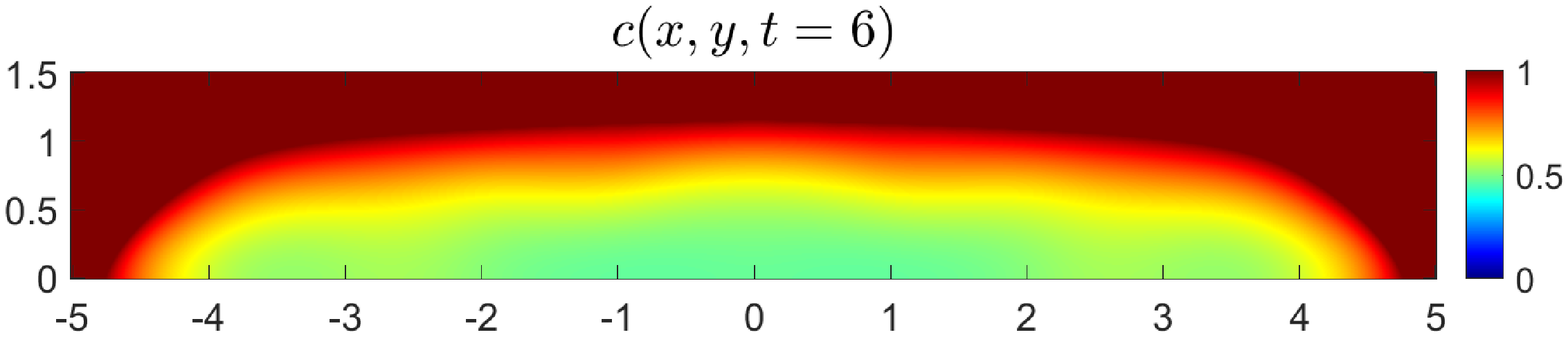}}
\caption{\sf Example 1: Time snapshots of the computed cell densities $n$ at different times and the computed oxygen concentration $c$ at
the final time.\label{fig51}}
\end{figure}

\paragraph{Example 1.} In the first example, we solve the system \eqref{3.1}--\eqref{3.8} subject to the following initial data:
\begin{equation*}
\begin{aligned}
&n(x,y,0)=\begin{cases}
1&\mbox{if }y>0.499-0.01\sin\big(\pi(x-1.5)\big),\\
0.5&\mbox{otherwise},
\end{cases}\\
&c(x,y,0)\equiv1,\quad u(x,y,0)=v(x,y,0)\equiv0,
\end{aligned}
\end{equation*}
which is prescribed in the domain $\Omega$ determined by
\begin{equation*}
f(x,y)=\begin{cases}
4.8+x-(0.9y+0.2)^2-0.1(0.9y+0.2)^{16}&\mbox{if }x\le0,\\
4.8-x-(0.9y+0.2)^2-0.1(0.9y+0.2)^{16}&\mbox{otherwise};
\end{cases}
\end{equation*}
see the upper left panel in Figure \ref{fig51}, where the shape of the drop and initial cell density are plotted.

\begin{figure}[ht!]
\centerline{\includegraphics[width=0.3\textwidth]{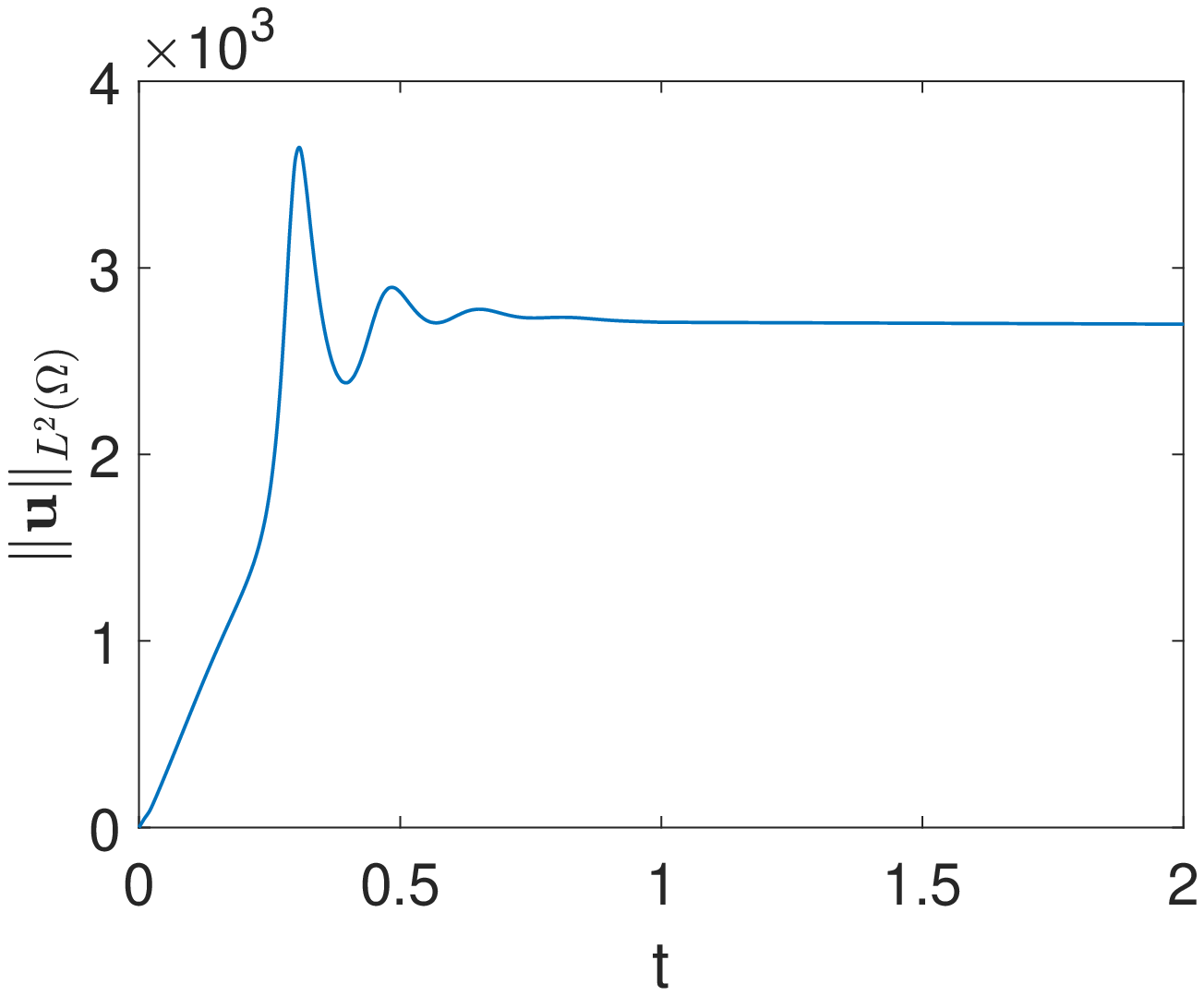}}
\caption{\sf Example 1: Time-evolution of the kinetic energy $\|\bm u\|_{L^2(\Omega)}$.\label{fig52}}
\end{figure}

\begin{figure}[ht!]
\centerline{\includegraphics[width=0.35\textwidth]{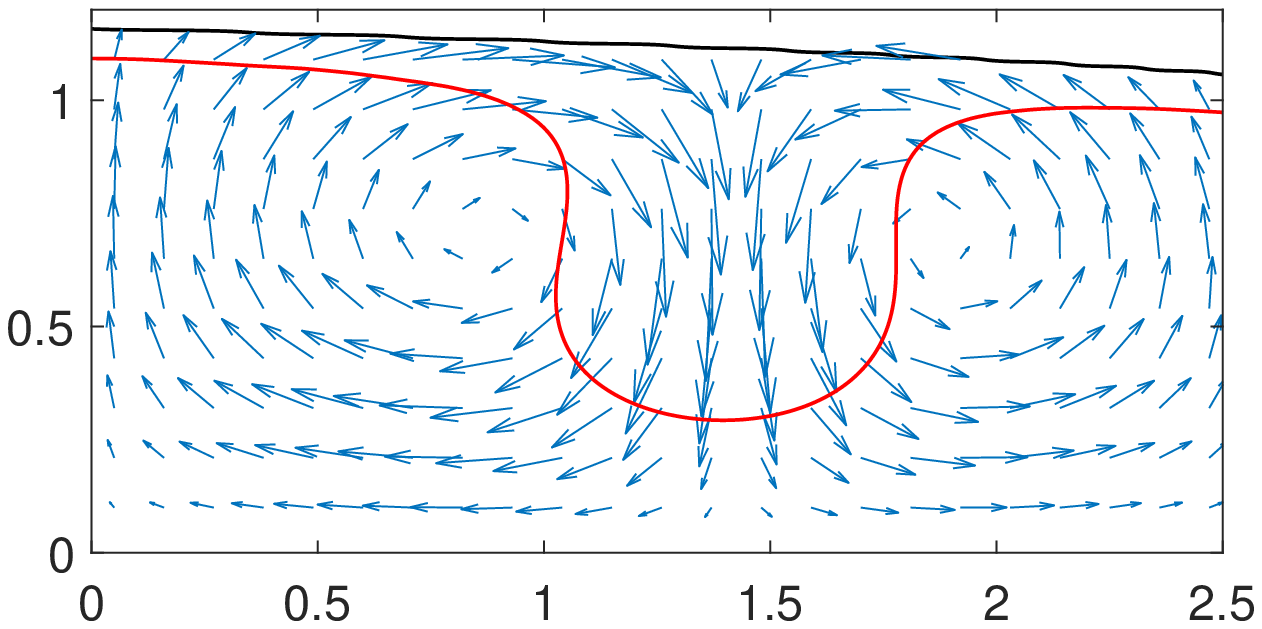}\hspace*{0.5cm}
            \includegraphics[width=0.35\textwidth]{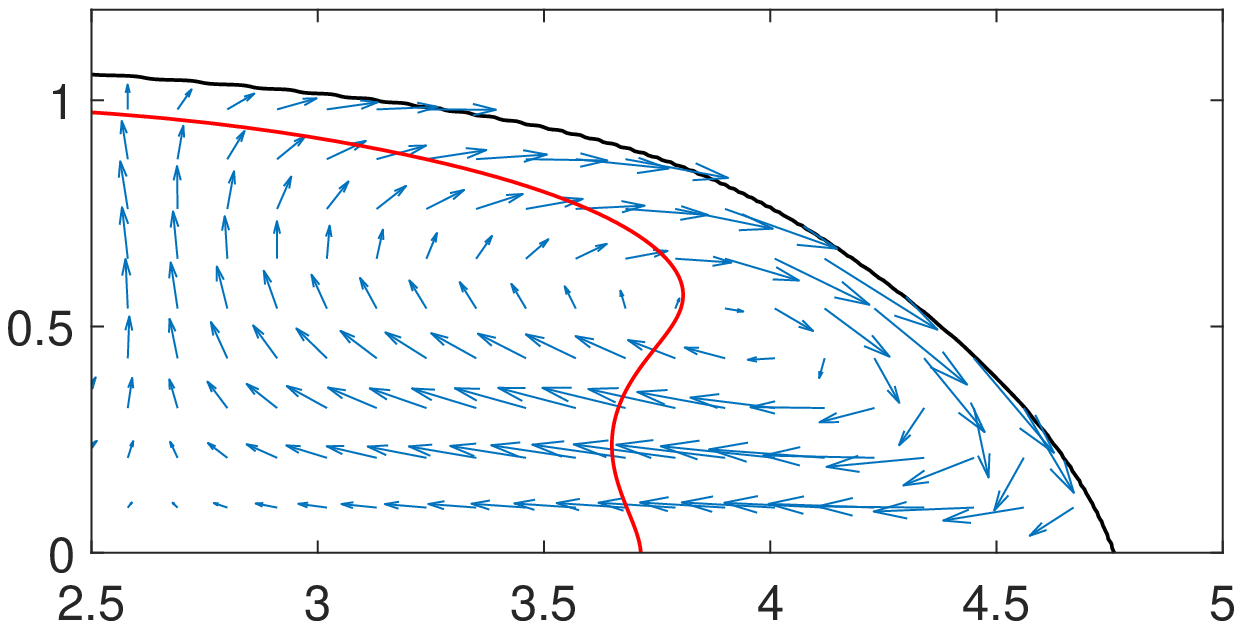}}
\caption{\sf Example 1: Velocity field $\bm u$ and the $n=0.7$ level set at time $t=1$.\label{fig53}}
\end{figure}

The time evolution of the cell density $n$ is also shown in Figure \ref{fig51}. As one can see, the bacteria first ($t=0.1$) aggregate along
the boundary $\Gamma$ as the concentration of the oxygen is high there, but then the gravity forces start dominating ($t=0.2$) and some of
the bacteria fall down forming the plumes ($t=0.3$) both at the corners and in the middle of the drop. Later on, the shape of the plumes
slightly changes ($t=1$) and by time $t=2$ the plumes are already stationary (compare with the cell density at the very large time $t=6$).
For the sake of brevity, we plot the $c$-component of the computed solution only at time $t=6$. The obtained results confirm the ability of
the proposed numerical method to capture stationary plumes in a stable manner. In order to numerically verify the stability of the plumes,
we plot the time evolution of the kinetic energy (Figure \ref{fig52}), which clearly converges to a constant value, and the velocity field
together with the $n=0.7$ cell density level set at time $t=1$ (Figure \ref{fig53}), which illustrate how the stationary plumes are
supported by the fluid.

We note that the obtained results are in a very good qualitative agreement with the results reported in \cite{chertock2012sinking}, where
the system \eref{1.1} was considered in a rectangular domain subject to the periodic boundary conditions in the horizontal direction. In
addition, the diffuse-domain based numerical method proposed here is capable of treating non-rectangular domains and resolving the
accumulation layers at the drop corners and creation of vortices there.

\begin{figure}[ht!]
\centerline{\hspace*{0.1cm}\includegraphics[trim=2.4cm 0.8cm 2.3cm 0.4cm,clip,width=0.49\textwidth]{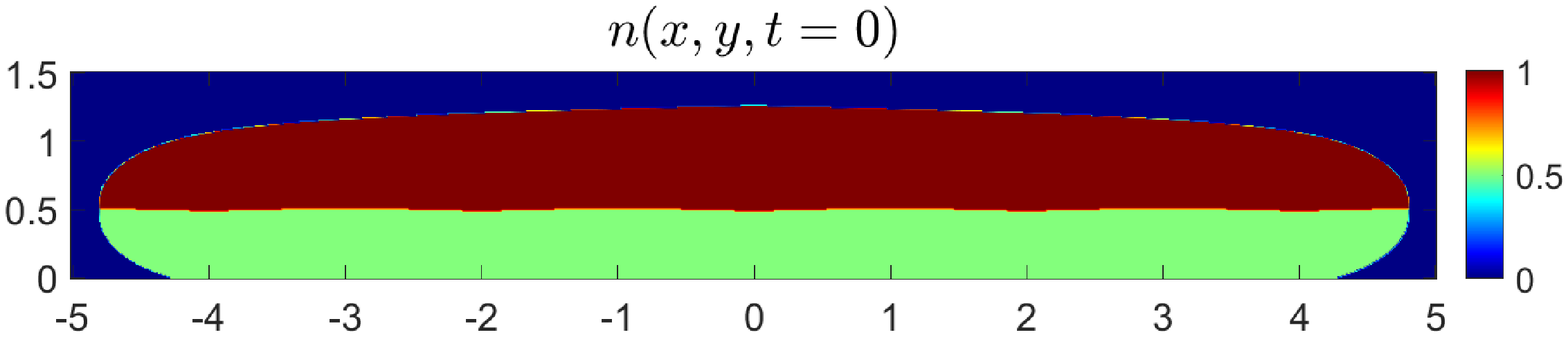}\hspace*{0.25cm}
            \includegraphics[trim=2.4cm 0.8cm 2.3cm 0.4cm,clip,width=0.49\textwidth]{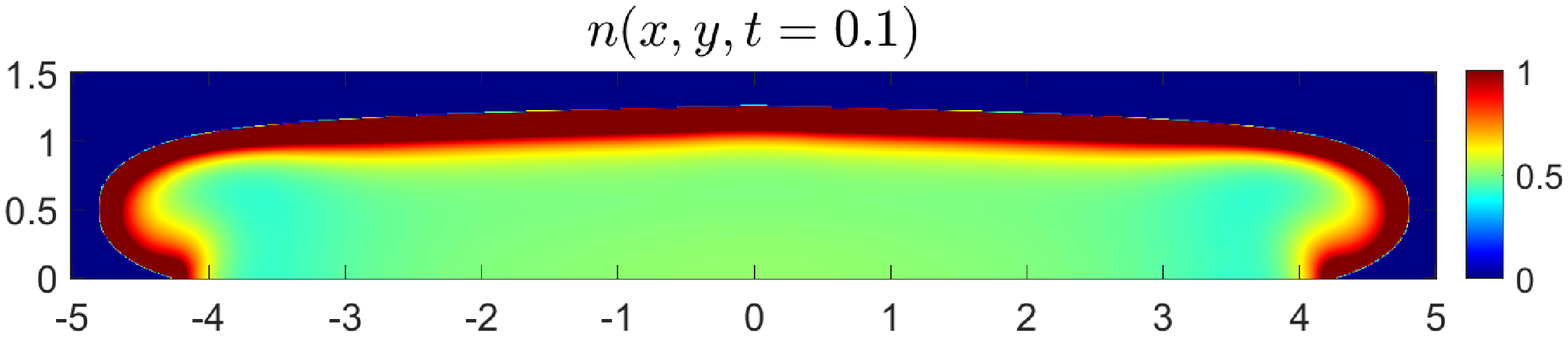}}
\vskip6pt
\centerline{\hspace*{0.1cm}\includegraphics[trim=2.4cm 0.8cm 2.3cm 0.4cm,clip,width=0.49\textwidth]{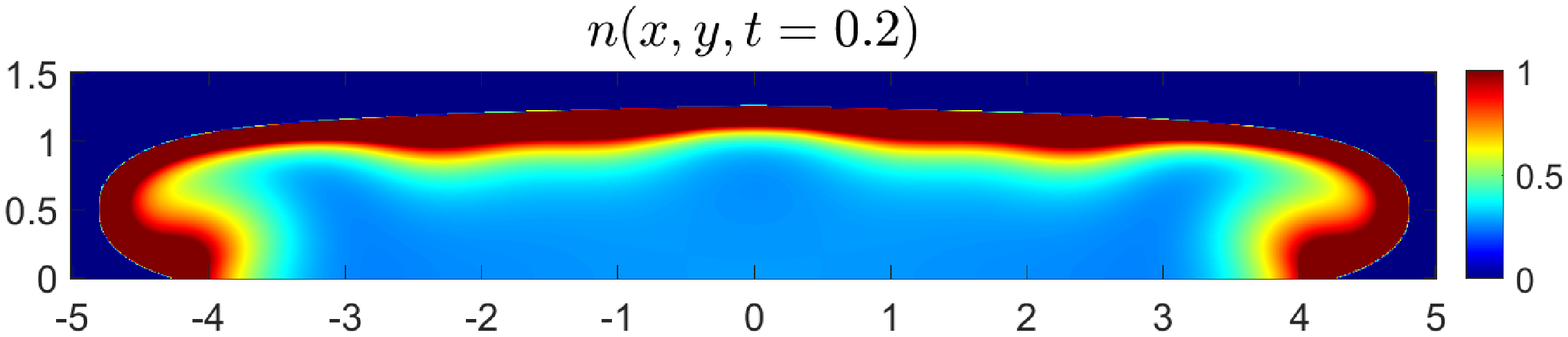}\hspace*{0.25cm}
            \includegraphics[trim=2.4cm 0.8cm 2.3cm 0.4cm,clip,width=0.49\textwidth]{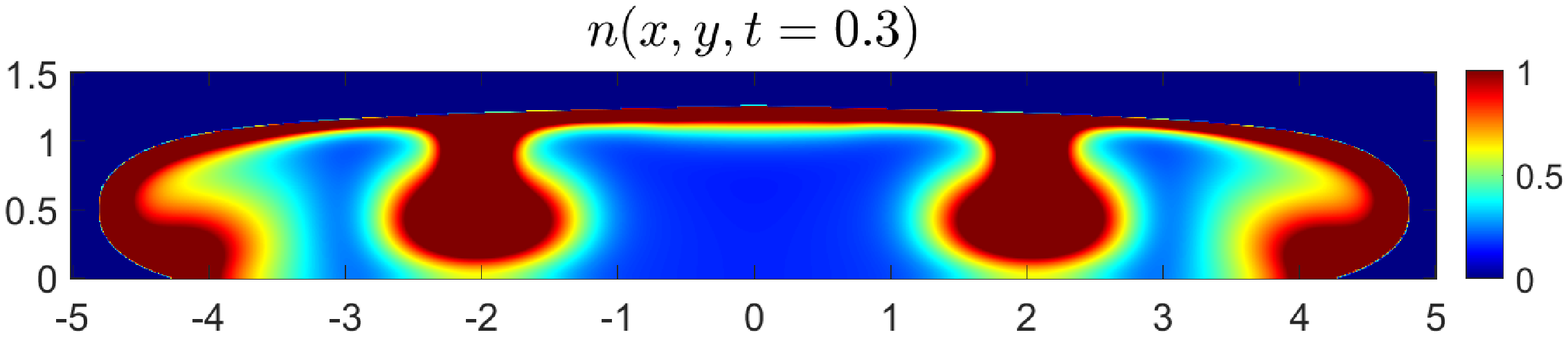}}
\vskip6pt
\centerline{\hspace*{0.1cm}\includegraphics[trim=2.4cm 0.8cm 2.3cm 0.4cm,clip,width=0.49\textwidth]{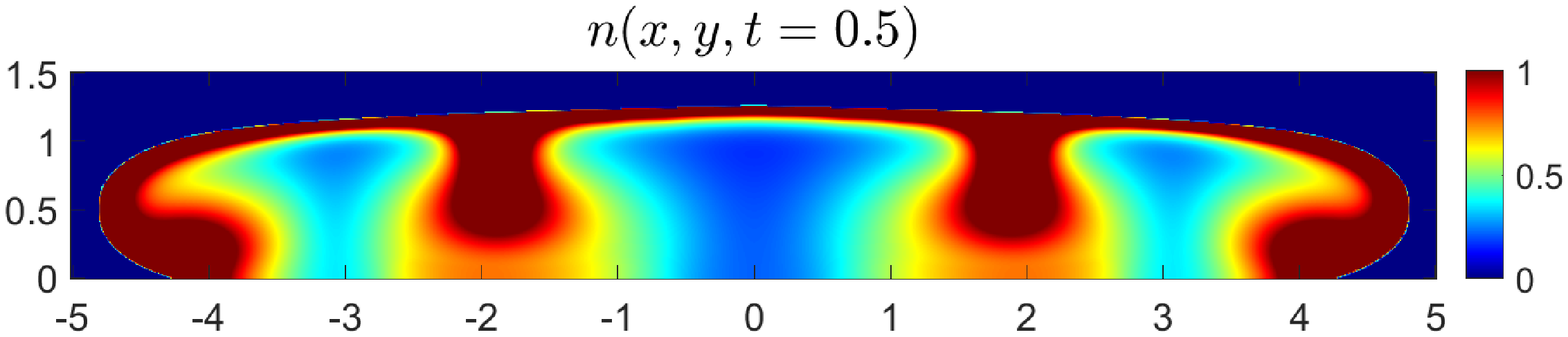}\hspace*{0.25cm}
            \includegraphics[trim=2.4cm 0.8cm 2.3cm 0.4cm,clip,width=0.49\textwidth]{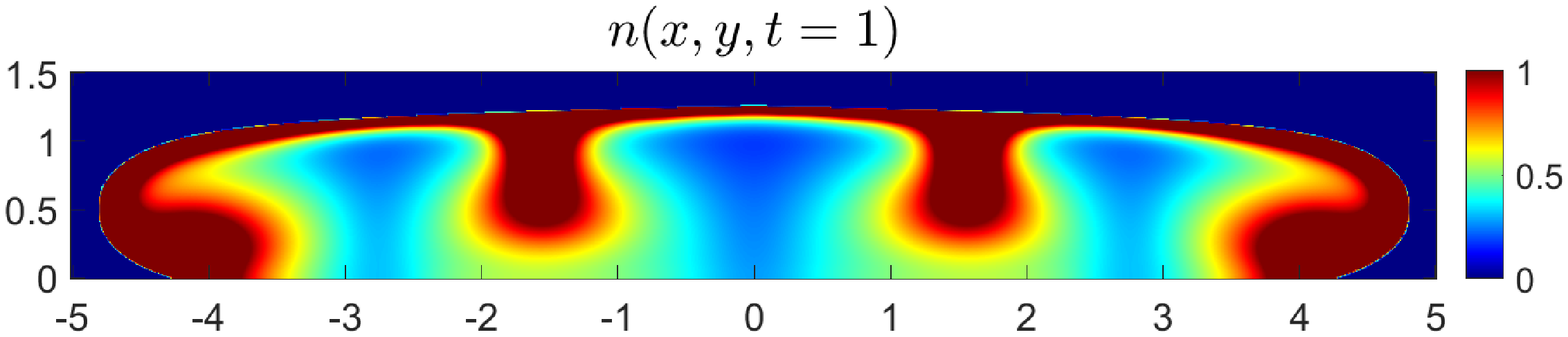}}
\vskip6pt
\centerline{\hspace*{0.1cm}\includegraphics[trim=2.4cm 0.8cm 2.3cm 0.4cm,clip,width=0.49\textwidth]{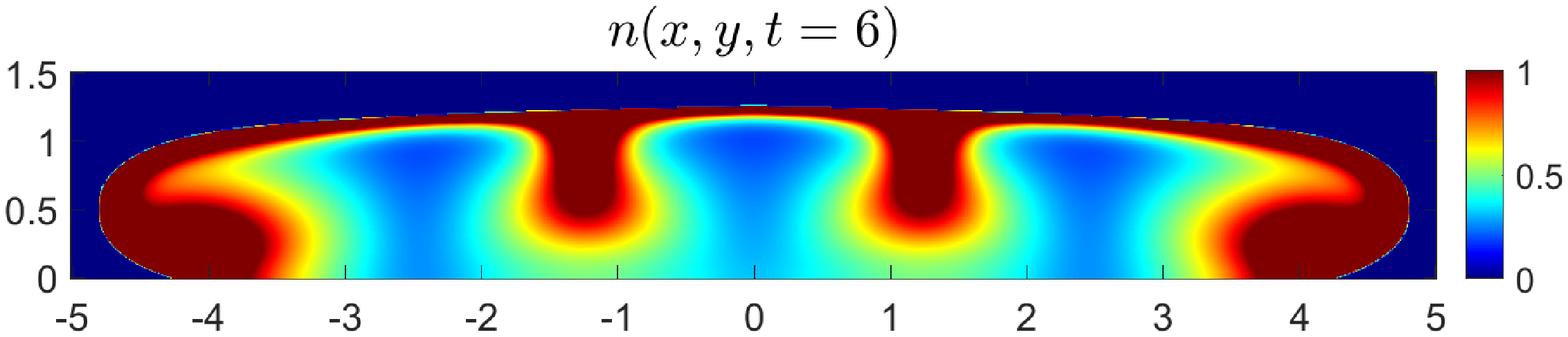}\hspace*{0.25cm}
            \includegraphics[trim=2.4cm 0.8cm 2.3cm 0.4cm,clip,width=0.49\textwidth]{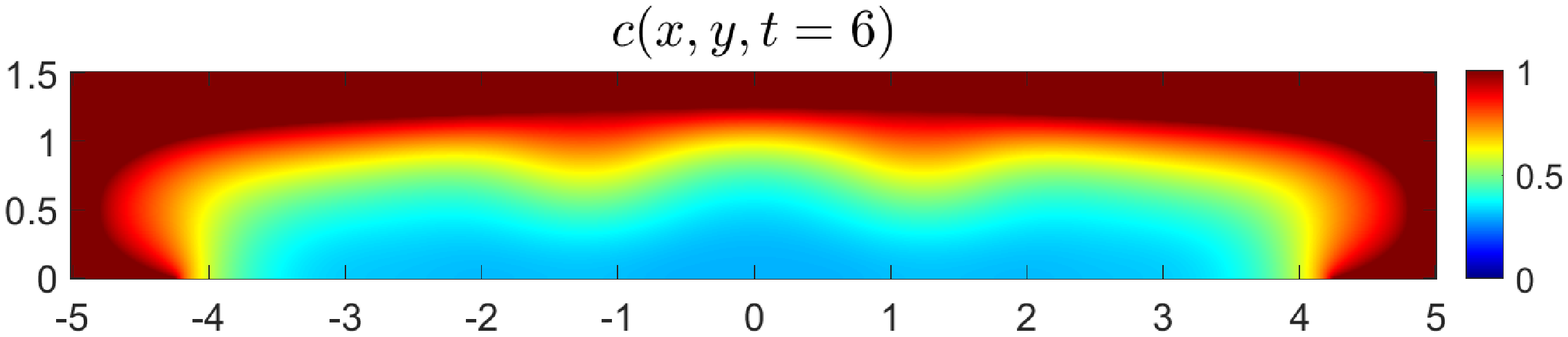}}
\caption{\sf Example 2: Time snapshots of the computed cell densities $n$ at different times and the computed oxygen concentration $c$ at
the final time.\label{fig54}}
\end{figure}

\paragraph{Example 2.} Next, we consider the same initial setting  as in Example 1 but in a sessile drop of a different shape determined by
\begin{equation*}
f(x,y)=\begin{cases}
4.8+x-|1.5y-0.75|^{2.5}-(1.5y-0.75)^{10}&\mbox{if }x\le0,\\
4.8-x-|1.5y-0.75|^{2.5}-(1.5y-0.75)^{10}&\mbox{otherwise}.
\end{cases}
\end{equation*}
Compared with the drop in the previous example, this one has rounded edges while still having a flat bottom interface; see the upper left
panel in Figure \ref{fig54}, where the shape of the drop and initial cell density are plotted.

\begin{figure}[ht!]
\centerline{\includegraphics[width=0.3\textwidth]{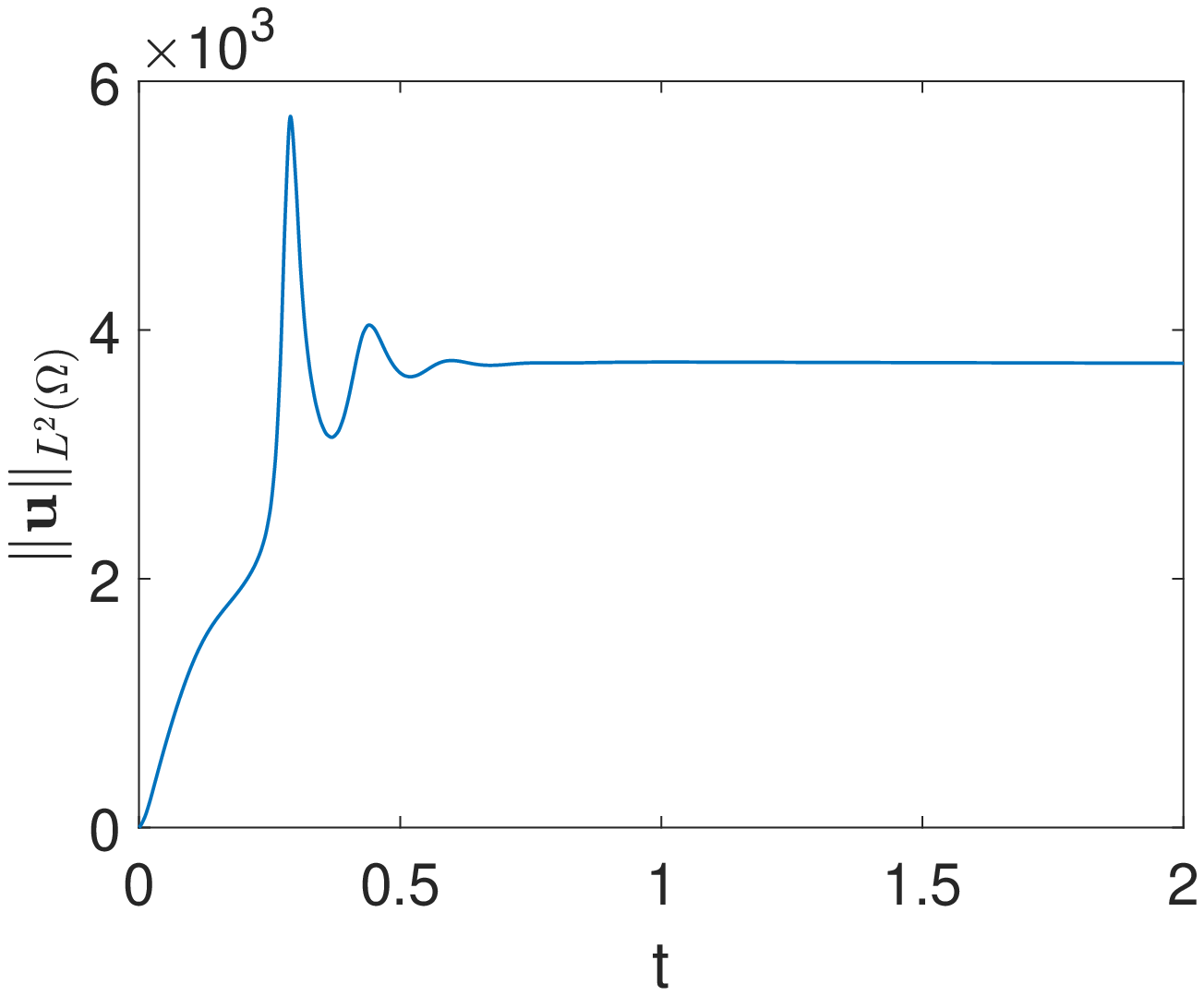}}
\caption{\sf Example 2: Time-evolution of the kinetic energy $\|\bm u\|_{L^2(\Omega)}$.\label{fig55}}
\end{figure}

\begin{figure}[ht!]
\centerline{\includegraphics[width=0.35\textwidth]{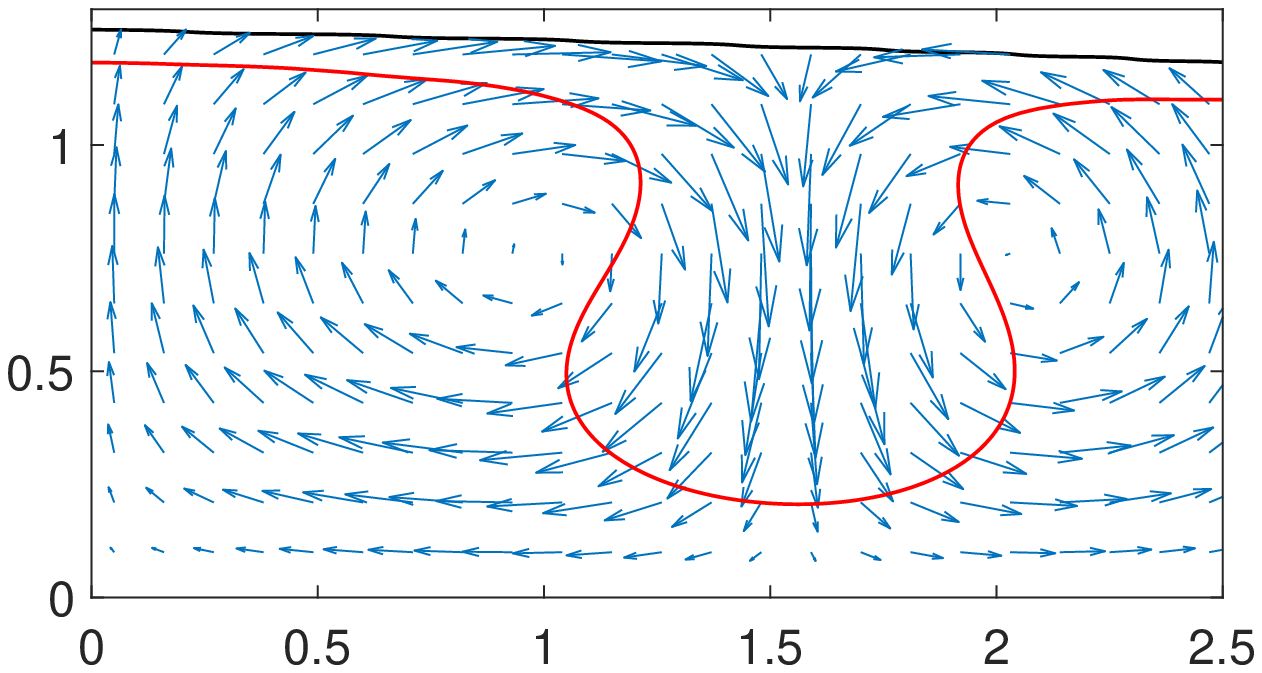}\hspace*{0.5cm}
            \includegraphics[width=0.35\textwidth]{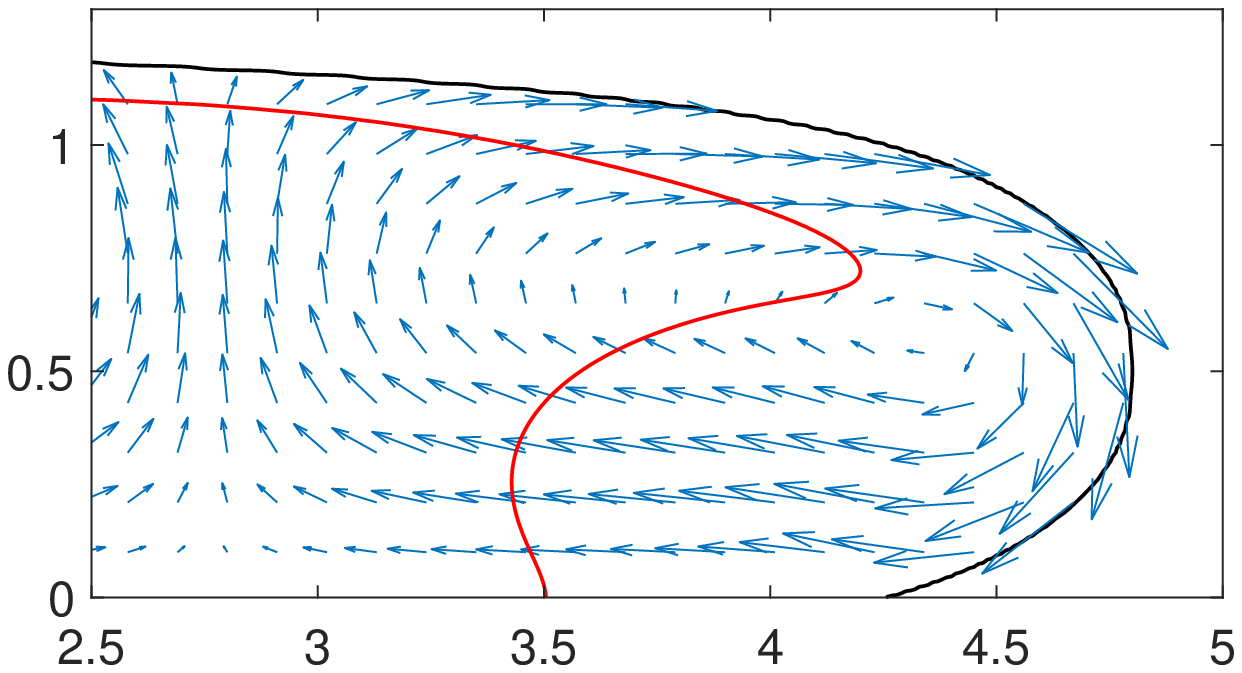}}
\caption{\sf Example 2: Velocity field $\bm u$ and the $n=0.7$ level set at time $t=1$.\label{fig56}}
\end{figure}

The time evolution of the cell density $n$ is also shown in Figure \ref{fig54}. As one can see, the evolution process is similar to the one
in Example 1 with the only exception that the structure of the aggregated area at the edges of the drop is different as the oxygen supply is
available underneath that part of the drop considered here. The solution converges to a stationary state containing bacteria plumes, which
can be seen in the bottom row of Figure \ref{fig54} (the oxygen concentration $c$ at the final time $t=6$ is also shown there). The
convergence towards the steady state is also confirmed by following the time evolution of the kinetic energy (Figure \ref{fig55}), which
clearly flattens by time $t=1$. As in Example 1, we also plot the the velocity field together with the $n=0.7$ cell density level set at
time $t=1$ (Figure \ref{fig56}), which illustrate how the stationary plumes are supported by the fluid.

Once again, we emphasize that the proposed diffuse-domain based numerical method is capable of numerically solving the fluid-chemotaxis
system in rather complicated domains.

\begin{figure}[ht!]
\centerline{\hspace*{0.1cm}\includegraphics[trim=3.6cm 0.9cm 3.2cm 0.5cm,clip,width=8.1cm,height=1.75cm]{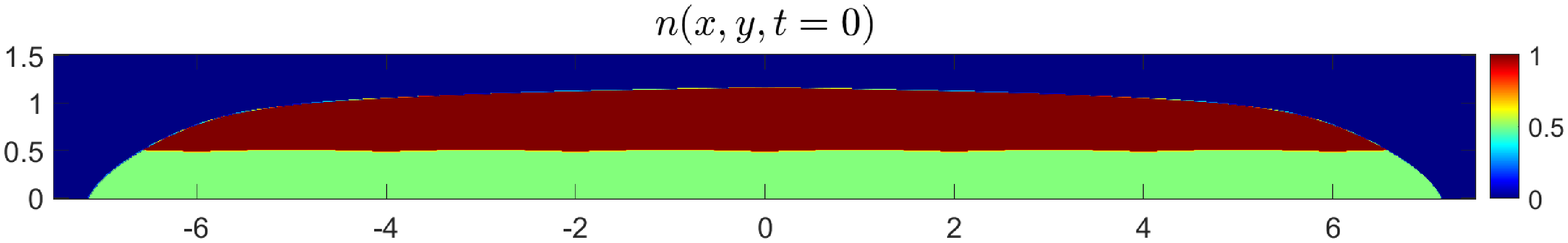}\hspace*{0.25cm}
           \includegraphics[trim=3.6cm 0.9cm 3.2cm 0.5cm,clip,width=8.1cm,height=1.75cm]{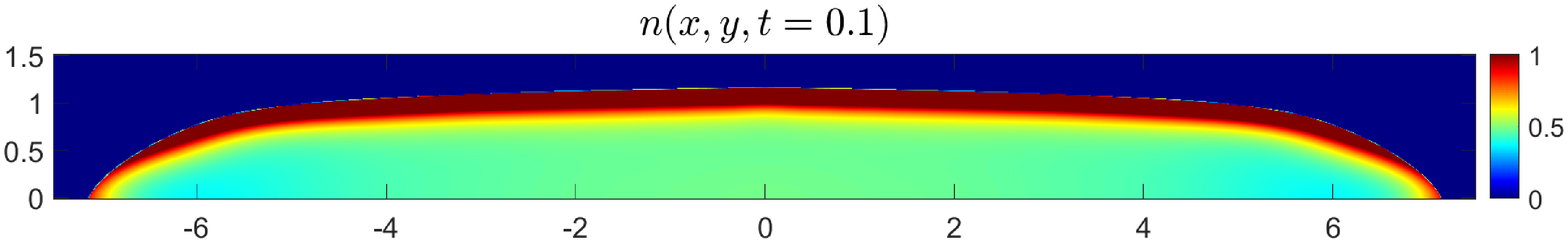}}
\vskip6pt
\centerline{\hspace*{0.1cm}\includegraphics[trim=3.6cm 0.9cm 3.2cm 0.5cm,clip,width=8.1cm,height=1.75cm]{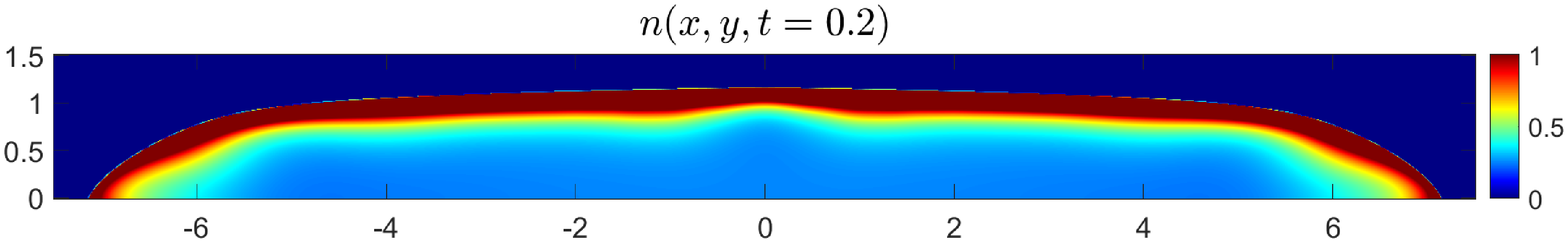}\hspace*{0.25cm}
            \includegraphics[trim=3.6cm 0.9cm 3.2cm 0.5cm,clip,width=8.1cm,height=1.75cm]{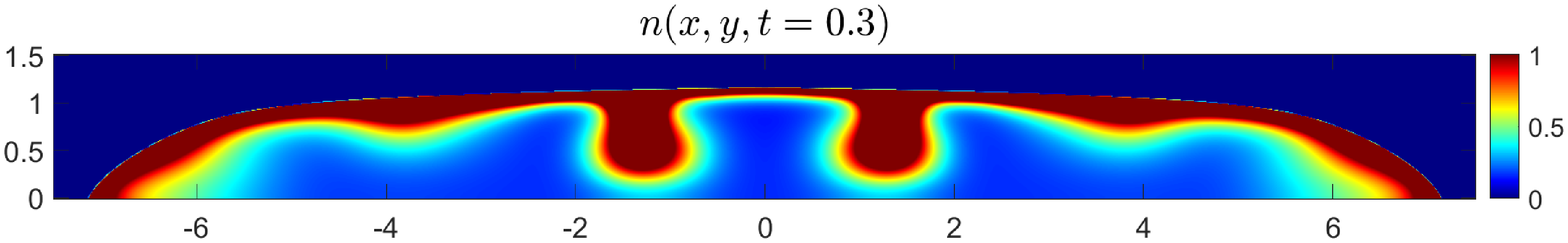}}
\vskip6pt
\centerline{\hspace*{0.1cm}\includegraphics[trim=3.6cm 0.9cm 3.2cm 0.5cm,clip,width=8.1cm,height=1.75cm]{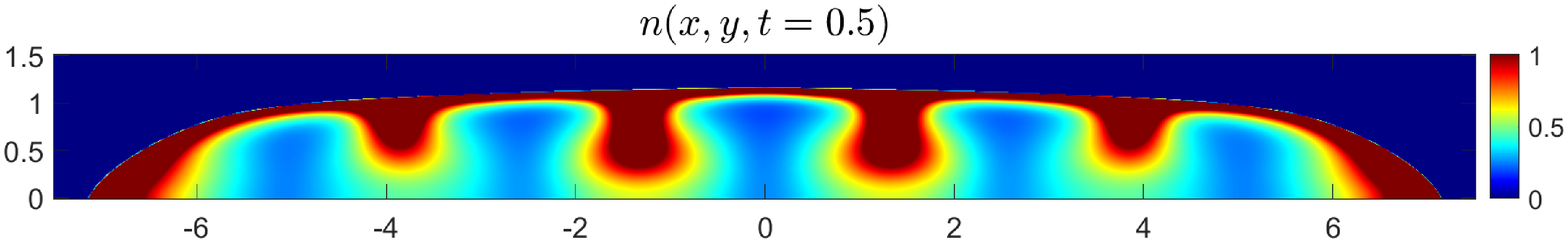}\hspace*{0.25cm}
            \includegraphics[trim=3.6cm 0.9cm 3.2cm 0.5cm,clip,width=8.1cm,height=1.75cm]{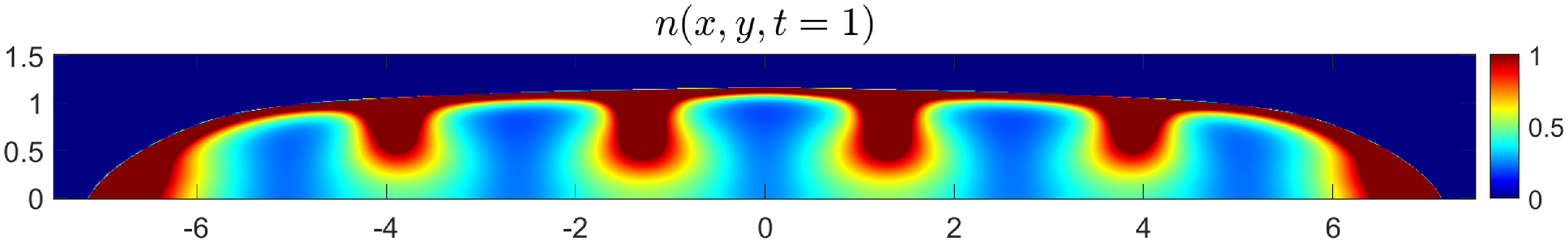}}
\vskip6pt
\centerline{\hspace*{0.1cm}\includegraphics[trim=3.6cm 0.9cm 3.2cm 0.5cm,clip,width=8.1cm,height=1.75cm]{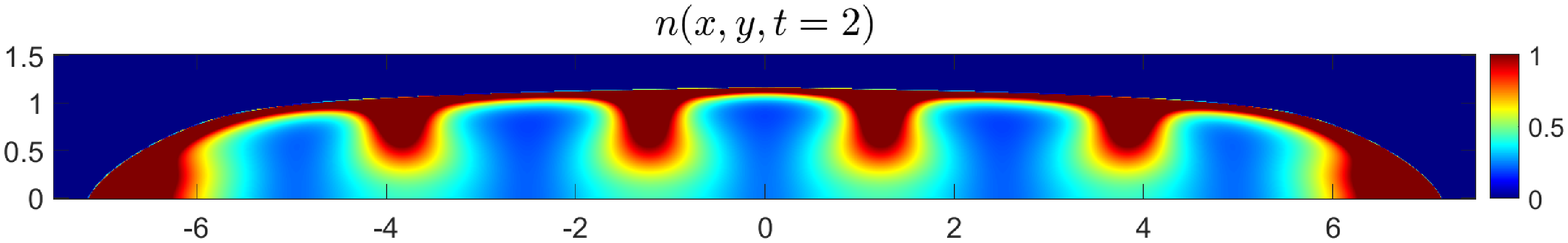}\hspace*{0.25cm}
            \includegraphics[trim=3.6cm 0.9cm 3.2cm 0.5cm,clip,width=8.1cm,height=1.75cm]{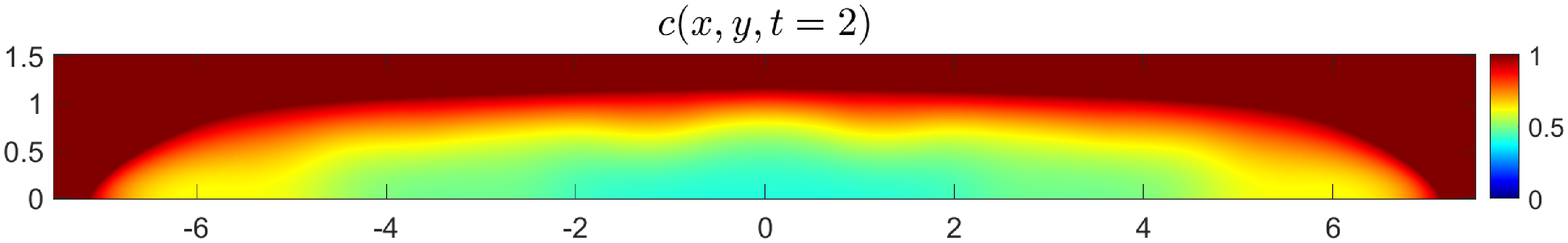}}
\caption{\sf Example 3: Time snapshots of the computed cell densities $n$ at different times and the computed oxygen concentration $c$ at
the final time.\label{fig57}}
\end{figure}

\paragraph{Example 3.} In the third example, we consider the same initial setting as in Example 1 but the drop is now longer. Its precise
shape is determined by
\begin{equation*}
f(x,y)=\begin{cases}
4.8+\frac{2}{3}x-(0.9y+0.2)^2-0.1(0.9y+0.2)^{16}&\mbox{if }x\le0,\\
4.8-\frac{2}{3}x-(0.9y+0.2)^2-0.1(0.9y+0.2)^{16}&\mbox{otherwise};
\end{cases}
\end{equation*}
see the upper left panel in Figure \ref{fig57}, where the shape of the drop and initial cell density are plotted. The time evolution of the
cell density $n$ as well as the profile of the oxygen concentration $c$ at the final time $t=2$ by which the solution reaches its steady
state, are also shown in Figure \ref{fig57}. As one can see, the proposed numerical method can handle longer drops and the only qualitative
difference between the steady states here and in Example 1 is in the number of plumes emerging during the evolution process.

\begin{figure}[ht!]
\centerline{\hspace*{0.1cm}\includegraphics[trim=3.6cm 0.9cm 3.2cm 0.5cm,clip,width=8.1cm,height=1.75cm]{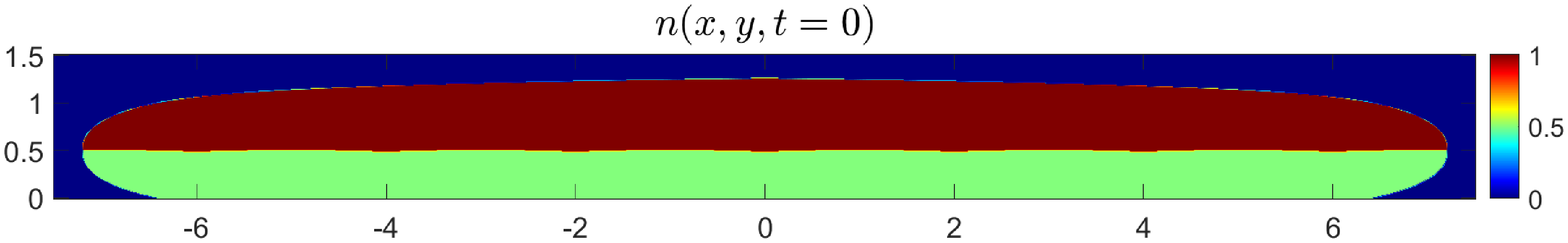}\hspace*{0.25cm}
           \includegraphics[trim=3.6cm 0.9cm 3.2cm 0.5cm,clip,width=8.1cm,height=1.75cm]{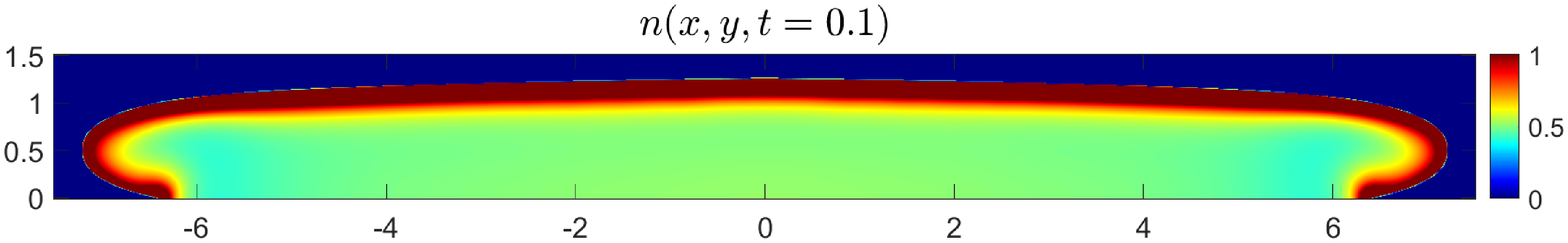}}
\vskip6pt
\centerline{\hspace*{0.1cm}\includegraphics[trim=3.6cm 0.9cm 3.2cm 0.5cm,clip,width=8.1cm,height=1.75cm]{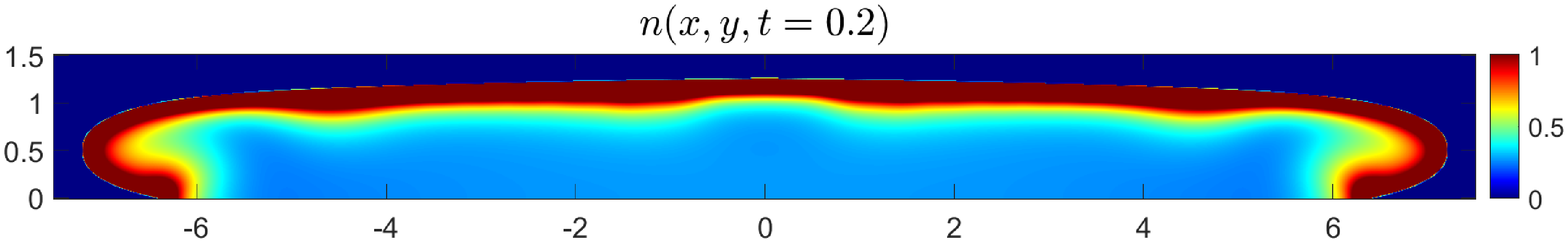}\hspace*{0.25cm}
            \includegraphics[trim=3.6cm 0.9cm 3.2cm 0.5cm,clip,width=8.1cm,height=1.75cm]{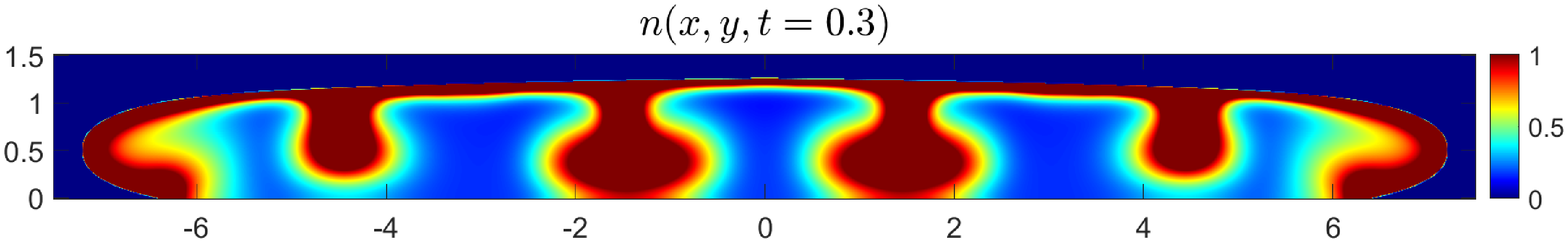}}
\vskip6pt
\centerline{\hspace*{0.1cm}\includegraphics[trim=3.6cm 0.9cm 3.2cm 0.5cm,clip,width=8.1cm,height=1.75cm]{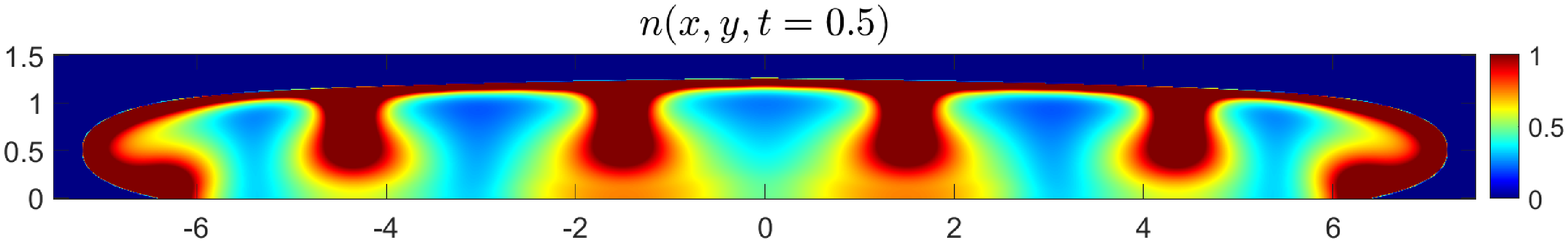}\hspace*{0.25cm}
            \includegraphics[trim=3.6cm 0.9cm 3.2cm 0.5cm,clip,width=8.1cm,height=1.75cm]{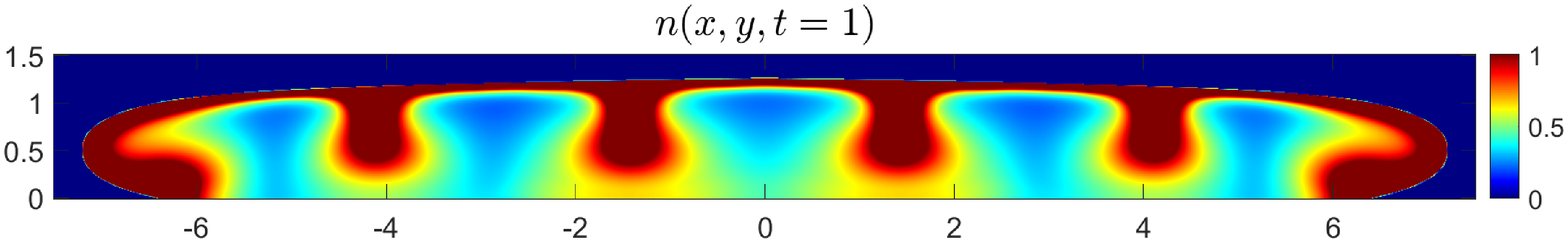}}
\vskip6pt
\centerline{\hspace*{0.1cm}\includegraphics[trim=3.6cm 0.9cm 3.2cm 0.5cm,clip,width=8.1cm,height=1.75cm]{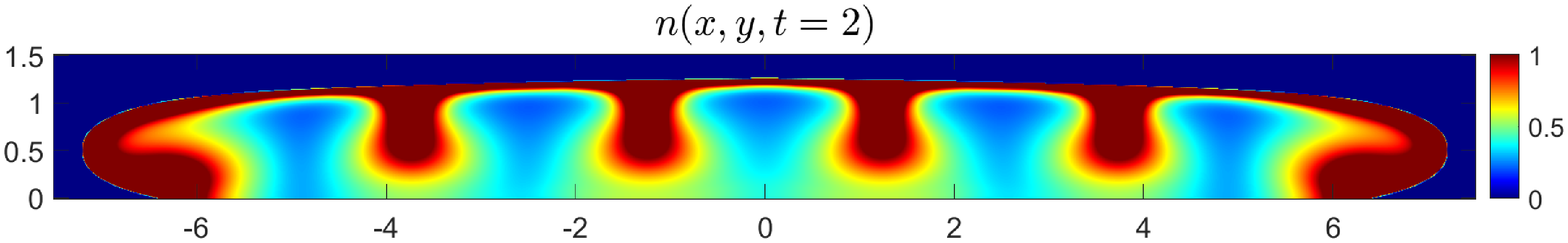}\hspace*{0.25cm}
            \includegraphics[trim=3.6cm 0.9cm 3.2cm 0.5cm,clip,width=8.1cm,height=1.75cm]{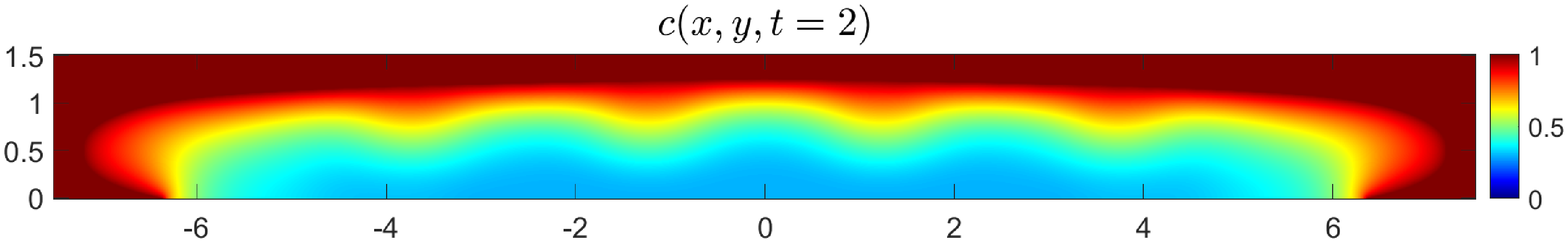}}
\caption{\sf Example 4: Time snapshots of the computed cell densities $n$ at different times and the computed oxygen concentration $c$ at
the final time.\label{fig58}}
\end{figure}

\paragraph{Example 4.} The final example of this section is a modification of Example 2 as we now take a longer drop determined by
\begin{equation*}
f(x,y)=\begin{cases}
4.8+\frac{2}{3}x-|1.5y-0.75|^{2.5}-(1.5y-0.75)^{10}&\mbox{if }x\le0,\\
4.8-\frac{2}{3}x-|1.5y-0.75|^{2.5}-(1.5y-0.75)^{10}&\mbox{otherwise};
\end{cases}
\end{equation*}
see the upper left panel in Figure \ref{fig58}, where the shape of the drop and initial cell density are plotted. The time evolution of $n$,
which converges to the steady state by $t=2$ together with the profile of $c$ at the final time can be also seen in Figure \ref{fig58}.
The obtained stationary solution contains two additional plumes compared with the solution reported in Example 2, but rather than this these
two solutions are qualitatively similar, which confirms the robustness of our numerical method.

\subsection{Mushroom-Shaped Plumes for High-Density Data}\label{sec52}
In this section, we choose the parameters $\beta=100$ and $\gamma=10000$, which correspond to a 10-times larger reference cell density
$n_r$; see \eref{2.2}. The evolution process will now be substantially faster so that we conduct the simulations for a shorter period and
take the final time $t=0.5$.

The goal of the simulations reported in Examples 5 and 6 below is to demonstrate the ability of the proposed diffuse-domain based numerical
method to handle more complicated bacteria propagation dynamics, which are expected to occur when the reference cell density $n_r$ is
larger.

\begin{figure}[ht!]
\centerline{\hspace*{0.1cm}\includegraphics[trim=2.4cm 0.8cm 2.3cm 0.4cm,clip,width=0.49\textwidth]{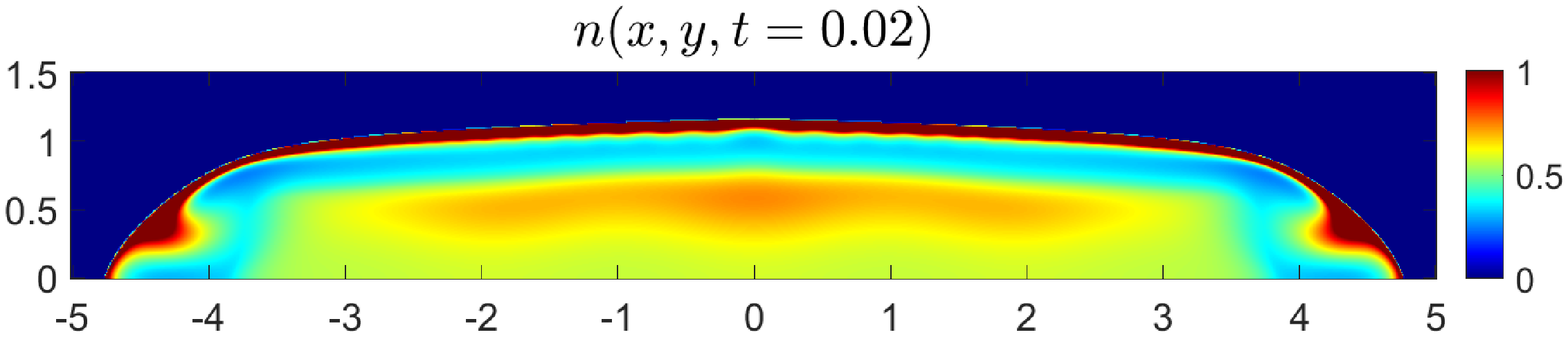}\hspace*{0.25cm}
            \includegraphics[trim=2.4cm 0.8cm 2.3cm 0.4cm,clip,width=0.49\textwidth]{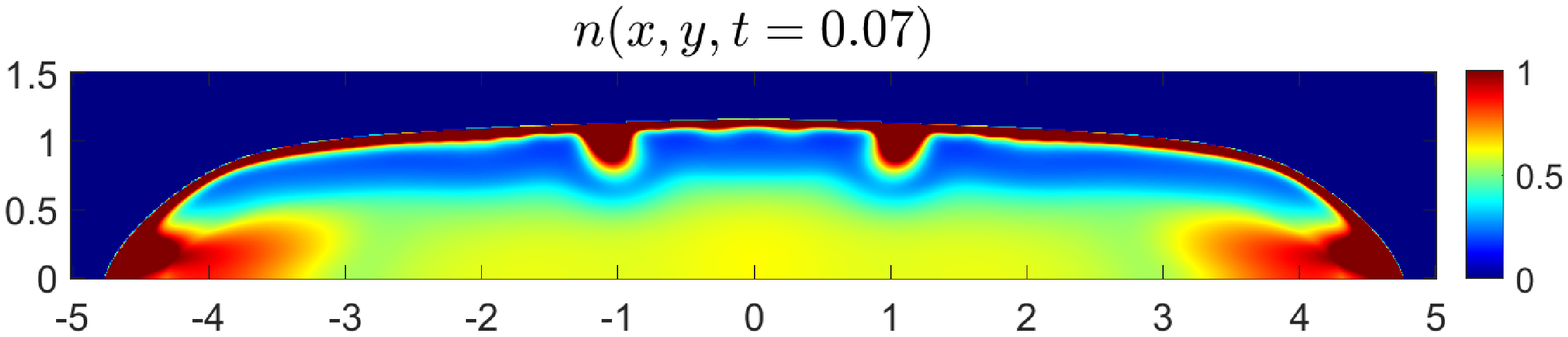}}
\vskip6pt
\centerline{\hspace*{0.1cm}\includegraphics[trim=2.4cm 0.8cm 2.3cm 0.4cm,clip,width=0.49\textwidth]{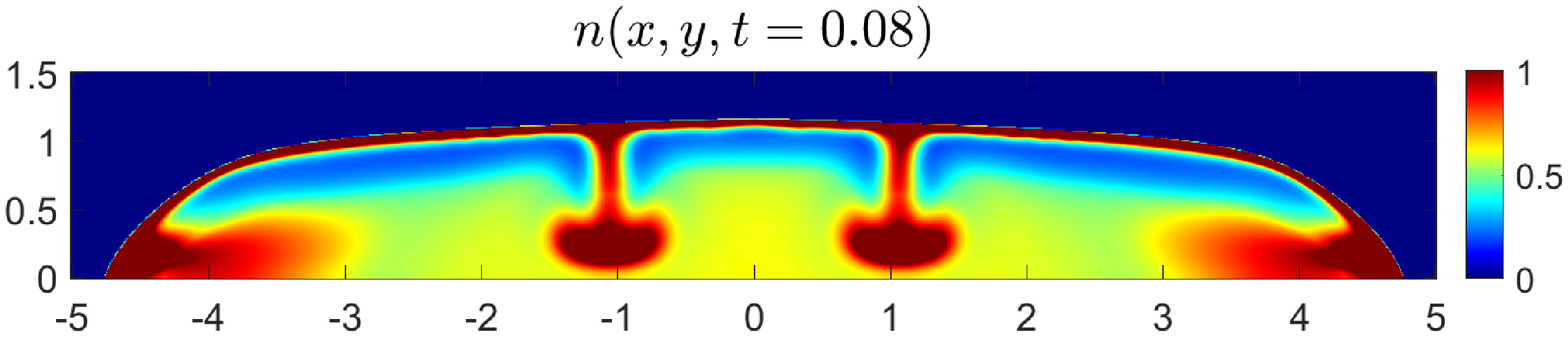}\hspace*{0.25cm}
            \includegraphics[trim=2.4cm 0.8cm 2.3cm 0.4cm,clip,width=0.49\textwidth]{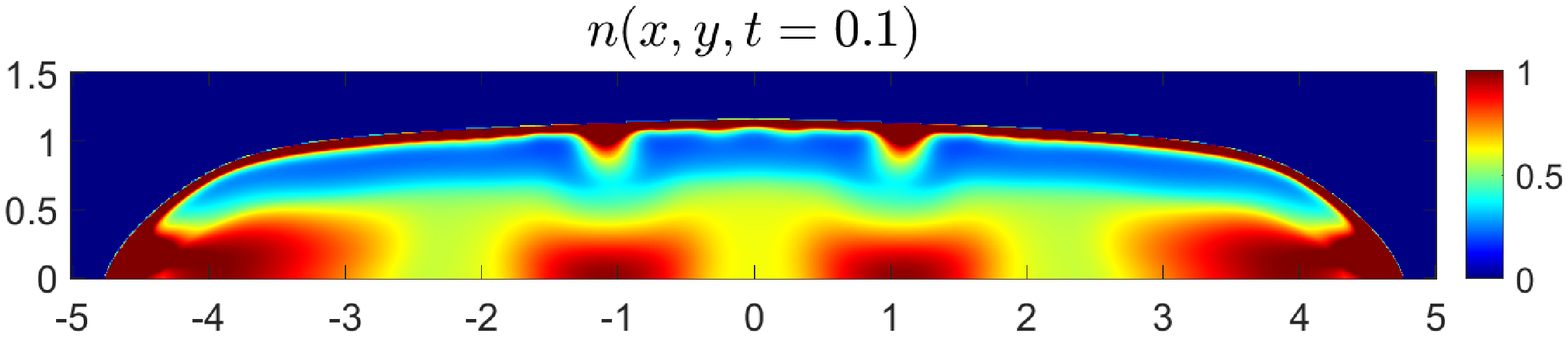}}
\vskip6pt
\centerline{\hspace*{0.1cm}\includegraphics[trim=2.4cm 0.8cm 2.3cm 0.4cm,clip,width=0.49\textwidth]{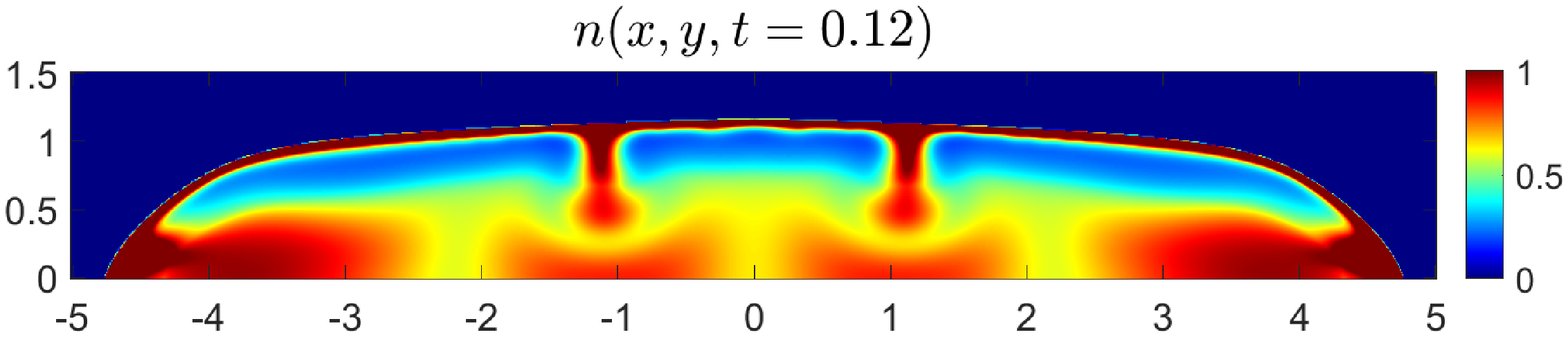}\hspace*{0.25cm}
            \includegraphics[trim=2.4cm 0.8cm 2.3cm 0.4cm,clip,width=0.49\textwidth]{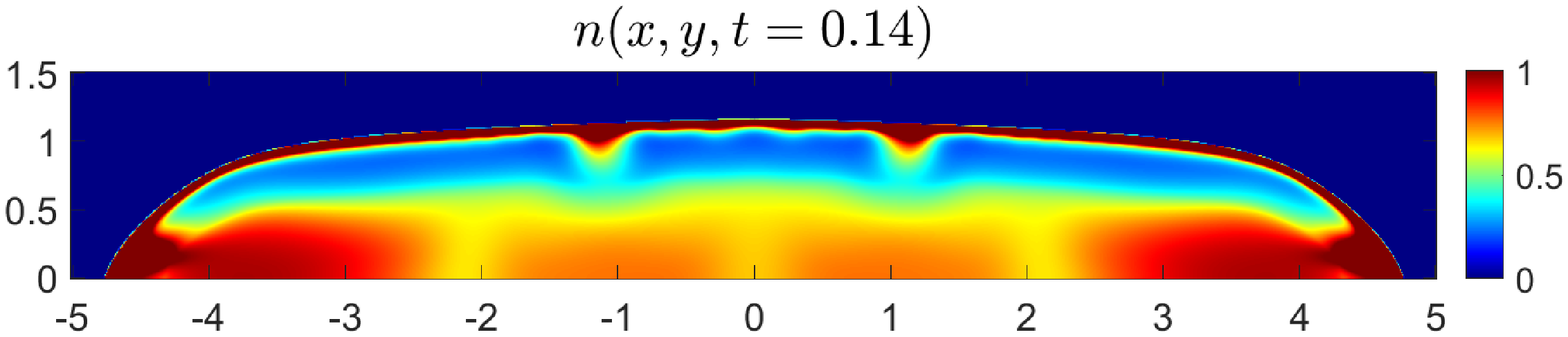}}
\vskip6pt
\centerline{\hspace*{0.1cm}\includegraphics[trim=2.4cm 0.8cm 2.3cm 0.4cm,clip,width=0.49\textwidth]{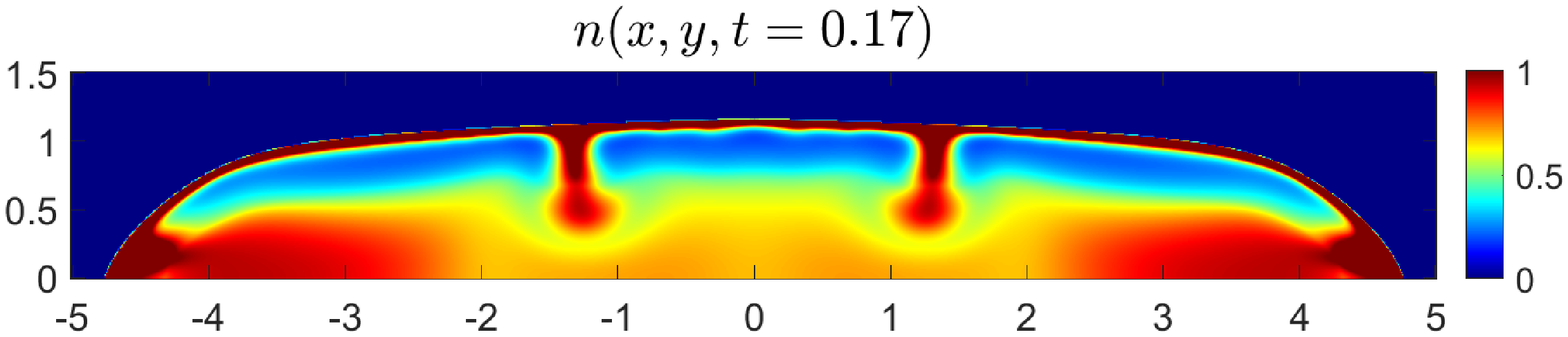}\hspace*{0.25cm}
            \includegraphics[trim=2.4cm 0.8cm 2.3cm 0.4cm,clip,width=0.49\textwidth]{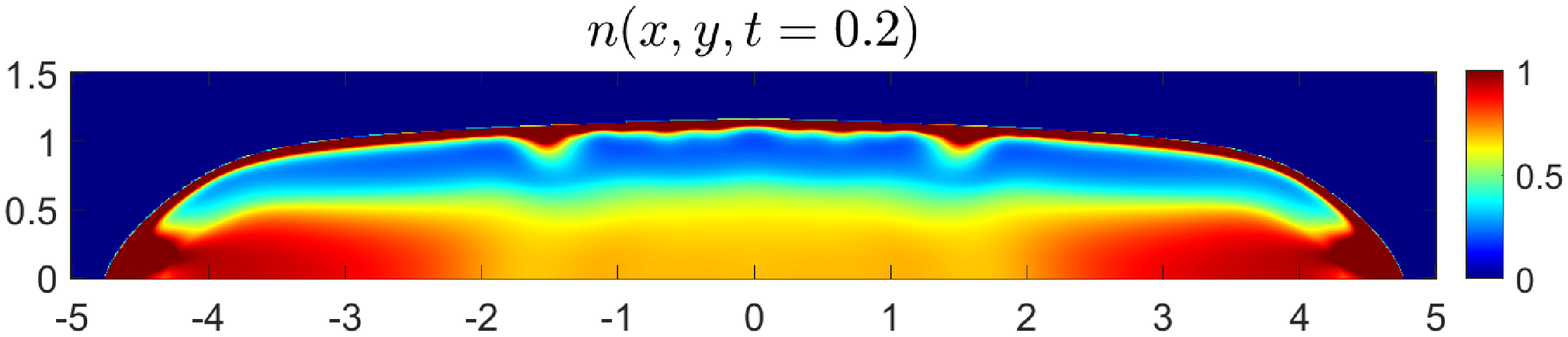}}
\vskip6pt
\centerline{\hspace*{0.1cm}\includegraphics[trim=2.4cm 0.8cm 2.3cm 0.4cm,clip,width=0.49\textwidth]{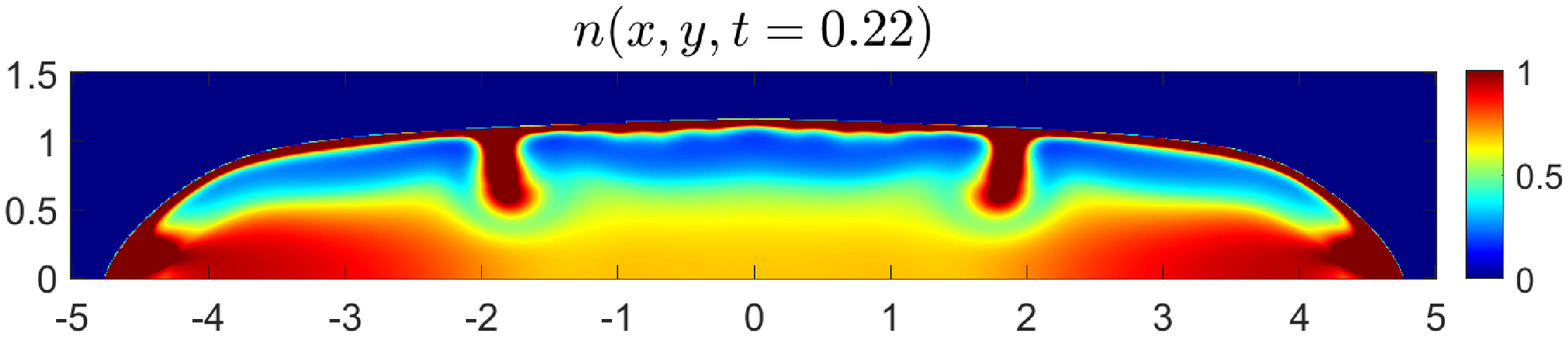}
\hspace*{0.25cm}\includegraphics[trim=2.4cm 0.8cm 2.3cm 0.4cm,clip,width=0.49\textwidth]{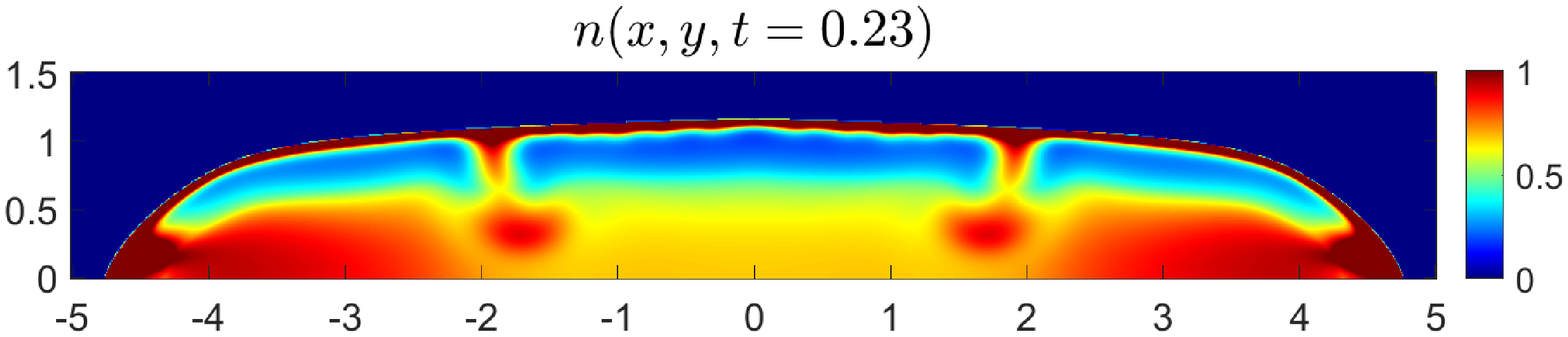}}
\vskip6pt
\centerline{\hspace*{0.1cm}\includegraphics[trim=2.4cm 0.8cm 2.3cm 0.4cm,clip,width=0.49\textwidth]{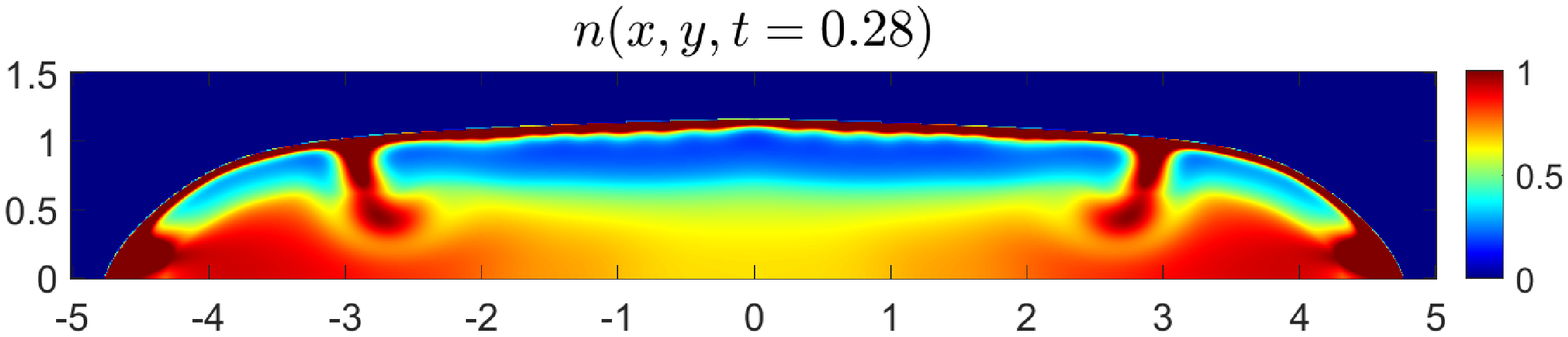}
\hspace*{0.25cm}\includegraphics[trim=2.4cm 0.8cm 2.3cm 0.4cm,clip,width=0.49\textwidth]{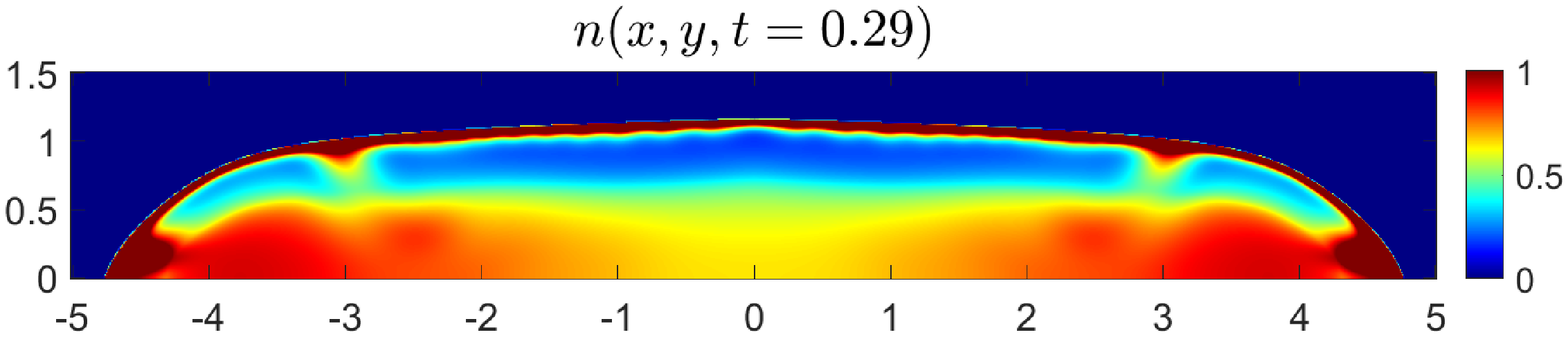}}
\vskip6pt
\centerline{\hspace*{0.1cm}\includegraphics[trim=2.4cm 0.8cm 2.3cm 0.4cm,clip,width=0.49\textwidth]{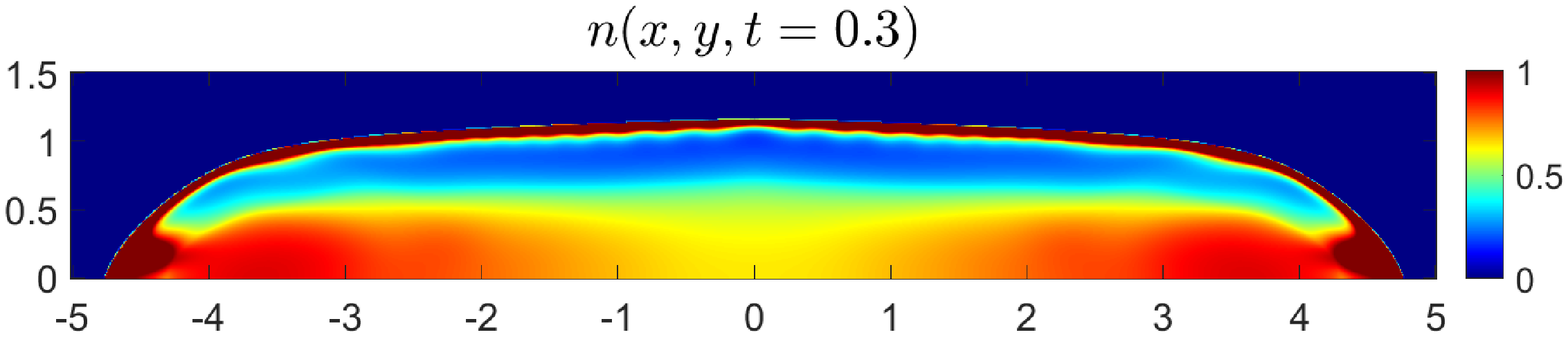}
\hspace*{0.25cm}\includegraphics[trim=2.4cm 0.8cm 2.3cm 0.4cm,clip,width=0.49\textwidth]{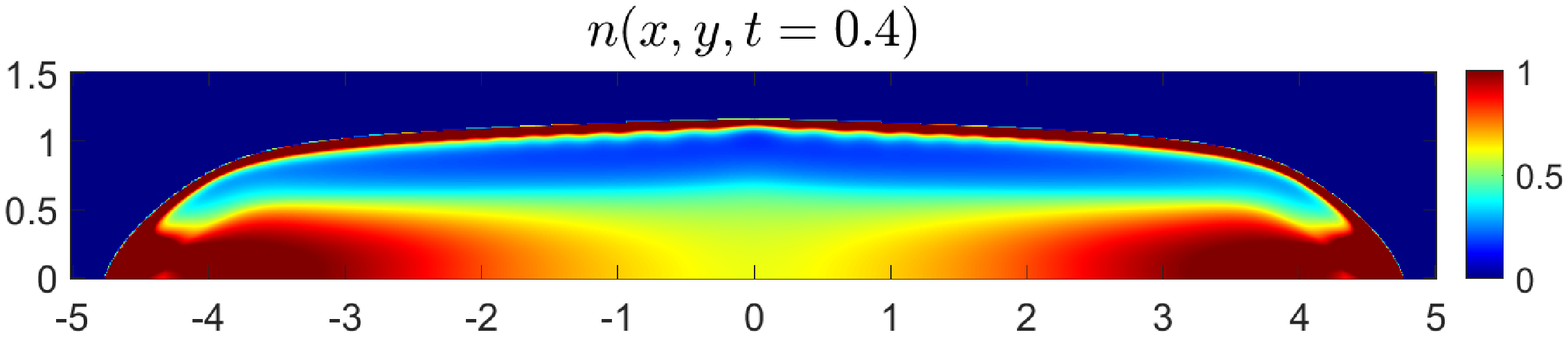}}
\vskip6pt
\centerline{\hspace*{0.1cm}\includegraphics[trim=2.4cm 0.8cm 2.3cm 0.4cm,clip,width=0.49\textwidth]{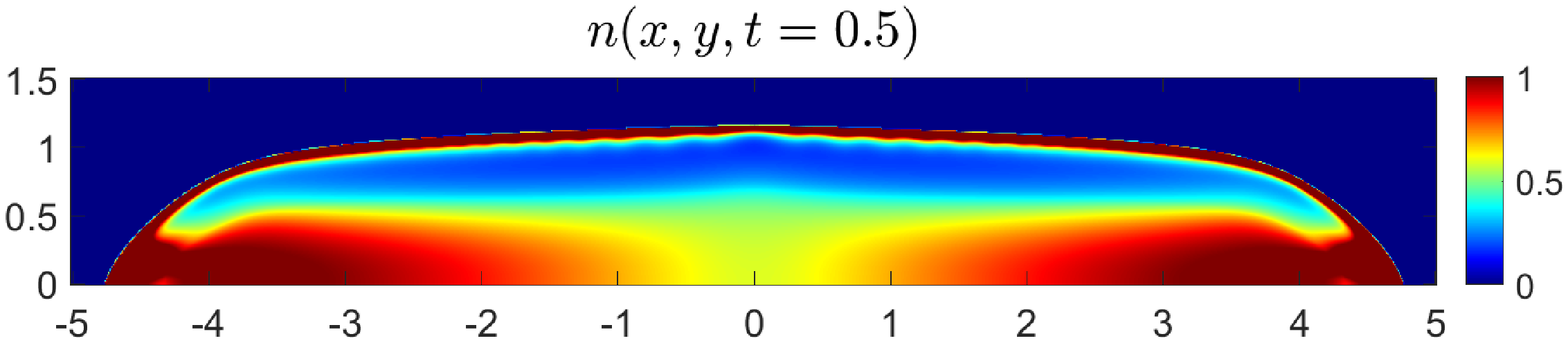}
\hspace*{0.25cm}\includegraphics[trim=2.4cm 0.8cm 2.3cm 0.4cm,clip,width=0.49\textwidth]{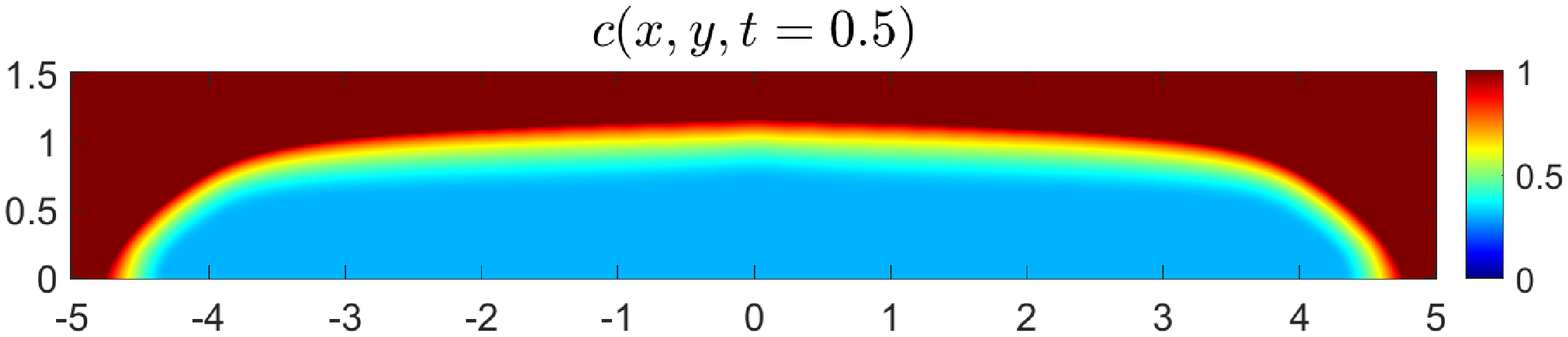}}
\caption{\sf Example 5: Time snapshots of the computed cell densities $n$ at different times and the computed oxygen concentration $c$ at
the final time.\label{fig59}}
\end{figure}

\begin{figure}[ht!]
\centerline{\includegraphics[width=0.3\textwidth]{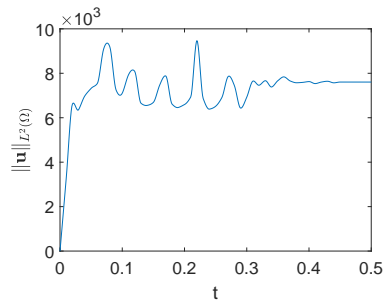}}
\caption{\sf Example 5: Time-evolution of the kinetic energy $\|\bm u\|_{L^2(\Omega)}$.\label{fig510}}
\end{figure}

\paragraph{Example 5.} In this example, we use precisely the same shape of the drop and initial data as in Example 1. The time evolution of
the computed cell density $n$ is shown in Figure \ref{fig59}. As one can see, compared with Example 1 heavier mushroom-shaped plumes are
formed by time $t=0.08$. Later on (by time $t=0.1$) these plumes are disintegrated and a part of the bacteria fall to the bottom of the drop
and become inactive due to the low oxygen concentration there. After that, smaller mushroom-shaped plumes are re-emerged and then
disintegrate several times. At the same time, one can observe the propagation of the bacteria along the top part of the drop towards its
corners. Eventually, the evolution process seems to converge to the steady state by the final time $t=0.5$ as confirmed by the stabilization
of the kinetic energy by then; see Figure \ref{fig510}. It may also be instructive to see the final time distribution of the oxygen
concentration $c$ (see the bottom right panel in Figure \ref{fig59}), which indicates that after falling down the bacteria in the lower part
of the drop remain inactive.

\begin{figure}[ht!]
\centerline{\hspace*{0.1cm}\includegraphics[trim=2.4cm 0.8cm 2.3cm 0.4cm,clip,width=0.49\textwidth]{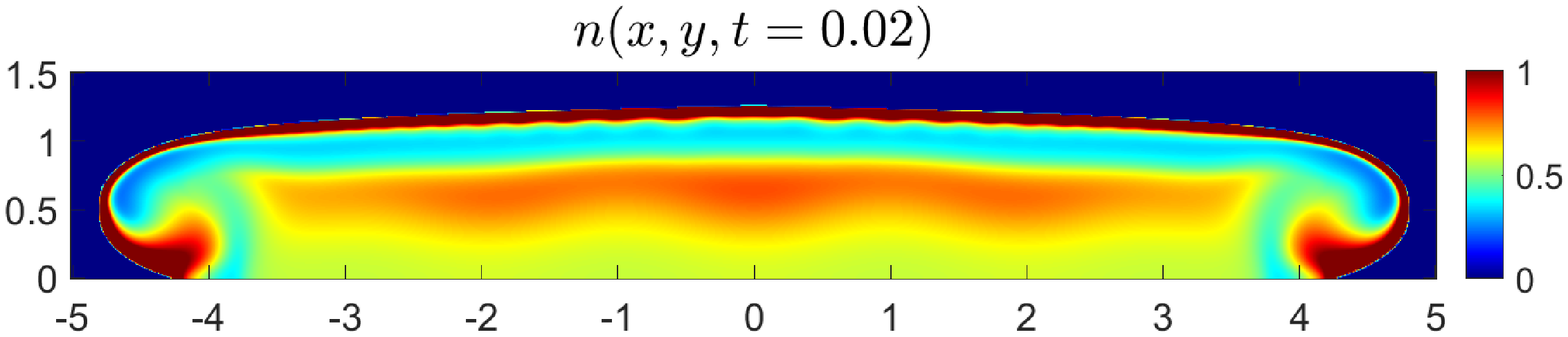}\hspace*{0.25cm}
            \includegraphics[trim=2.4cm 0.8cm 2.3cm 0.4cm,clip,width=0.49\textwidth]{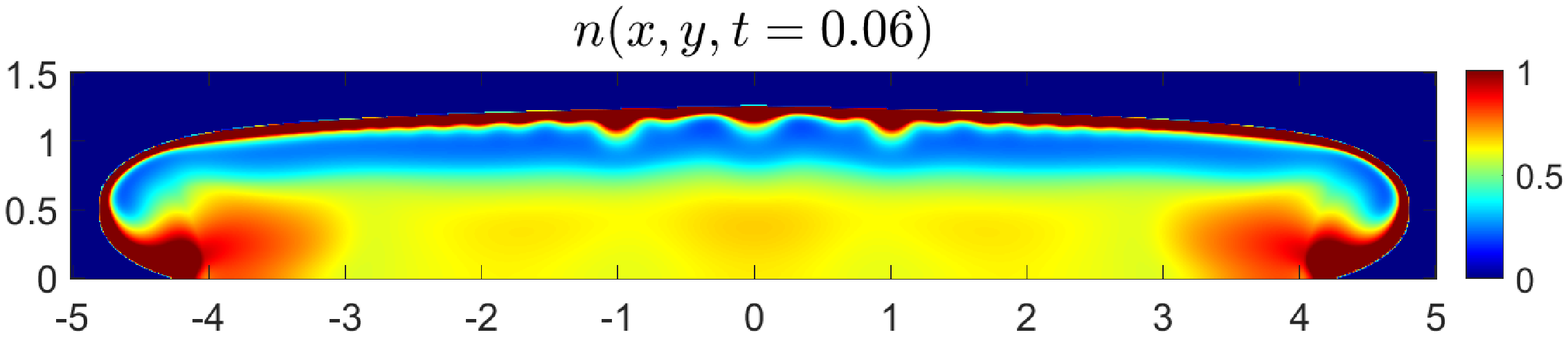}}
\vskip6pt
\centerline{\hspace*{0.1cm}\includegraphics[trim=2.4cm 0.8cm 2.3cm 0.4cm,clip,width=0.49\textwidth]{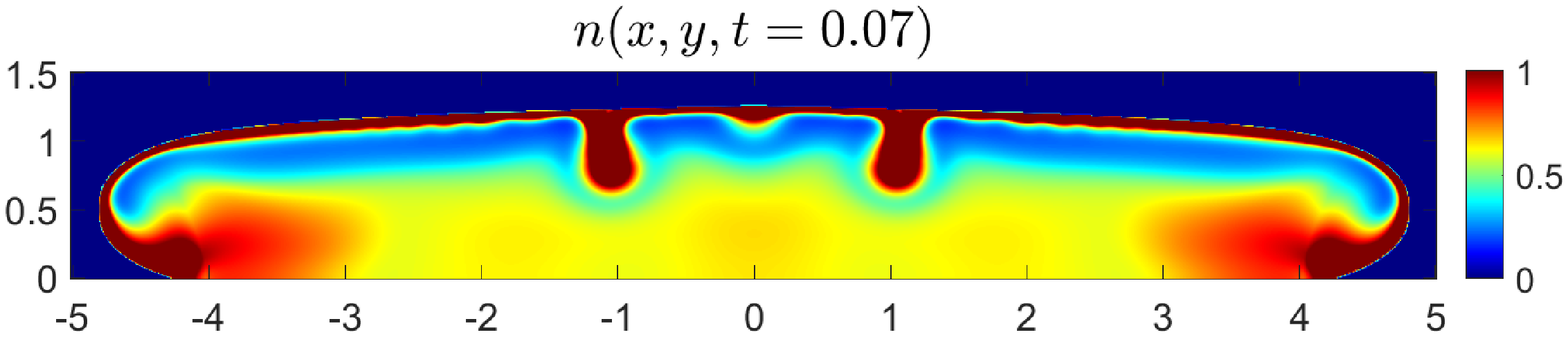}\hspace*{0.25cm}
            \includegraphics[trim=2.4cm 0.8cm 2.3cm 0.4cm,clip,width=0.49\textwidth]{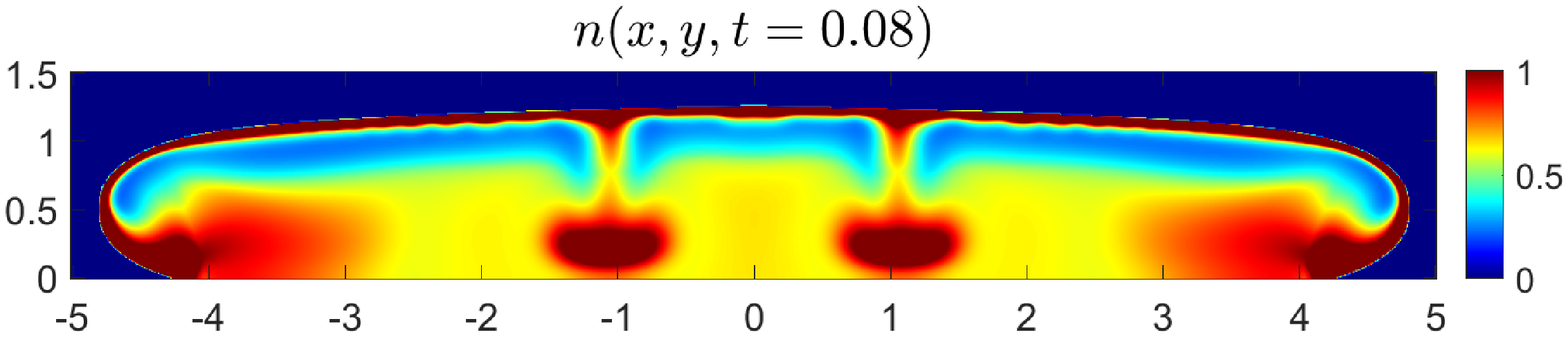}}
\vskip6pt
\centerline{\hspace*{0.1cm}\includegraphics[trim=2.4cm 0.8cm 2.3cm 0.4cm,clip,width=0.49\textwidth]{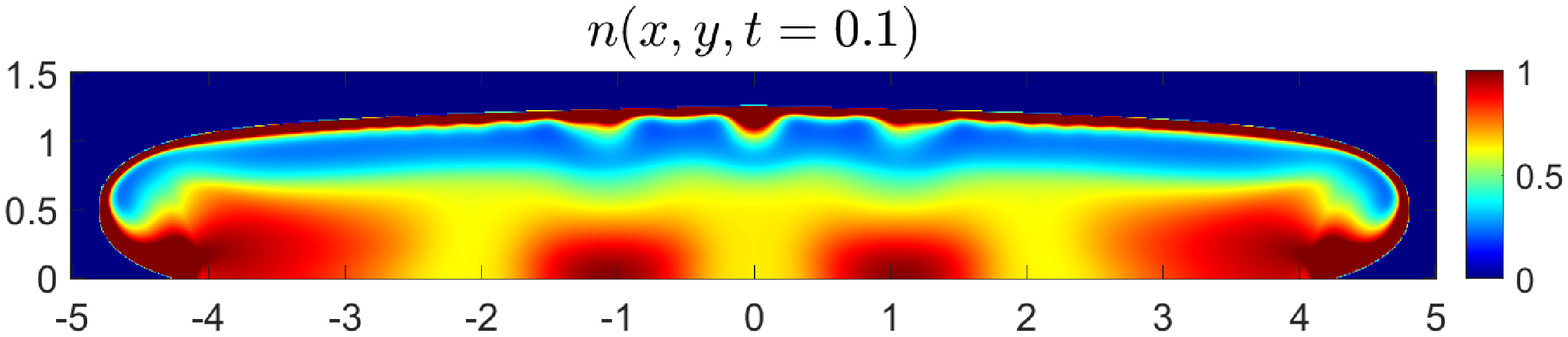}\hspace*{0.25cm}
            \includegraphics[trim=2.4cm 0.8cm 2.3cm 0.4cm,clip,width=0.49\textwidth]{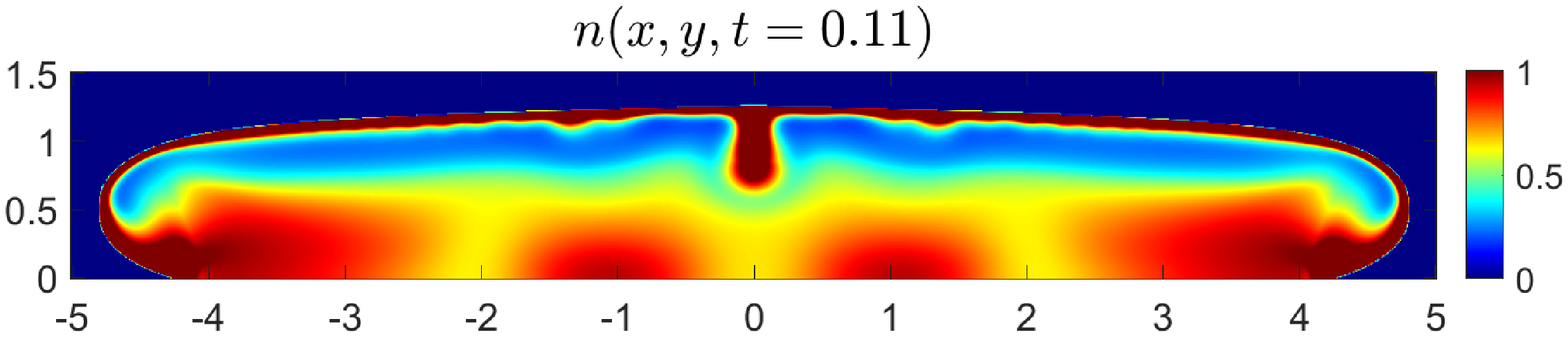}}
\vskip6pt
\centerline{\hspace*{0.1cm}\includegraphics[trim=2.4cm 0.8cm 2.3cm 0.4cm,clip,width=0.49\textwidth]{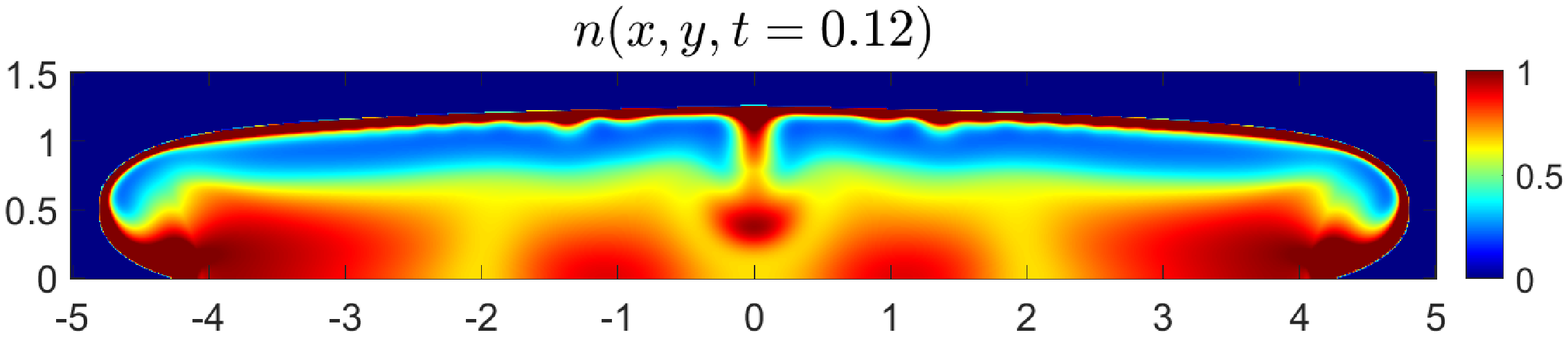}\hspace*{0.25cm}
            \includegraphics[trim=2.4cm 0.8cm 2.3cm 0.4cm,clip,width=0.49\textwidth]{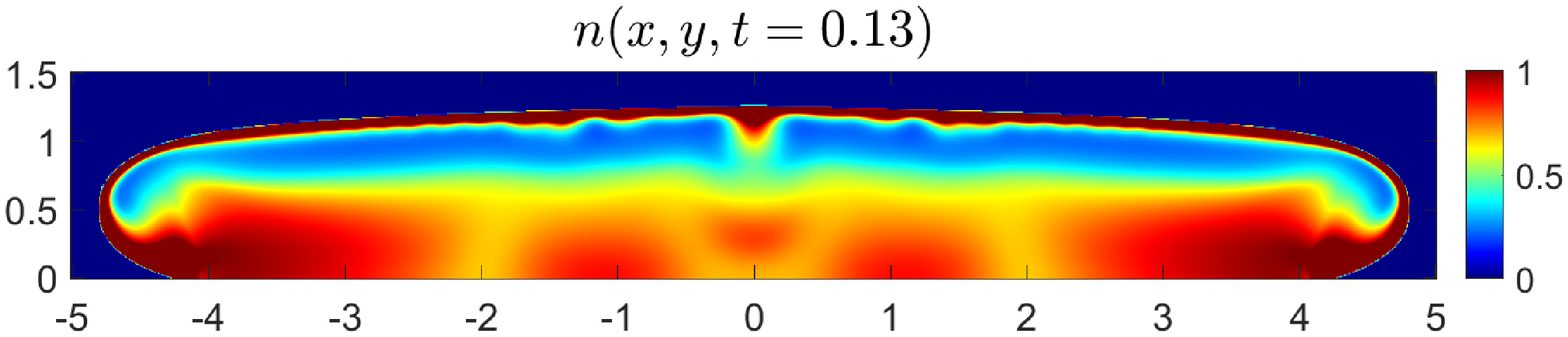}}
\vskip6pt
\centerline{\hspace*{0.1cm}\includegraphics[trim=2.4cm 0.8cm 2.3cm 0.4cm,clip,width=0.49\textwidth]{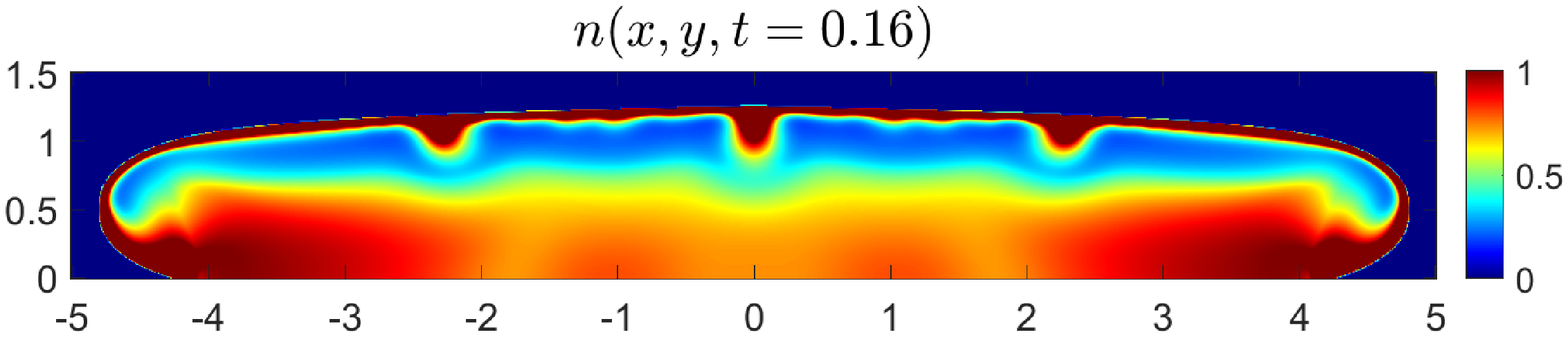}\hspace*{0.25cm}
            \includegraphics[trim=2.4cm 0.8cm 2.3cm 0.4cm,clip,width=0.49\textwidth]{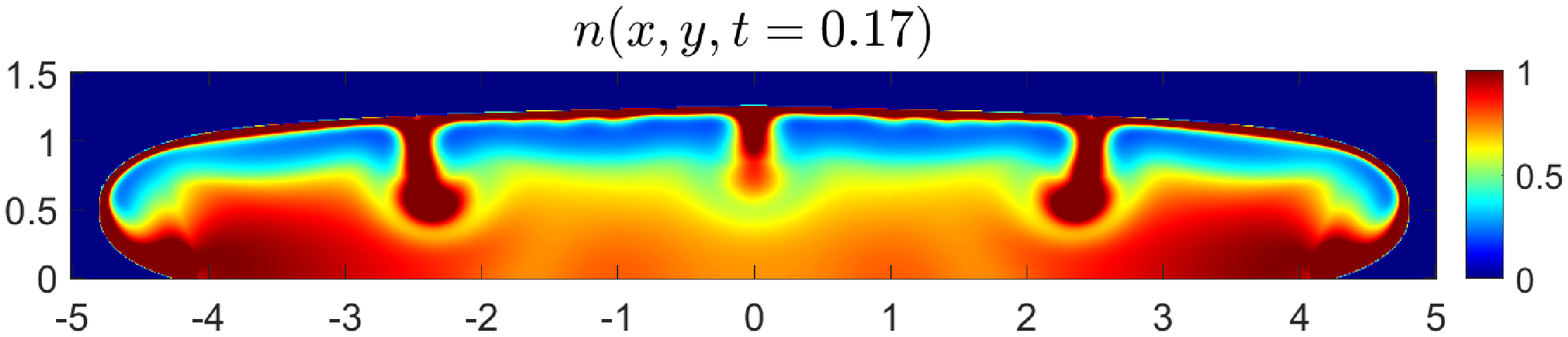}}
\vskip6pt
\centerline{\hspace*{0.1cm}\includegraphics[trim=2.4cm 0.8cm 2.3cm 0.4cm,clip,width=0.49\textwidth]{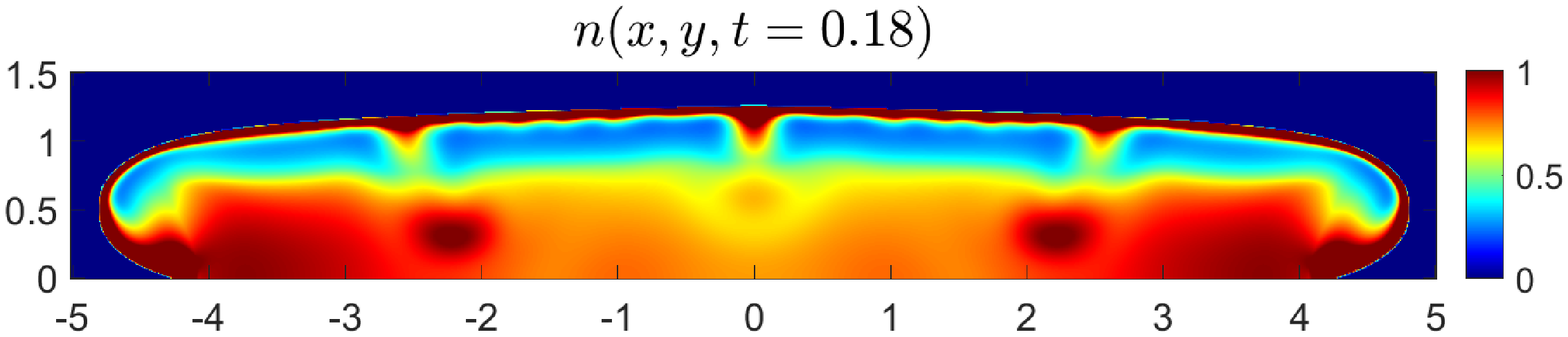}\hspace*{0.25cm}
            \includegraphics[trim=2.4cm 0.8cm 2.3cm 0.4cm,clip,width=0.49\textwidth]{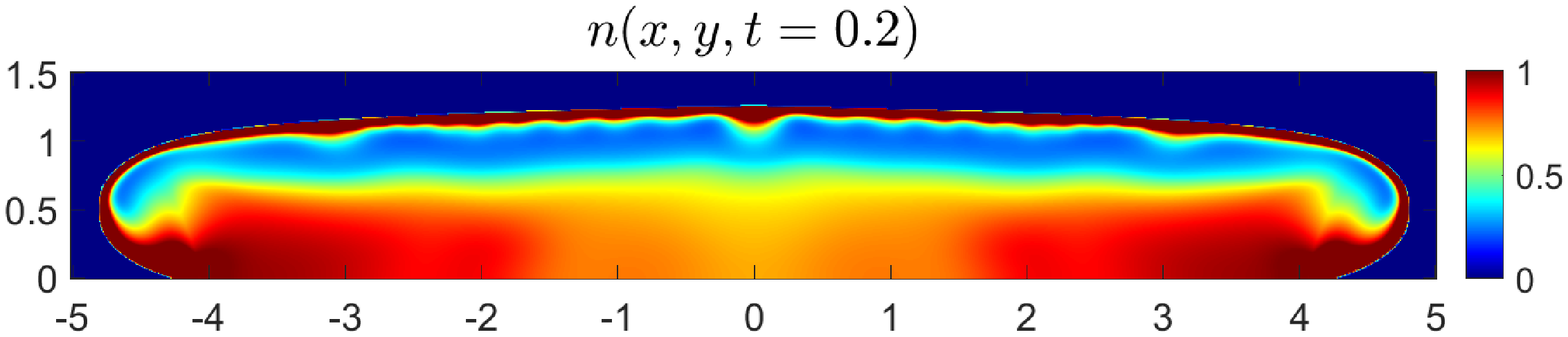}}
\vskip6pt
\centerline{\hspace*{0.1cm}\includegraphics[trim=2.4cm 0.8cm 2.3cm 0.4cm,clip,width=0.49\textwidth]{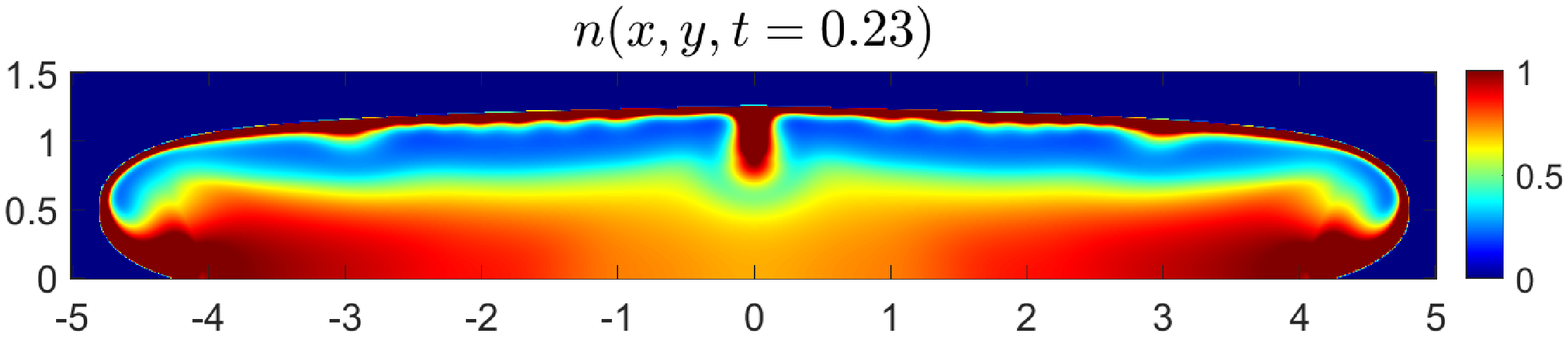}
\hspace*{0.25cm}\includegraphics[trim=2.4cm 0.8cm 2.3cm 0.4cm,clip,width=0.49\textwidth]{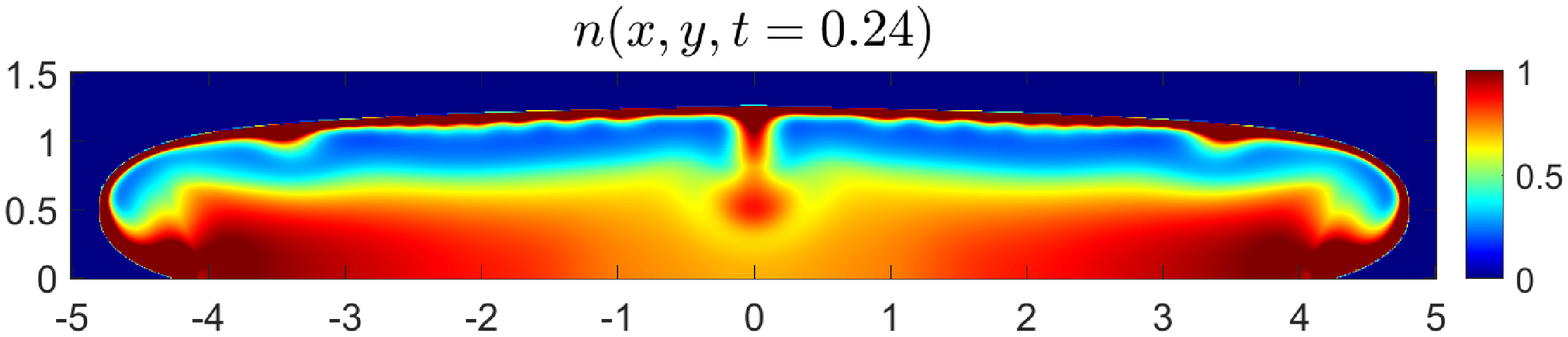}}
\vskip6pt
\centerline{\hspace*{0.1cm}\includegraphics[trim=2.4cm 0.8cm 2.3cm 0.4cm,clip,width=0.49\textwidth]{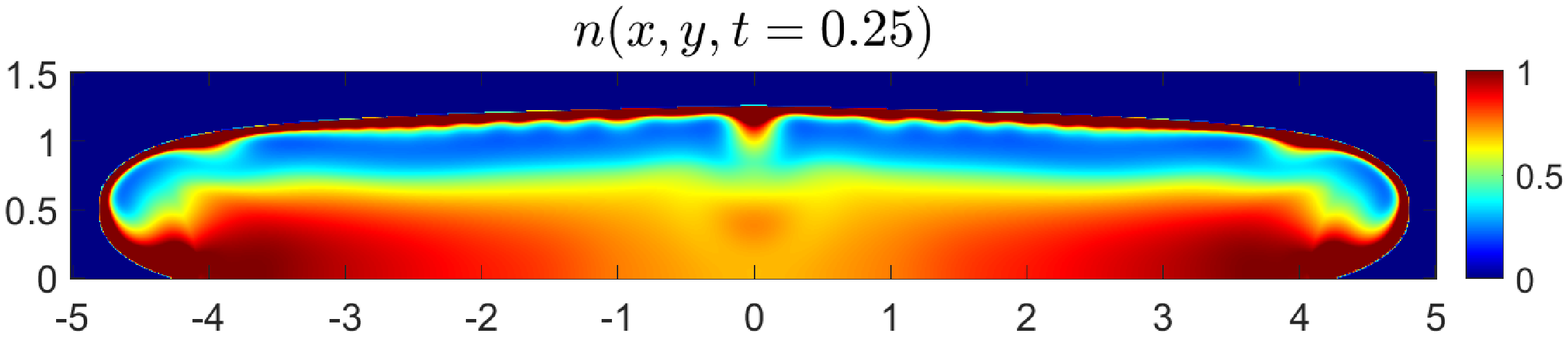}
\hspace*{0.25cm}\includegraphics[trim=2.4cm 0.8cm 2.3cm 0.4cm,clip,width=0.49\textwidth]{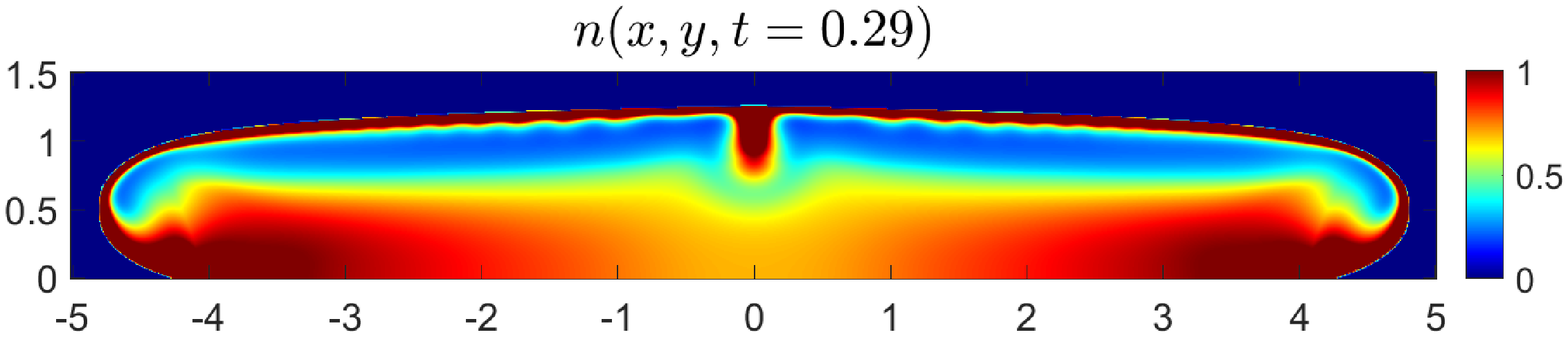}}
\vskip6pt
\centerline{\hspace*{0.1cm}\includegraphics[trim=2.4cm 0.8cm 2.3cm 0.4cm,clip,width=0.49\textwidth]{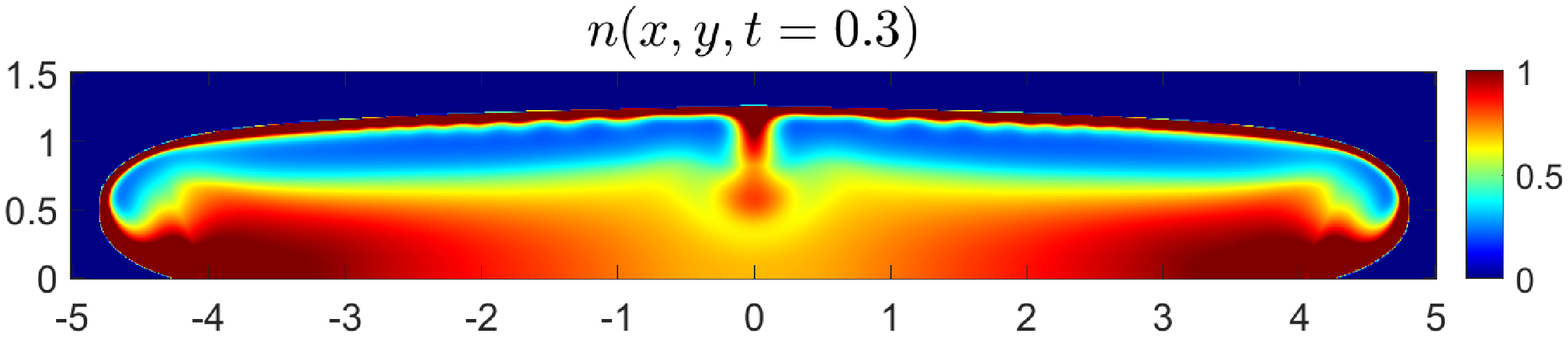}
\hspace*{0.25cm}\includegraphics[trim=2.4cm 0.8cm 2.3cm 0.4cm,clip,width=0.49\textwidth]{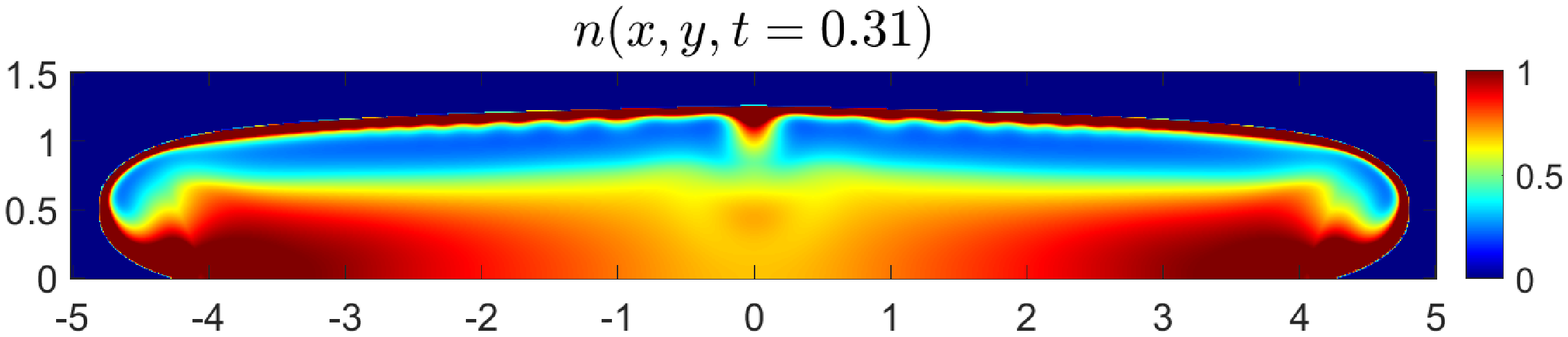}}
\vskip6pt
\centerline{\hspace*{0.1cm}\includegraphics[trim=2.4cm 0.8cm 2.3cm 0.4cm,clip,width=0.49\textwidth]{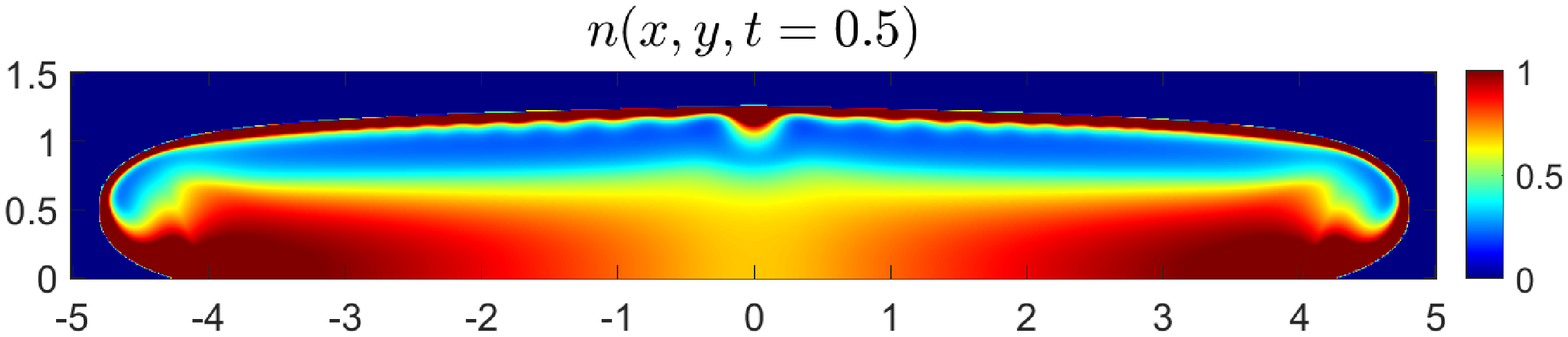}
\hspace*{0.25cm}\includegraphics[trim=2.4cm 0.8cm 2.3cm 0.4cm,clip,width=0.49\textwidth]{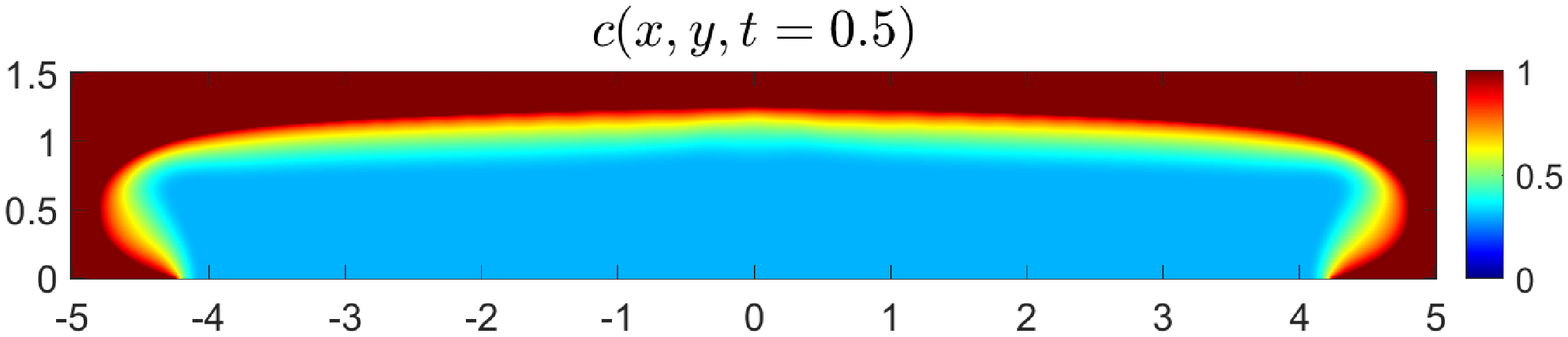}}
\caption{\sf Example 6: Time snapshots of the computed cell densities $n$ at different times and the computed oxygen concentration $c$ at
the final time.\label{fig511}}
\end{figure}

\begin{figure}[ht!]
\centerline{\includegraphics[width=0.3\textwidth]{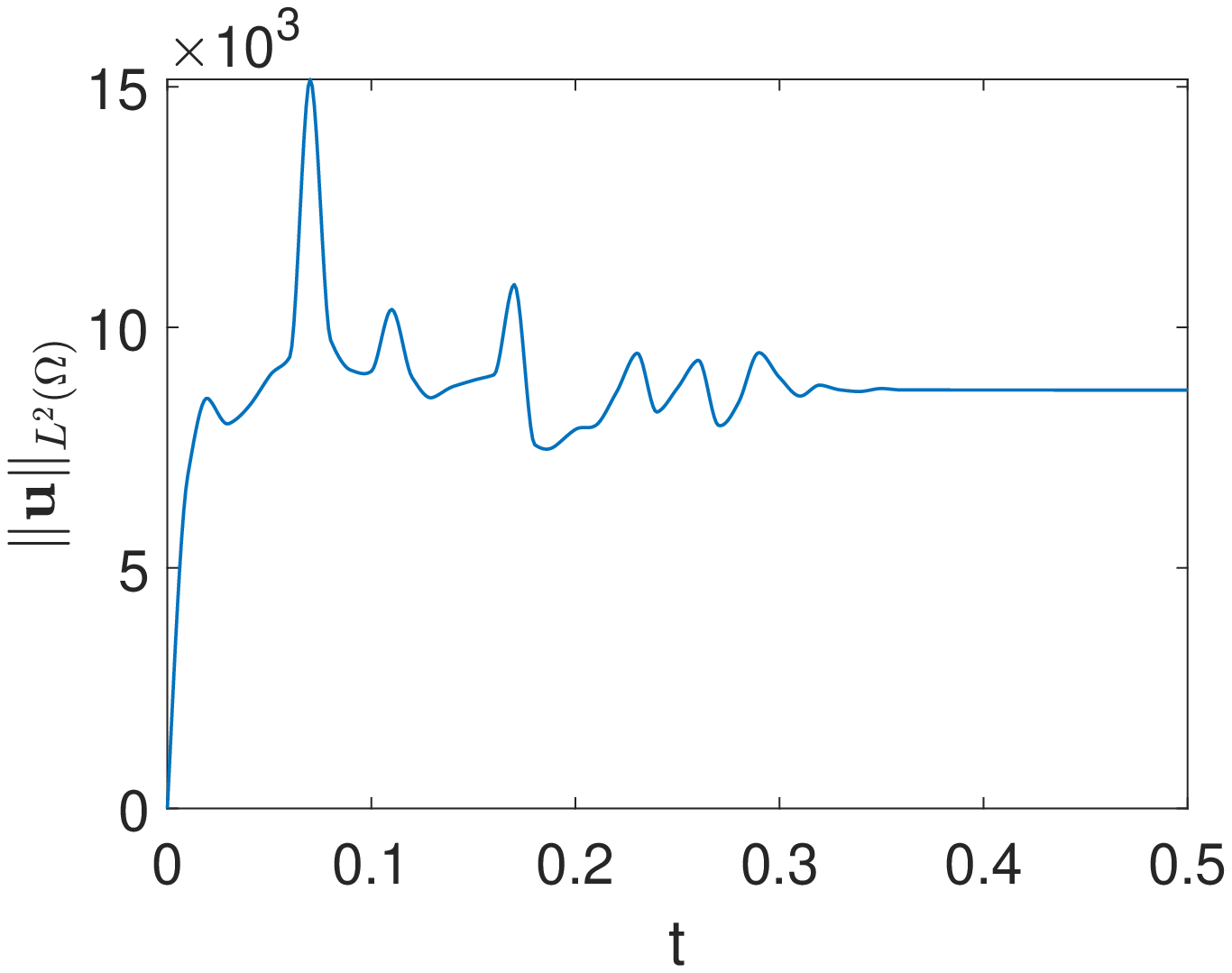}}
\caption{\sf Example 6: Time-evolution of the kinetic energy $\|\bm u\|_{L^2(\Omega)}$.\label{fig512}}
\end{figure}

\paragraph{Example 6.} In the next example, we use precisely the same shape of the drop and initial data as in Example 2. Time snapshots of
the computed cell density $n$ at different times are shown in Figure \ref{fig511}. In principle, the time evolution is quite similar to what
was observed in Example 5, but due to the difference in the shape of the drops, several distinctive features can be seen. While the
mushroom-type plums formed at about $t=0.08$ are qualitatively similar to those in Figure \ref{fig59}, the solution at later time develops a
different symmetry: by the time $t=0.16$--0.17 three plums (one in the center of the drop and two plums propagating to the sides) are
formed. Later on they keep disintegrating and re-appearing and by $t=0.24$, one can see only one plum, which remained in the center of the
drop as the other two plums practically merged with the top boundary cell layer. After that, the remaining plum keeps disintegrating and
re-emerging until the solution reaches its steady state. This convergence is confirmed by the stabilization of the kinetic energy (see
Figure \ref{fig512}) and also by the final time oxygen distribution (see the bottom right panel in Figure \ref{fig511}).

\subsection{Plumes in Sessile Drops Surrounded by Oxygen}\label{sec53}
In this section, we consider the sessile drop surrounded by oxygen. The shape of the drop is determined by the function
\begin{equation*}
f(x,y)=\begin{cases}
4.8+x-|1.5y-0.95|^{2.5}-(1.5y-0.95)^{10}&\mbox{if }x\le0,\\
4.8-x-|1.5y-0.95|^{2.5}-(1.5y-0.95)^{10}&\mbox{otherwise},
\end{cases}
\end{equation*}
representing the original domain $\Omega=\{(x,y)~|~f(x,y)>0,y>0.1\}$, for which we compute the signed distance function $d(x,y)$ to
$$
\partial\Omega=\{(x,y)~|~f(x,y)=0,y>0.1\}~\bigcup~\{(x,y)~|~f(x,y)>0,y=0.1\}
$$
needed to obtain the diffuse-domain function $\phi(x,y)$ in \eref{3.7}. In order to implement the proposed diffuse-domain based method,
$\Omega$ is imbedded into $\widetilde\Omega=[-5,5]\times[0,1.5]$.

Unlike the drops considered up to now, here we model the drop surrounded by oxygen. Therefore, the boundary conditions
\begin{equation*}
\bm\nu\!\cdot\!\bm u=0,\quad\bm\nu\!\cdot\!\nabla(\bm u\!\cdot\!\bm\tau)=0,\quad\left(\alpha n\nabla c-\nabla n\right)\!\cdot\!\bm\nu=0,
\quad c=1,\quad\forall(x,y)\in\partial\Omega,
\end{equation*}
which were used along the top portion $\Gamma$ in \eref{2.8}, are now set along the entire boundary $\partial\Omega$.

We solve the system \eref{3.1}--\eref{3.4}, \eref{3.7}, \eref{3.8} subject to the boundary conditions
$$
\bm\nu\!\cdot\!\nabla\bm u=\bm0,\quad\left(\alpha n\nabla c-\nabla n\right)\!\cdot\!\bm\nu=0,\quad c=1,\quad
\forall(x,y)\in\partial\widetilde\Omega
$$
instead of the previously used \eref{3.6}--\eref{3.5}, and the following initial data:
\begin{equation*}
\begin{aligned}
&n(x,y,0)=\begin{cases}
1&\mbox{if }y>0.599-0.01\sin\big(\pi(x-1.5)\big),\\
0.5&\mbox{otherwise},
\end{cases}\\
&c(x,y,0)\equiv1,\quad u(x,y,0)=v(x,y,0)\equiv0.
\end{aligned}
\end{equation*}

\begin{figure}[ht!]
\centerline{\includegraphics[trim=2.4cm 0.8cm 2.3cm 0.4cm,clip,width=0.49\textwidth]{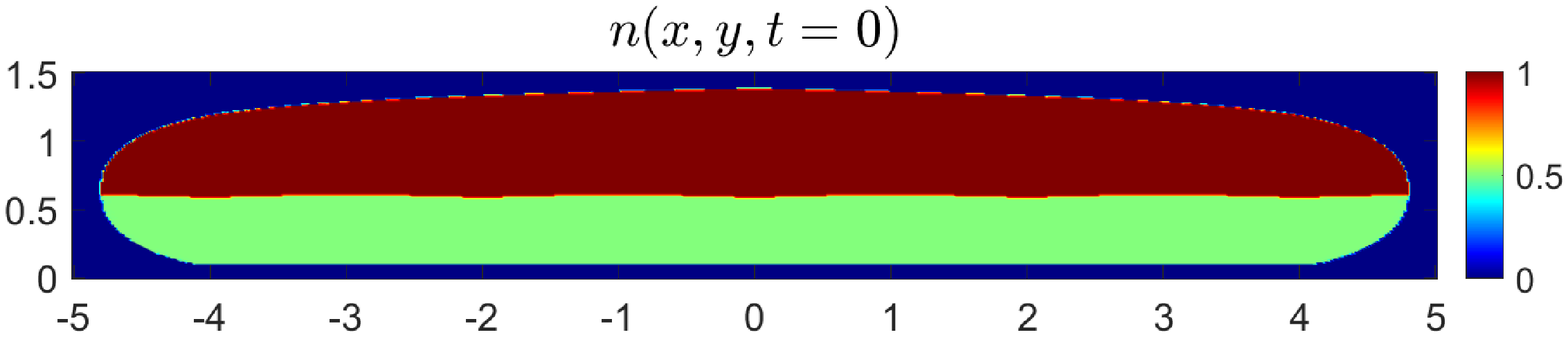}\hspace*{0.25cm}
	    \includegraphics[trim=2.4cm 0.8cm 2.3cm 0.4cm,clip,width=0.49\textwidth]{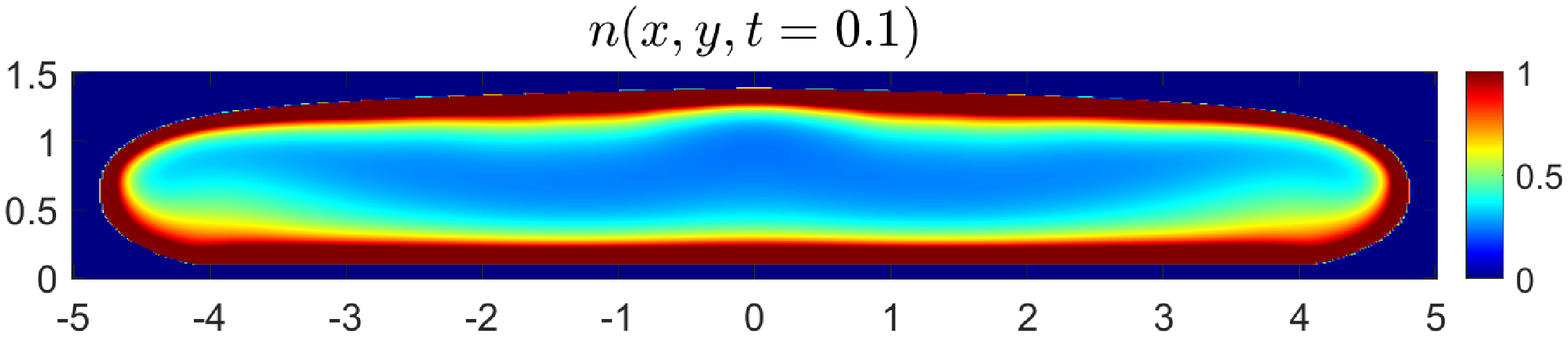}}

\centerline{\hspace*{0.1cm}\includegraphics[trim=2.4cm 0.8cm 2.3cm 0.4cm,clip,width=0.49\textwidth]{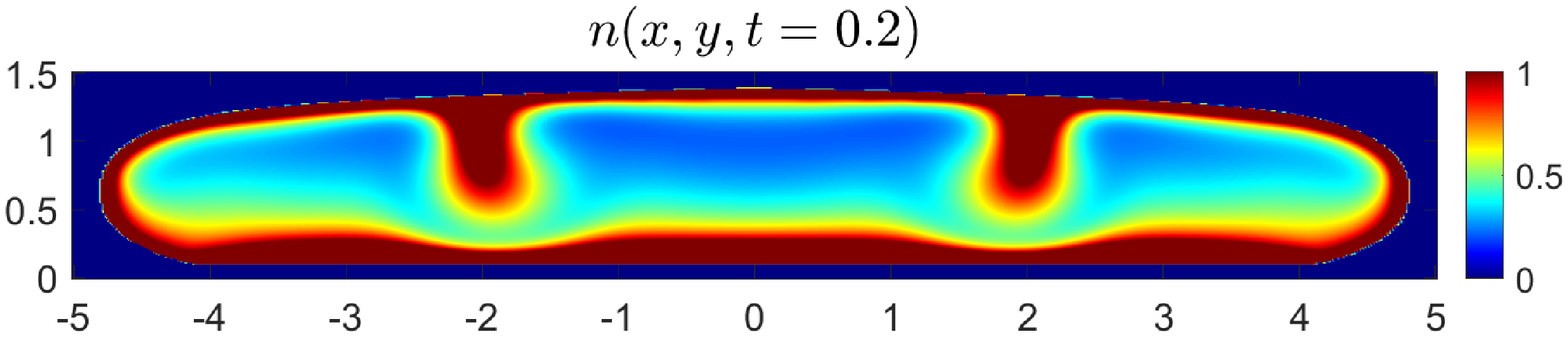}\hspace*{0.25cm}
            \includegraphics[trim=2.4cm 0.8cm 2.3cm 0.4cm,clip,width=0.49\textwidth]{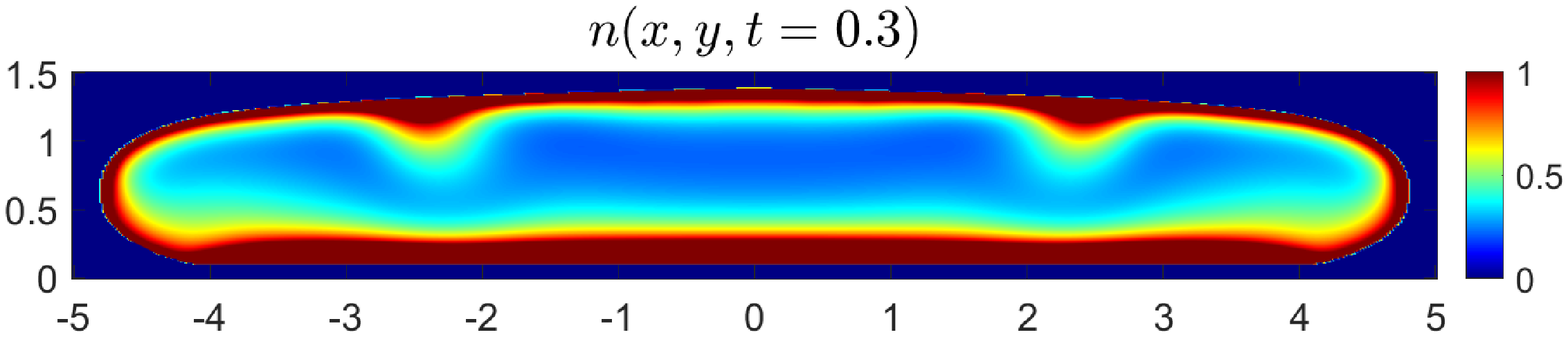}}
\vskip6pt
\centerline{\hspace*{0.1cm}\includegraphics[trim=2.4cm 0.8cm 2.3cm 0.4cm,clip,width=0.49\textwidth]{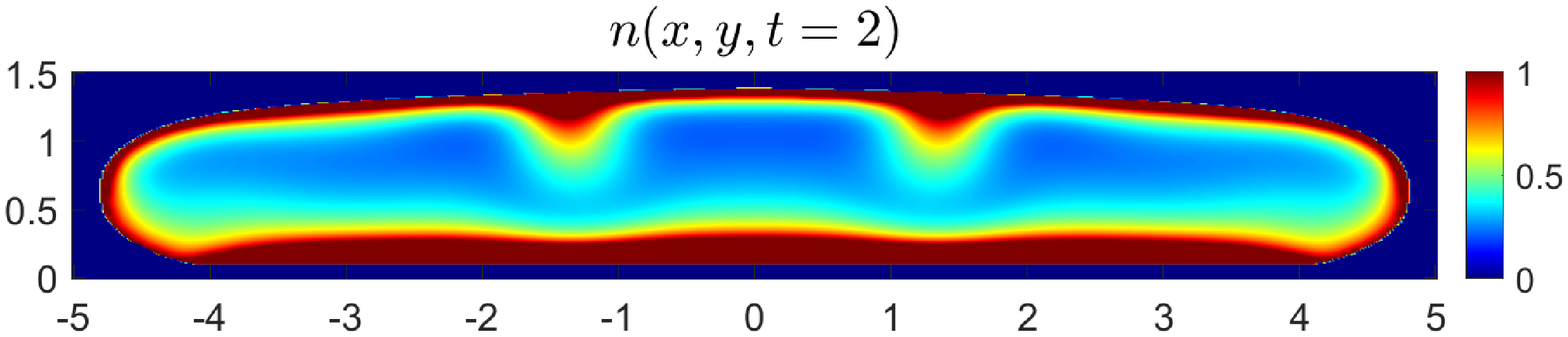}\hspace*{0.25cm}
            \includegraphics[trim=2.4cm 0.8cm 2.3cm 0.4cm,clip,width=0.49\textwidth]{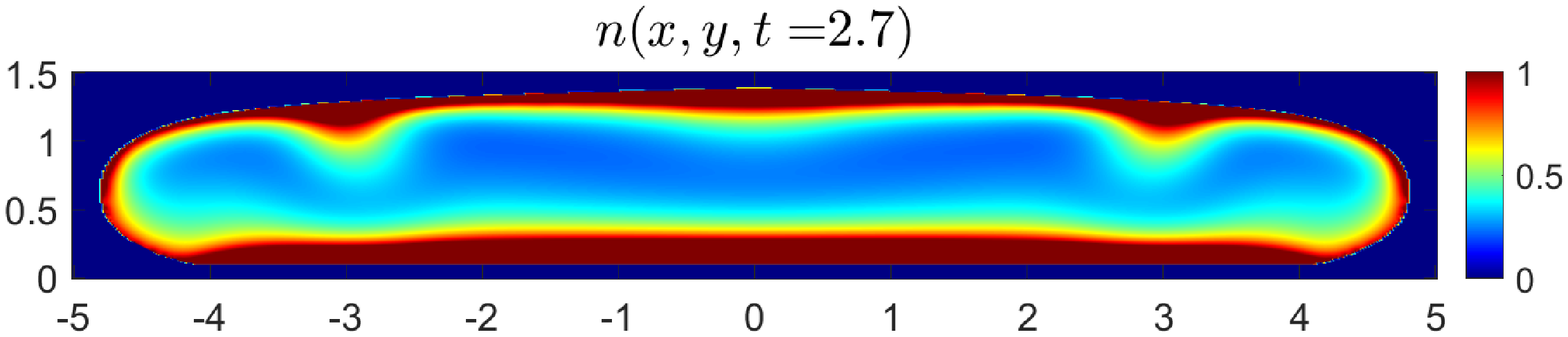}}
\vskip6pt
\centerline{\hspace*{0.1cm}\includegraphics[trim=2.4cm 0.8cm 2.3cm 0.4cm,clip,width=0.49\textwidth]{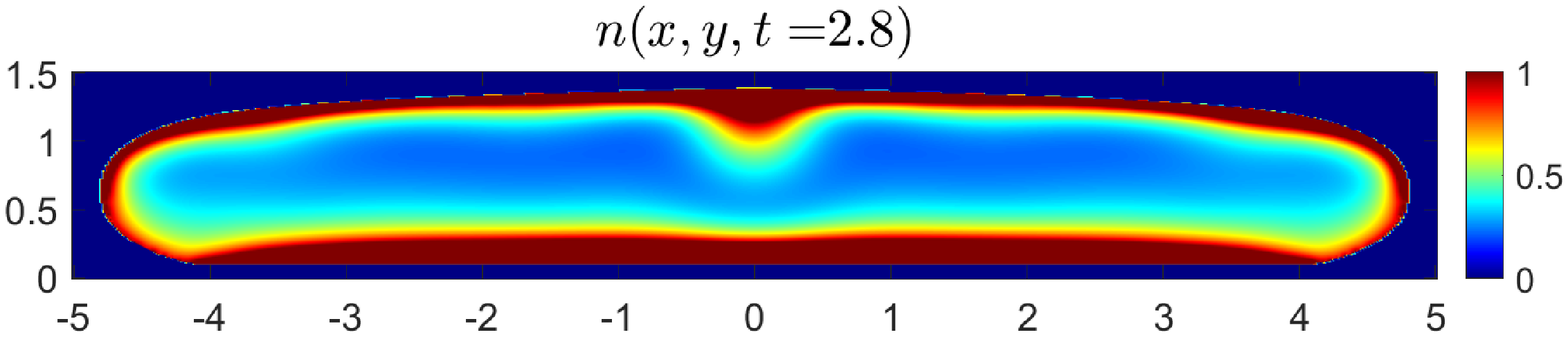}
\hspace*{0.25cm}\includegraphics[trim=2.4cm 0.8cm 2.3cm 0.4cm,clip,width=0.49\textwidth]{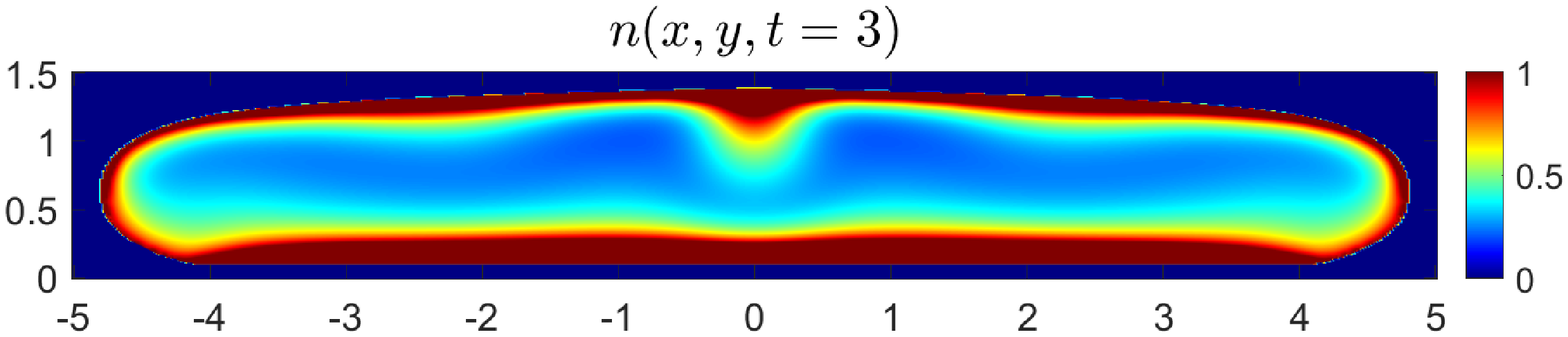}}
\vskip6pt
\centerline{\hspace*{0.1cm}\includegraphics[trim=2.4cm 0.8cm 2.3cm 0.4cm,clip,width=0.49\textwidth]{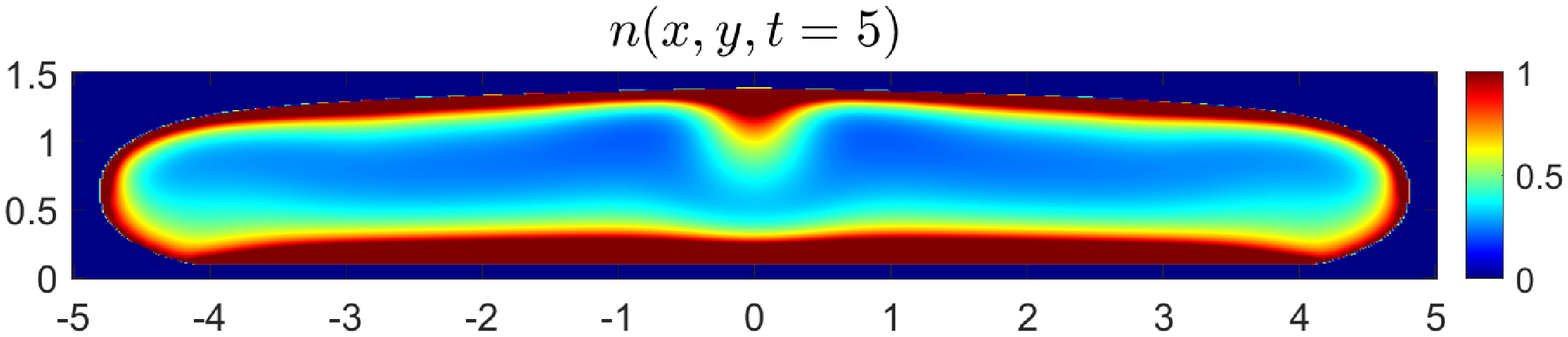}
\hspace*{0.25cm}\includegraphics[trim=2.4cm 0.8cm 2.3cm 0.4cm,clip,width=0.49\textwidth]{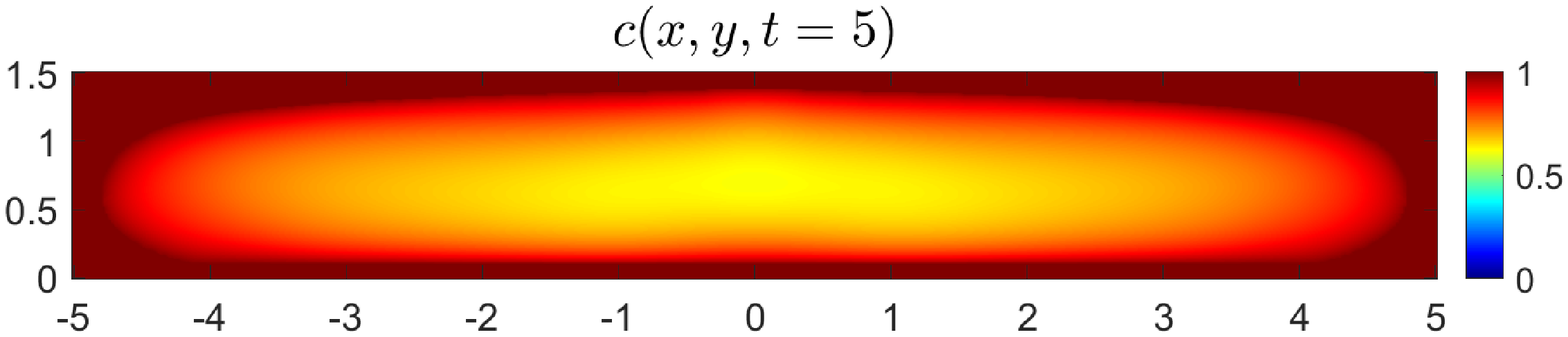}}
\caption{\sf Example 7: Time snapshots of the computed cell densities $n$ at different times and the computed oxygen concentration $c$ at
the final time.\label{fig513}}
\end{figure}

\begin{figure}[ht!]
\centerline{\includegraphics[width=0.3\textwidth]{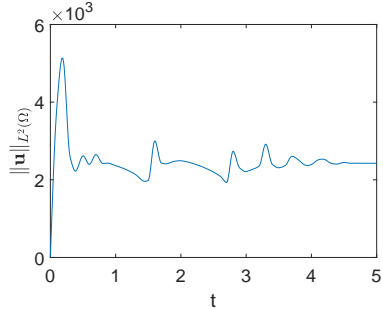}}
\caption{\sf Example 7: Time-evolution of the kinetic energy $\|\bm u\|_{L^2(\Omega)}$.\label{fig514}}
\end{figure}

\paragraph{Example 7.} In this example, we take $\beta=20$ and $\gamma=2000$ and we compute the solutions until the final time $t=5$. In
Figure \ref{fig513}, the computed cell densities at different times are plotted along with the oxygen concentration, which is shown at the
final time only. The major difference between this example and Examples 1--6 is that the oxygen is now accessible around the entire boundary
of the drop. Therefore, the bacteria immediately start propagating along the boundary towards the lower part of the drop (this can be seen
even at a small time $t=0.1$). At the same time, the gravity causes the formation of the plums (see, e.g., the solution at $t=0.2$). These
plums are unstable and later on more plums are formed. At larger times, a small plum at the center of the drop is merged and it seems to be
stable as the solution converges to its steady state by the final time; see also Figure \eref{fig514}, where the kinetic energy is depicted.

\begin{figure}[ht!]
\centerline{\hspace*{0.1cm}\includegraphics[trim=2.4cm 0.8cm 2.3cm 0.4cm,clip,width=0.49\textwidth]{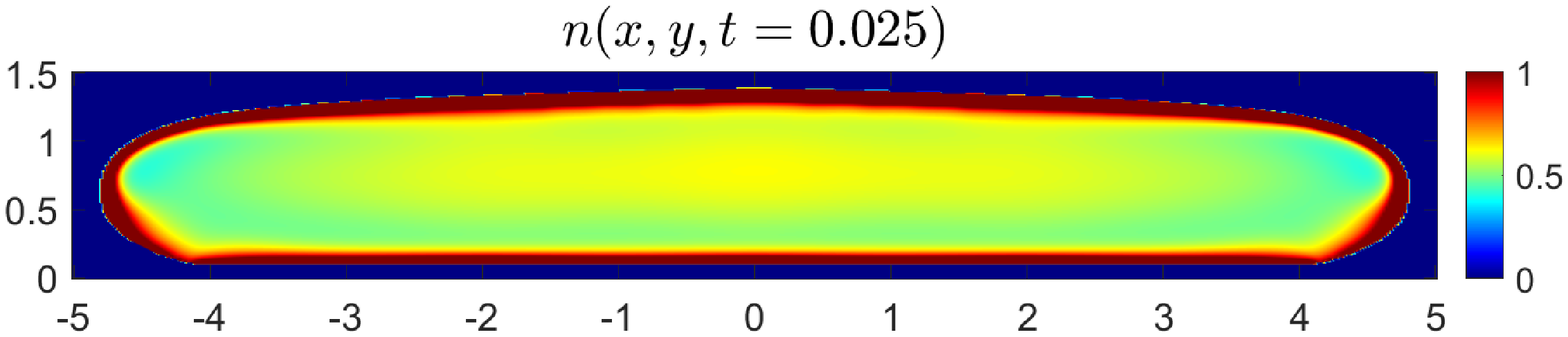}\hspace*{0.25cm}
            \includegraphics[trim=2.4cm 0.8cm 2.3cm 0.4cm,clip,width=0.49\textwidth]{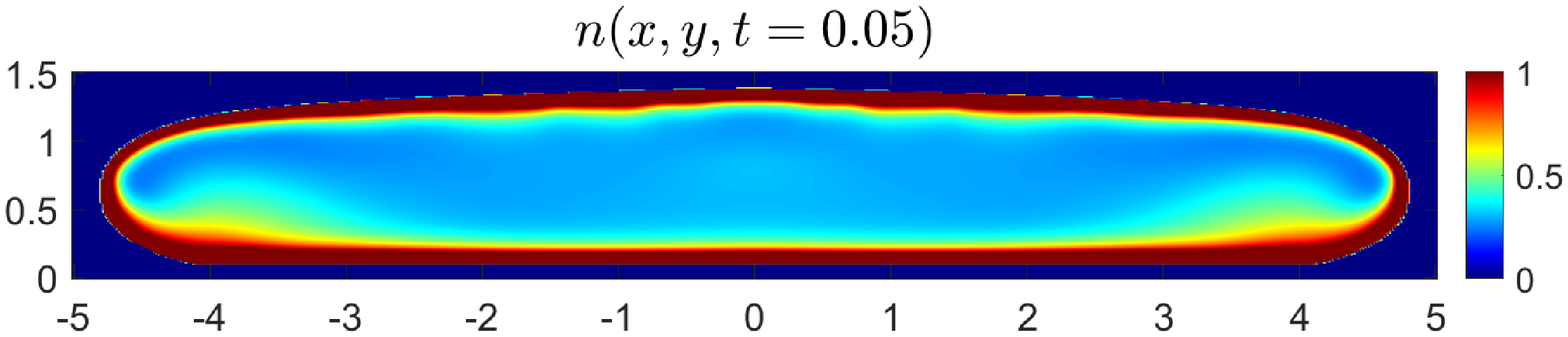}}
\vskip6pt
\centerline{\hspace*{0.1cm}\includegraphics[trim=2.4cm 0.8cm 2.3cm 0.4cm,clip,width=0.49\textwidth]{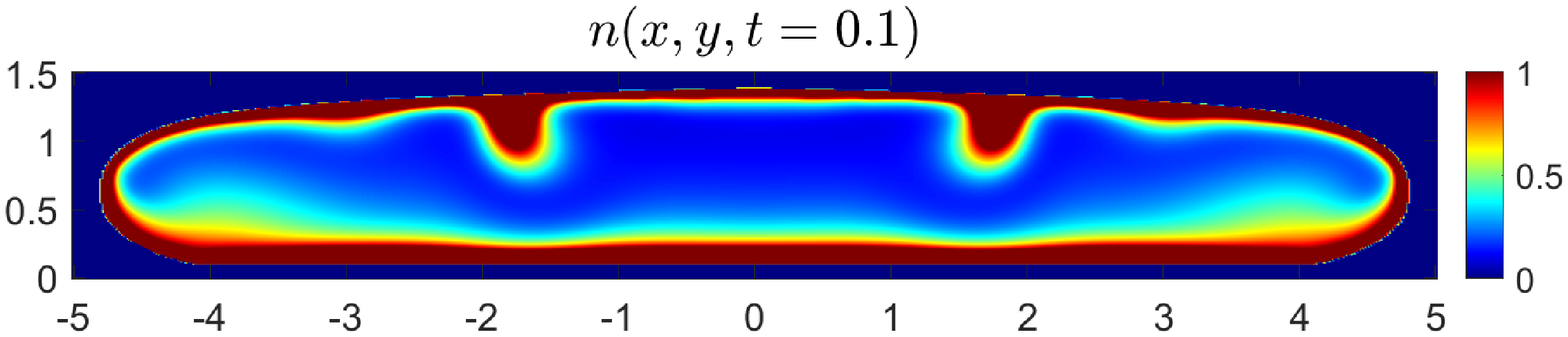}\hspace*{0.25cm}
            \includegraphics[trim=2.4cm 0.8cm 2.3cm 0.4cm,clip,width=0.49\textwidth]{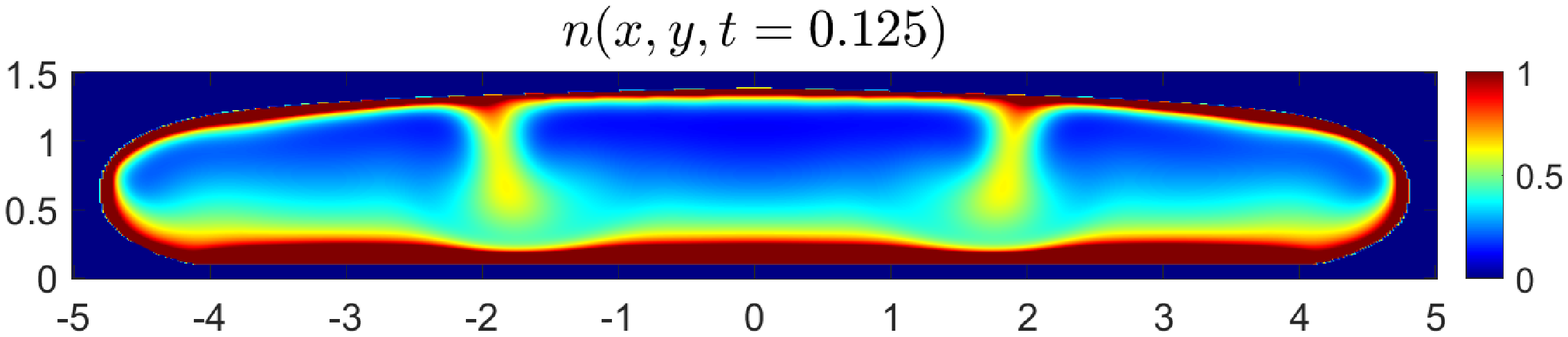}}
\vskip6pt
\centerline{\hspace*{0.1cm}\includegraphics[trim=2.4cm 0.8cm 2.3cm 0.4cm,clip,width=0.49\textwidth]{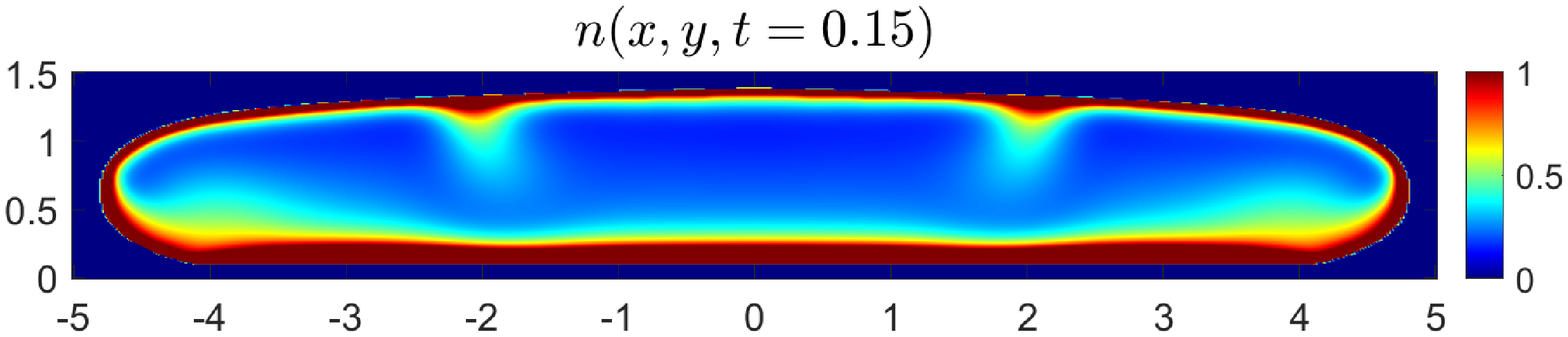}\hspace*{0.25cm}
            \includegraphics[trim=2.4cm 0.8cm 2.3cm 0.4cm,clip,width=0.49\textwidth]{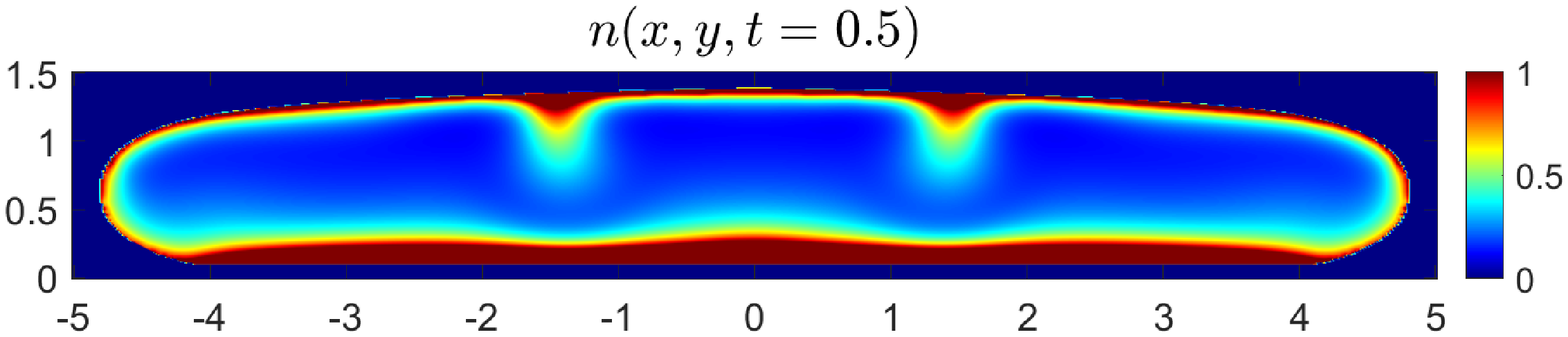}}
\vskip6pt
\centerline{\hspace*{0.1cm}\includegraphics[trim=2.4cm 0.8cm 2.3cm 0.4cm,clip,width=0.49\textwidth]{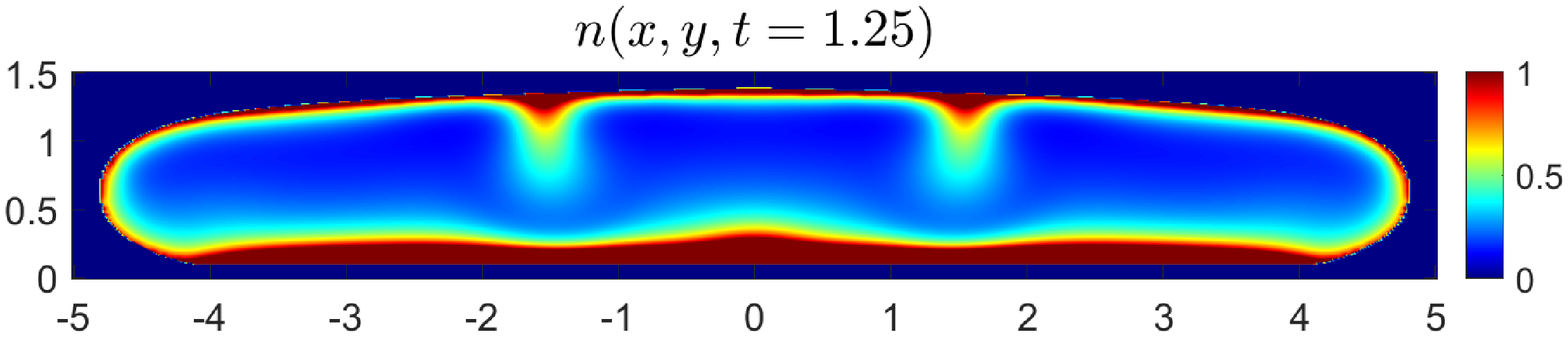}
\hspace*{0.25cm}\includegraphics[trim=2.4cm 0.8cm 2.3cm 0.4cm,clip,width=0.49\textwidth]{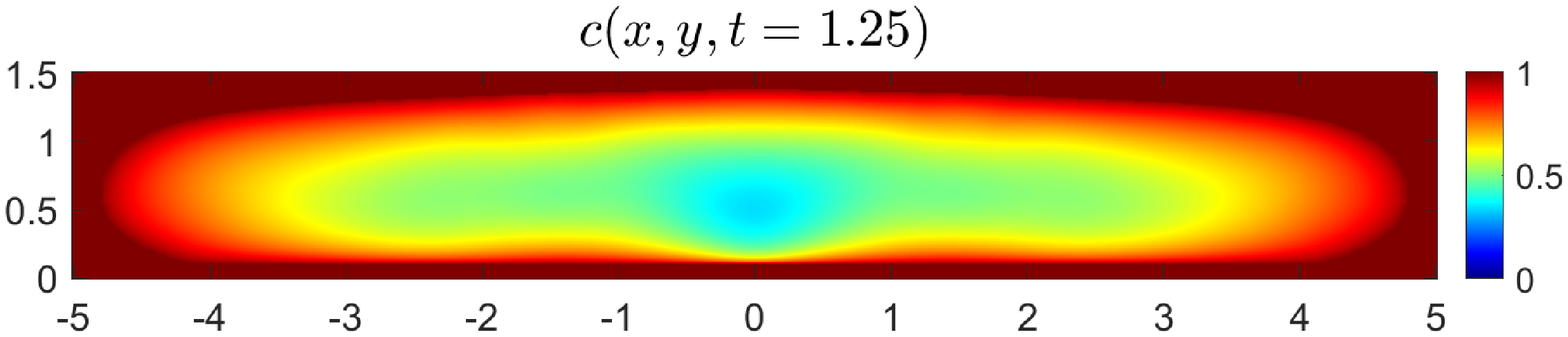}}
\caption{\sf Example 8: Time snapshots of the computed cell densities $n$ at different times and the computed oxygen concentration $c$ at
the final time.\label{fig516}}
\end{figure}

\begin{figure}[ht!]
\centerline{\includegraphics[width=0.3\textwidth]{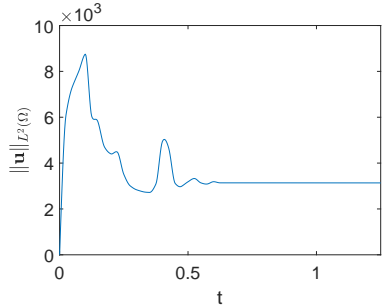}}
\caption{\sf Example 8: Time-evolution of the kinetic energy $\|\bm u\|_{L^2(\Omega)}$.\label{fig517}}
\end{figure}

\paragraph{Example 8.} The final example is similar to Example 7 with the only exception that here we take $\beta=40$ and $\gamma=4000$.
These values correspond to a twice larger reference cell density $n_r$, which leads to a faster dynamics. Indeed, as one can clearly see
from Figures \ref{fig516} and \ref{fig517}, the solution converges to its steady state substantially faster than in Example 7. It should
also be observed that the obtained steady state qualitatively different from the one reported in Figure \ref{fig513}: the steady state now
contains two plums (not only one plum as in the previous example) and there is a slightly larger concentration of bacteria in the internal
part of the drop (this can be clearly seen from the final time oxygen distribution shown in the lower right panel of Figure \ref{fig516}).

\section{Conclusion}\label{sec6}
In this paper, we have introduced a new positivity preserving and high-resolution method for the coupled chemotaxis-fluid system in a
sessile drop. Our method is based on the diffuse-domain approach, which is implemented to derive a chemotaxis-fluid diffuse-domain (cf-DD)
model. We have shown that the obtained cf-DD system converges to the original chemotaxis-fluid system as the thickness of the diffuse-domain
interface shrinks to zero. In order to numerically solve the resulting cf-DD system, we have developed a second-order hybrid finite-volume
finite-difference method, which preserves non-negativity of the computed cell density.

We have tested the proposed diffuse-domain based on a number of numerical experiments, in which we have not only demonstrated the ability of
our method to handle complex computational domains, but also systematically studied bacteria collective behavior in sessile droplets of a
variety of different shapes. It has been observed that when the amount of bacteria is moderate, stable stationary plumes are formed inside
the droplet and substantial amount of the bacteria will aggregate in the corners of the droplet while creating vortices there. When the
amount of bacteria is increased, the finger-like plumes flare out into mushroom-shaped plumes, which are, however, unstable and disintegrate
in time. At the same time, the solutions converge to nontrivial steady states in all of the studied examples. Based on the obtained
numerical results, we conjecture that the evolution of bacteria is related to both the total amount of bacteria in the droplet and the shape
of the droplet. The main goal of the presented simulations is to demonstrate that the proposed numerical method can provide one with a
valuable insight on the bacteria collective behavior in complex geometries, whose detailed study is left for the future work.

\section*{Acknowledgments}
The work of A. Chertock was supported in part by NSF grants DMS-1818684 and DMS-2208438. The work of A. Kurganov was supported in part by
NSFC grants 12111530004 and 12171226, and by the fund of the Guangdong Provincial Key Laboratory of Computational Science and Material
Design (No. 2019B030301001). The work of Zhen Zhang supported in part by the NSFC grants 11731006 and 12071207, and the Natural Science
Foundation of Guangdong Province (2021A1515010359).

%
%

\bibliography{refs}
\bibliographystyle{siam}

\end{document}